\documentclass[12pt,a4paper,leqno]{article}
\usepackage{amscd,amsmath,amssymb,amsfonts}
\input xy
\xyoption{all}

\usepackage{epsf}
\newcommand{\myfig}[2]{
\begin{center}
{\epsfxsize=#2\hsize \epsfbox{#1}}\nobreak
\end{center}}

\usepackage{epic,eepic}
\usepackage{graphicx}

\def\noi{\noindent}

\def\Spec{\mathrm{Spec}}
\def\dim{\mathrm{dim}}
\def\qz{{\mathbb Q/\mathbb Z}}
\def\qzl{\mathbb Q_\ell/\mathbb Z_\ell}
\def\funct{\raisebox{1pt}{$\shortmid$}\hspace{-5pt}\rightsquigarrow}
\def\zn{{\mathbb Z/n\mathbb Z}}
\def\zln{{\z/\ell^n\z}}
\def\zlm{{\z/\ell^m\z}}
\def\z{{\mathbb Z}}
\def\q{{\mathbb Q}}
\def\et{{\acute{e}t}}
\def\qzl{\mathbb Q_\ell/\mathbb Z_\ell}
\def\ker{\mathrm{ker}}
\def\coker{\mathrm{coker}}
\def\et{{\mbox{\scriptsize \'{e}t}}}
\def\Hom{\mathrm{Hom}}
\def\Cb{\overline{C}}
\def\Hb{\overline{H}}

\def\cC{{\mathcal C}}

\def\cP{{\mathcal P}}
\def\cO{{\mathcal O}}
\def\cSP{{\mathcal{SP}}}
\def\cM{{\mathcal{M}}}
\def\cV{{\mathcal V}}

\def\ch{\mathrm{char}}

\def\oH{{\overline{H}}}
\def\ocH{{\overline{\mathcal H}}}
\def\ocC{{\overline{\mathcal C}}}
\def\WH{{H^W}}

\def\v{\vspace{0,5cm}}

\begin{document}

\vspace{20cm}

\centerline{\large \textbf{Hasse principles for higher-dimensional fields}}

\vspace{1,0cm}

\centerline{\textbf{by Uwe Jannsen}}

%\vspace{0,5cm}

%\centerline{October 26, 2014}

\vspace{1,5cm}

\noi
dedicated to J\"urgen Neukirch

\vspace{0,5cm}

\begin{abstract}
For rather general excellent schemes $X$, K. Kato defined complexes of Gersten-Bloch-Ogus
type involving the Galois cohomology groups of all residue fields of $X$. For arithmetically
interesting schemes, he developed a fascinating web of conjectures on some of these complexes,
which generalize the classical Hasse principle for Brauer groups over global fields,
and proved these conjectures for low dimensions. We prove Kato's conjecture over global fields
in any dimension, for coefficients prime to the characteristic. In particular this gives a
cohomological Hasse principle for function fields $F$ over a global field $K$, involving the
corresponding function fields $F_v$ over the completions $K_v$ of $K$. We also prove
a similar conjecture of Kato over finite fields, again with coefficients prime to the characteristic,
and a generalization to arbitrary finitely generated fields.
We get conditional results for the $p$-torsion cases, where $p$ is the characteristic, assuming
some conditions on resolution of singularities.
\end{abstract}

%\vspace{0,4cm}
\noi\small{
{2010 Mathematics Subject Classification$:$ 14F20, 14G25, 11R34}\\
{Keywords$:$ \'etale cohomology, varieties over global fields, function fields, Galois cohomology}
}
\vspace{1cm}

\centerline{\noindent\textbf{\S 0 Introduction}}

\vspace{1,0cm}

In this paper we prove some conjectures of K. Kato \cite{Ka} which
were formulated to generalize the classical exact sequence of Brauer groups
for a global field $K$,
$$
0 \longrightarrow Br(K) \longrightarrow \mathop{\textstyle\bigoplus}\limits_{v} Br(K_v) \longrightarrow \qz \longrightarrow 0\,,
\leqno(0.1)
$$
to function fields $F$ over $K$ and varieties $X$ over $K$.
In the above sequence, which is also called the Hasse-Brauer-Noether sequence,
the sum runs over all places $v$ of $K$, and $K_v$ is the completion of $K$
with respect to $v$. The injectivity of the restriction map into the sum of local Brauer
groups is called the Hasse principle.

Kato's generalization does not concern Brauer groups but rather the following
cohomology groups. Let $L$ be any field and let $n>0$ be an integer. Write $n = mp^r$
with $p = \ch(L)$ the characteristic of $L$ and $p \nmid m$ (so that $n=m$ for $\ch(L)=0$).
Define the following Galois cohomology groups for $i,j\in \z$
$$
H^i(L,\zn(j)) := H^i(L,\mu^{\otimes j}_m) \oplus H^{i-j}(L,W_r\Omega^j_{L,log})
\leqno(0.2)
$$
where $\mu_m$ is the Galois module of $m$-th roots of unity (in the separable
closure $L^{sep}$ of $L$) and $W_r\Omega^j_{L,log}$ is the logarithmic part of
the de Rham-Witt sheaf $W_r\Omega^j_L$ \cite{Il1} I 5.7 (an \'etale sheaf, regarded
as a Galois module). It is a fact that $Br(L)[n] = H^2(L,\zn(1))$, where
$A[n] = \{\,x\in A \,\mid\,nx = 0\}$ denotes the $n$-torsion in an abelian group $A$,
so the $n$-torsion of the sequence (0.1) can be identified with an exact sequence
$$
0 \longrightarrow H^2(K,\zn(1)) \longrightarrow \mathop{\textstyle\bigoplus}\limits_{v} H^2(K_v,\zn(1)) \longrightarrow \zn \longrightarrow 0\,.
\leqno(0.3)
$$
In fact, this sequence is often used for the Galois cohomology of number fields,
independently of Brauer groups; it is closely related to class field theory and
Tate-Poitou duality.

For the generalization let $F$ be the function field in $d$ variables over a global fields $K$
and assume $F/K$ is primary, i.e., that $K$ is separably closed in $F$. For each place $v$ of $K$, let $F_v$ be the
corresponding function field over $K_v$: If $F = K(V)$, the function field of
a geometrically integral variety $V$ over $K$, then $F_v = K_v(V\times_KK_v)$.
Then Kato conjectured:

\v

\noi\textbf{Conjecture 1} The following restriction map is injective:
$$
\alpha_n: \; H^{d+2}(F,\zn(d+1)) \longrightarrow \mathop{\textstyle\bigoplus}\limits_{v} H^{d+2}(F_v,\zn(d+1))\,.
$$

\medskip

\noi
Note that this generalizes the injectivity in (0.3), which is the case $d=0$ and $F = K$.
On the other hand it is known that the corresponding restriction map for Brauer groups is not
in general injective for $d\geq 1$: If $X$ is a smooth projective curve over a number
field which has a $K$-rational point, then for $F=K(X)$ the kernel of Br$(F) \rightarrow \prod_v$ Br$(F_v)$
is isomorphic to the Tate-Shafarevich group of the Jacobian Jac$(X)$.
Kato proved Conjecture 1 for $d=1$ \cite{Ka}. Here we prove:

\v

\noi\textbf{Theorem 0.1} Conjecture 1 is true if $n$ is invertible in $K$.

\v
The proof uses three main ingredients. First we prove the analogue for infinite coefficients.
For a field $L$, a prime $\ell$ and integers $i$ and $j$ we let
$$
H^i(L,\qzl(j)) = \mathop{\mbox{lim}}\limits_{\rightarrow} H^i(L,\zln(j))\,,\leqno(0.4)
$$
where the inductive limit is taken via the obvious monomorphisms $\zln(j) \hookrightarrow \z/\ell^{n+1}(j)$.
Then we prove (see Theorem 2.10):

\v

\noi\textbf{Theorem 0.2} Let $K$ be a global field, let $\ell$ be a prime invertible in $K$, and let $F$ be
a function field in $d$ variables over $K$ such that $F/K$ is primary.
Then the restriction map
$$
\alpha_{\ell^\infty}: \; H^{d+2}(F,\qzl(d+1)) \longrightarrow \mathop{\textstyle\bigoplus}\limits_{v} H^{d+2}(F_v,\qzl(d+1))\,.
$$
is injective.

\v
For number fields and $d=2$ this result was already proved in \cite{Ja4}. Concerning the case of finite
coefficients, i.e., the original Conjecture 1, we use the following. For any field $L$,
any prime $\ell$ and any integer $t\geq 0$, there is a symbol map
$$
h^t_{L,\ell}: K^M_t(L)/\ell \longrightarrow H^t(L,\z/\ell\z(t)),
$$
where $K^M_t(L)$ denotes the $t$-th Milnor $K$-group of $L$ (\cite{Mi1} and \cite{BK} \S2).
Extending an earlier conjecture of Milnor \cite{Mi1} for $\ell = 2 \neq \ch(L)$,
Bloch and Kato stated the following conjecture:

\medskip
\noi BK$(L,t,\ell)$:  The map $h^t_{L,\ell}$ is an isomorphism.
\medskip

\noi
This conjecture was proved in the recent years. In fact, for $\ell =\ch(L)$ it was proved by
Bloch, Gabber and Kato \cite{BK}, and for $\ell\neq \ch(L)$ it is classical for $t=1$
(Kummer theory), was proved for $t=2$ by Merkurjev and Suslin \cite{MS}), for $\ell =2$
by Voevodsky \cite{V1}, and for arbitrary $t$ and $\ell$ by work of Rost and Voevodsky
(\cite{SJ}, \cite{V2}, \cite{V3}, \cite{HW}).

\v
\noi
Property BK$(F,d+1,\ell)$ for all $\ell$ dividing $n$ allows to deduce Theorem 0.1 from
Theorem 0.2 for all $\ell$ dividing $n$ as follows. One has the exact cohomology sequence
\begin{footnotesize}
$$
H^{d+1}(F,\qzl(d+1)) \mathop{\longrightarrow}\limits^{\ell^m} H^{d+1}(F,\qzl(d+1))
\rightarrow H^{d+2}(F,\z/\ell^m(d+1)) \mathop{\longrightarrow }\limits^{i} H^{d+2}(F,\qzl(d+1)),
$$
\end{footnotesize}\noi
and it follows from BK$(F,t,\ell)$ that $H^{d+1}(F,\qzl(d+1))$ is divisible. Therefore $i$
is injective, and this shows that the injectivity of $\alpha_{\ell^\infty}$ in Theorem 0.1
implies the injectivity of $\alpha_{\ell^m}$ in Conjecture 1.
It should be noted that Kato did in fact use BK($K,2,\ell$),
i.e., the Merkurjev-Suslin theorem, in his proof of Conjecture 1 for $d=1$.

\v
\noi
Finally the proof of Theorem 0.2 uses weights, i.e., Deligne's proof of the Weil conjectures, and some results on
resolution of singularities, to control the weights. Over number fields the required resolution of singularities holds
by work of Hironaka. For $\ell$ invertible in $K$ we observe that a weaker form of resolution suffices.
More precisely we use alterations, as introduced by de Jong, but in a refined version established by Gabber, see \cite{Il2}.

\v

As in the classical case $d=0$ and the case of $d=1$ (see the appendix to \cite{Ka}), and the case
of $d=2$ in \cite{Ja4}, Theorem 0.1 has applications to quadratic forms over $F$, see \cite{CTJ}:

\v

\noi\textbf{Corollary 0.3} If $F$ is a finitely generated field of characteristic zero, then
the Pythagoras number of $F$ is finite. More precisely, if $F$ is of transcendence degree $d$
over $\mathbb Q$, then any sum of squares in $F$ is a sum of $2^{d+1}$ squares, provided $d\geq 2$.

\v
\noi
The proof uses the following instance of Theorem 0.1, which only needs the Milnor conjecture,
i.e., the theorem of Voevodsky in \cite{V1}.

\v
\noi\textbf{Corollary 0.4} The restriction map
$$
H^{d+2}(F,\z/2\z) \longrightarrow \mathop{\textstyle\bigoplus}\limits_{v}\,H^{d+2}(F_v,\z/2\z)\,.
$$
is injective.

\v
It should be mentioned that the finiteness of the Pythagoras number, with the weaker bound $2^{d+2}$,
can be obtained by some more elementary means, still using the Milnor conjecture \cite{P}.

\v
Kato also stated a conjecture on the cokernel of the above restriction maps,
in the following way. Let $L$ be a global or local field, let $X$ be any
variety over $L$, and let $n$ be an integer. Then in \cite{Ka} Kato defined
a certain homological complex $C^{2,1}(X,\zn)$ of Galois cohomology groups:
\begin{multline*}
\cdots \longrightarrow \mathop{\textstyle\bigoplus}\limits_{x\in X_a} H^{a+2}(k(x),\zn(a+1))\longrightarrow
\mathop{\textstyle\bigoplus}\limits_{x\in X_{a-1}} H^{a+1}(k(x),\zn(a))\longrightarrow \cdots \\
\cdots \longrightarrow \mathop{\textstyle\bigoplus}\limits_{x\in X_1} H^{3}(k(x),\zn(2))
\longrightarrow \mathop{\textstyle\bigoplus}\limits_{x \in X_0} H^{2}(k(x),\zn(1)).
\end{multline*}
Here $X_a$ denotes the set of points $x \in X$ of dimension $a$, the term involving $X_a$ is placed in degree $a$.
and $k(x)$ denotes the residue field of $x$.
A complex of the same shape can also be defined via the method of Bloch
and Ogus, and it is shown in \cite{JSS} that these two definitions
agree up to (well-defined) signs (see also \S4 for a discussion of more
general complexes $C^{a,b}(X,\zn)$).

\v

Now let $K$ be a global field and let $X$ be a variety over $K$. Then there
are obvious maps of complexes $C^{2,1}(X,\zn) \rightarrow C^{2,1}(X_v,\zn)$
for each place $v$ of $K$, where $X_v=X\times_KK_v$, and these induce a
map of complexes
$$
\beta_{X,n}: C^{2,1}(X,\zn) \longrightarrow \mathop{\textstyle\bigoplus}\limits_{v}\;C^{2,1}(X_v,\zn)\,.
$$
Then Kato conjectured the following \cite{Ka}.

\v

\noi\textbf{Conjecture 2} Let $K$ be a global field, let $n > 0$ be an integer,
and let $X$ be a connected smooth projective variety over $K$. Then the above map
induces isomorphisms
$$
H_a(C^{2,1}(X,\zn)) \mathop{\longrightarrow}\limits^{\sim} \mathop{\textstyle\bigoplus}\limits_{v}\;H_a(C^{2,1}(X_v,\zn))
$$
for $a>0$, and an exact sequence
$$
0 \longrightarrow H_0(C^{2,1}(X,\zn)) \longrightarrow \mathop{\textstyle\bigoplus}\limits_{v}\; H_0(C^{2,1}(X_v,\zn)) \longrightarrow \zn \longrightarrow 0\,.
$$

\v

\noi
Note that we obtain a sequence of the shape (0.3) for $X = \Spec(K)$, where
the complexes are concentrated in degree zero. Kato \cite{Ka} proved Conjecture 2 for $d=1$.
Here we prove:

\v
\noi
\textbf{Theorem 0.5} Conjecture 2 is true if $K$ is a number field or if $n$ is invertible in $K$.
More precisely, in this case there is an exact sequence of complexes

$$
0 \rightarrow C^{2,1}(X,\zn) \longrightarrow \mathop{\textstyle\bigoplus}\limits_{v}\;C^{2,1}(X_v,\zn) \longrightarrow C'(X,\zn) \rightarrow 0
$$

\noi
with $H_0(C'(X,\zn)) = \zn$, and $H_a(C'(X,\zn)) = 0$ for $a > 0$.

\v

Again this version is deduced from a version with infinite coefficients by using the property
$B(L,d+1,\ell)$ (for all residue fields of $X$ and all $\ell$ dividing $n$),
and the version with infinite coefficients is proved using weight arguments and some results on
resolution of singularities. The second case was proved earlier in \cite{KS}, by a different method.

\v

Our techniques also allow to get results on another conjecture of Kato,
over finite fields. For any variety over a finite field $k$ and any natural number $n$,
Kato considered a complex $C^{1,0}(X,\zn)$ which is of the form
\begin{multline*}
\cdots \longrightarrow \mathop{\textstyle\bigoplus}\limits_{x\in X_a} H^{a+1}(k(x),\zn(a))\longrightarrow
\mathop{\textstyle\bigoplus}\limits_{x\in X_{a-1}} H^{a}(k(x),\zn(a-1))\longrightarrow \cdots \\
\cdots \longrightarrow \mathop{\textstyle\bigoplus}\limits_{x\in X_1} H^{2}(k(x),\zn(1))
\longrightarrow \mathop{\textstyle\bigoplus}\limits_{x \in X_0} H^{1}(k(x),\zn)
\end{multline*}
with the term involving $X_a$ placed in (homological) degree $a$
(this is another special case of the general complexes $C^{a,b}(X,\zn)$).
Kato conjectured the following (where the case $a=0$ is easy):

\v

\noi
\textbf{Conjecture 3} If $X$ is connected, smooth and projective over a finite field $k$,
then one has
$$
H_a(C^{1,0}(X,\zn))=\left\{\begin{array}{ccc} 0 & , & a\neq 0,\\
\zn & , & a=0\, .\end{array}\right.
$$

\v

\noi
For $\dim(X) = 1$ this conjecture amounts to the exact sequence (0.3) with $K = k(X)$, and for $\dim(X)=2$
the conjecture follows from \cite{CTSS} for $n$ invertible in $k$, and from \cite{Gro2}
and \cite{Ka} if $n$ is a power of $\ch(k)$. S. Saito \cite{Sa} proved that $H_3(C^{2,1}(X,\qzl)) = 0$
for $\dim(X)=3$ and $\ell \neq \ch(k)$. For $X$ of any dimension
Colliot-Th\'el\`ene \cite{CT1} (for $\ell\neq\ch(k)$) and Suwa \cite{Su}
(for $\ell=\ch(k))$ proved that $H_a(C^{1,0}(X,\qzl))=0$ for $0<a\leq 3$.
Here we prove the following:

\v

\noi
\textbf{Theorem 0.6} Conjecture 3 holds if resolution of singularities holds over $k$, or if $n$ is invertible in $k$.

\v

\noi
The second result was proved earlier in \cite{KS}, by a different method. Moreover, the first
result also follows with the methods of \cite{JS2}.

\v
By another technique, resolution of singularities in low dimension is used in \cite{JS2}
to show unconditionally that $H_a(C^{1,0}(X,\zn))=0$ for $X$ smooth projective of any dimension, any $n$, and $0<a\leq 4$.
Finally, Kato also formulated an arithmetic analogue of Conjecture 3, for regular
flat proper schemes of $\Spec(\z)$, and in \cite{JS} some results on this are
obtained using Theorem 0.2.

\v
Our method of proof is the same for Theorem 0.5 and 0.6. In fact, under certain conditions, which are always fulfilled in our cases,
Kato defined more general complexes $C^{r,s}(X,\zn)$ of the form
\begin{multline*}
\cdots \longrightarrow \mathop{\textstyle\bigoplus}\limits_{x\in X_a} H^{r+a}(k(x),\zn(s+a))\longrightarrow
\mathop{\textstyle\bigoplus}\limits_{x\in X_{a-1}} H^{r+a-1}(k(x),\zn(s+a-1))\longrightarrow \cdots \\
\cdots \longrightarrow \mathop{\textstyle\bigoplus}\limits_{x\in X_1} H^{r+1}(k(x),\zn(s+1))
\longrightarrow \mathop{\textstyle\bigoplus}\limits_{x \in X_0} H^{r}(k(x),\zn(s))
\end{multline*}

We show that for $n$ invertible in $K$ there is a canonical quasi-isomorphism between the complex $C'(X,\zn)$ in Theorem 0.5
and the complex $C^{0,0}(\overline{X},\zn)_{G_K}$ obtained from the Kato complex $C^{0,0}(\overline{X},\zn)$
by taking coinvariants under the absolute Galois group $G_K$, where $\overline{X}=X\times_K\overline{K}$.
On the other hand, for a finite field $k$ one has a canonical isomorphism $C^{0,0}(\overline{X},\zn)_{G_k} \cong C^{1,0}(X,\zn)$
for a variety $X$ over a finite field $k$. Therefore Theorems 0.5 and 0.6 follow from the following more general result:

\v
\noi
\textbf{Theorem 0.7} Let $K$ be a finitely generated field with algebraic closure $\overline{K}$, let
$X$ be a smooth projective variety over $K$, and let $n$ be natural number. Then
$$
H_a(C^{0,0}(X\times_K \overline{K},\zn)_{G_K}) = \left\{\begin{array}{cl}
\zn & ,~~a=0\\
0 & ,~~a\neq 0\, .
\end{array}\right. \,,
$$
if resolution of singularities holds over $K$, or if $n$ is invertible in $K$.

\v
\noi
Again the second case was proved earlier in \cite{KS}.

\v\noi

This paper had a rather long evolution time. Theorem 0.2 for number fields was obtained
in 1990, rather short after the proofs of Theorem 0.2 and Theorem 0.4 for number fields and
$\qzl$-coefficients and $d=2$ in \cite{Ja4}. In 1996, right after the appearance of \cite{GS},
it became clear to me how to obtain Theorem 0.5 (for number fields and infinite coefficients),
but a first account was only written in 2004. Meanwhile I had also noticed that these methods
allow a proof of Theorem 0.7, i.e., a proof of Kato's conjecture over finite fields, with
infinite coefficients, assuming resolution of singularities. Part of the delay was caused by the long time to
complete the comparison of Kato's original complexes with the complexes of
Gersten-Bloch-Ogus type used here, which was recently accomplished \cite{JSS}.
Most recently, inspired by the results of Kerz and Saito \cite{KS}, I saw how to use the
results of Gabber-Illusie with my method to avoid resolution of singularities for
coefficients away from the characteristic.

\v
\noi
I dedicate this paper to my teacher and friend J\"urgen Neukirch, who
helped and inspired me in so many ways by his support and enthusiasm. I also
thank Jean-Louis Colliot-Th\'el\`ene for his long lasting interest in this work,
for the discussions on the rigidity theorems 2.9 and 4.14, and for the proof of
Theorem 2.10. Moreover I thank Wayne Raskind, Florian Pop, Tam\'{a}s Szamuely
and Thomas Geisser for their interest and useful hints and discussions. In establishing
the strategy for proving Theorems 0.7 and 0.5 I profited from an incomplete preprint
by Michael Spie\ss.
My contact with Shuji Saito started with the subject of this paper, and I thank
him for all these years of a wonderful collaboration and the countless inspirations
I got from our discussions.

\vspace{15mm}

\centerline{\noindent\textbf{\S 1 First reductions and a Hasse
principle for global fields}}

\vspace{1,0cm}

Let $K$ be a global field, and let $F$ be a function field of
transcendence degree $d$ over $K$. We assume that $K$ is
separably closed in $F$. For every place $v$ of $K$, let $K_v$
be the completion of $K$ at $v$, and let $F_v$ be the
corresponding function field over $K_v$: there exists a
geometrically irreducible variety $V$ of dimension $d$ over $K$,
such that $F=K(V)$, and then $F_v=K_v(V_v)$, where $V_v=V\times_K
K_v$ (this is integral, since $F/K$ is primary and $K_v/K$ is separabel,
see [EGA IV],2, (4.3.2) and (4.3.5)).
This definition does not depend on the choice of $V$.

\medskip

Fix a prime $\ell\neq \mbox{char}(K)$. We want to study the map
$$
\mbox{res: } H^{d+2}(F,\mathbb Q_\ell/\mathbb Z_\ell (d+1))\rightarrow
\prod\limits_v
~H^{d+2}(F_v,\mathbb Q_\ell/\mathbb Z_\ell(d+1))
$$
induced by the restrictions from $F$ to $F_v$. For this it will be
useful to first replace the completions $K_v$ by the
Henselizations. For each place $v$ of $K$, denote by $K_{(v)}$ the
Henselization of $K$ at $v$. It can be regarded as a subfield of a
fixed separable closure $\overline K$ of $K$, equal to the fixed
field of a decomposition group $G_v$ at $v$. For $V$ as above, let
$F_{(v)}=K_{(v)}(V\times_K K_{(v)})$ be the corresponding function
field over $K_{(v)}$. Since $K_{(v)}$ is separably algebraic over
$K$ and linearly disjoint from $F,~ F_{(v)}$ is equal to the
composite $FK_{(v)}$ in a fixed separable closure $\overline F$
of $F$. We obtain a diagram of fields
$$
\xymatrix{ & \overline F\ar@{-}[d]\\
\overline K\ar@{-}[r]\ar@{-}[d] & F\overline K\ar@{-}[d]\\
K_{(v)} \ar@{-}[r]\ar@{-}[d] & F_{(v)}\ar@{-}[d]\\
K\ar@{-}[r] & F}\leqno(1.1)
$$
which identifies $G_K= \mbox{Gal}(\overline K/K)$ with $
\mbox{Gal}(F\overline K/F)$, and $G_{K_{(v)}
}=\mbox{Gal}(\overline K/K_{(v)} )$ with \\
$\mbox{Gal}(F\overline
K/FK_{(v)} )$.

\vspace{0,5cm}

\noindent\textbf{Proposition 1.2} Let $M$ be a discrete
$\ell$-primary torsion $G_F$-module. The restriction map
$$
H^{d+2}(F,
M)\rightarrow \prod\limits_v H^{d+2}(F_{(v)}, M)
$$
has image in the direct sum
$\mathop{\textstyle\bigoplus}\limits_{v}
H^{d+2}(F_{(v)}, M)$. There is a commutative diagram
$$
\begin{array}{cccc}
f: & H^{d+2}(F, M) & \longrightarrow &
\bigoplus\limits_v H^{d+2}(F_{(v)}, M)\\
& \downarrow & & \downarrow\\
g: & H^2(K, H^d(F\overline K, M)) & \longrightarrow &
\bigoplus\limits_v
H^2(K_{(v)}, H^d(F\overline K, M)) ,
\end{array}
$$
in which the horizontal maps are induced by the restrictions, and the vertical maps by the Hochschild-Serre spectral sequences. This diagram is functorial in $M$ and induces canonical isomorphisms
$$
\ker(f)\stackrel{\sim}{\longrightarrow} \ker (g)\quad\quad \mbox{and} \quad\quad
\mbox{coker}(f)\stackrel{\sim}{\longrightarrow} \mbox{coker}(g) \cong H^d(F\overline{K},M)(-1)_{G_K}\;
.
$$
Here $N(n)$ denotes the $n$-fold Tate twist of a $\ell$-primary discrete torsion $G_K$-module $N$,
and $N_{G_K}$ denotes its cofixed module, i.e., the maximal quotient on which $G_K$ acts trivially.
\vspace{0,5cm}

\noindent\textbf{Proof} Diagram (1.1) gives Hochschild-Serre
spectral sequences
$$
\begin{array}{rcl}
E_2^{p,q}(K)& = & H^p(K, H^q(F\overline K, M))\Rightarrow
H^{p+q}(F, M)\\
E_2^{p,q}(K_{(v)}) & = & H^p(K_{(v)}, H^q(F\overline K,
M))\Rightarrow H^{p+q}(F_{(v)}, M)
\end{array}
$$
Moreover, for each $v$ we obtain a natural map $E(K)\rightarrow
E(K_{(v)})$ between the above spectral sequences which gives the
restriction maps for $K\subset K_{(v)}$ on the $E_2$-terms and the
restriction maps for $F\subset F_{(v)}$ on the abutment,
respectively. On the other hand, the field $F\overline K$ has
cohomological dimension $d$, so that
$E_2^{p,q}(K)=0=E_2^{p,q}(K_{(v)})$ for $q>d$. This gives a
commutative diagram
$$
\begin{array}{ccc}
H^{d+2}(F, M) & \longrightarrow &
H^{d+2}(F_{(v)}, M)\\
\downarrow & & \downarrow\\
H^2(K, H^d(F\overline K, M)) & \longrightarrow & H^2(K_{(v)},
H^d(F\overline K, M)) ,
\end{array}
$$
where the vertical maps are edge morphisms of the spectral
sequences. If $v$ is not a real archimedean place, or if $\ell\neq
2$, we have $cd_\ell (K_{(v)})\leq 2$ and, hence,
$E^{p,q}_2(K_{(v)})=0$ for $p>2$, and the right vertical edge
morphism is an isomorphism. This already shows the first claim of
the proposition, since the restriction map
$$
H^{2}(K, N)\rightarrow\prod\limits_v H^{2}(K_{(v)}, N)
$$
is known to have image in the direct sum
$\mathop{\textstyle\bigoplus}\limits_{v}$ for any
torsion $G_K$-module $N$. If $K$ has no real archimedean
valuations (or if $\ell\neq 2)$, then $cd_\ell(K)=2$, the left hand
edge morphism is an isomorphism as well, and the second claim
follows.

\medskip

So let now $K$ be a number field. Here we need the following lemma.

\vspace{0,5cm}

\noindent\textbf{Lemma 1.3} The above maps between the
spectral sequences induce

\begin{itemize}

\item[(a)] surjections for all $r\geq 2$ and all $p+q=d+1$
$$
E^{p,q}_r
(K)\twoheadrightarrow\textstyle\bigoplus\limits_
{v\mid\infty}E^{p,q}_r(K_{(v)})\; ,
$$

\item[(b)] surjections for all $r\geq 2$
$$
E^{2,d}_r
(K)\twoheadrightarrow\textstyle\bigoplus\limits_
{v \mid\infty}E^{2,d}_r(K_{(v)})\; ,
$$

\item[(c)] isomorphisms between the kernels and between the cokernels of the maps
$$
E^{2,d}_r
(K)\rightarrow\textstyle\bigoplus\limits_v E^{2,d}_r(K_{(v)})\quad\mbox{and}\quad E^{2,d}_{r+1}
(K)\rightarrow\textstyle\bigoplus\limits_v E^{2,d}_{r+1}(K_{(v)})\; ,
$$
for all $r\geq 2$,

\item[(d)] isomorphisms
$$
E^{p,q}_r
(K)\stackrel{\sim}{\longrightarrow}\textstyle\bigoplus\limits_{v \mid\infty}E^{p,q}_r(K_{(v)})
$$
for all $r\geq 2$ and all $(p, q)\neq (2, d)$ with $p+q\geq
d+2$.

\end{itemize}

\medskip

\noindent\textbf{Proof:} By induction on $r$. Recall that
$E^{p,q}_r(K)=0=E^{p,q}_r(K_{(v)})$ for all $q>d$ and all $r\geq
2$. Hence, for $r=2$ the claims (a), (b) and (d) follow from the
following well-known facts of global Galois cohomology: the maps
$$
\begin{array}{ccc}
H^1(K, N) & \rightarrow &
\bigoplus\limits_{ v \mid\infty}
H^1(K_{(v)},
N)\\
H^2(K, N) & \rightarrow &
\bigoplus\limits_{ v \mid S}
H^2(K_{(v)}, N)
\end{array}
$$
are surjective for any torsion $G_K$-module $N$ and any finite set
$S$ of places, and the maps
$$
H^i(K, N)  \stackrel{\sim}{\longrightarrow}
\textstyle\bigoplus\limits_{ v \mid\infty}
H^i(K_{(v)}, N)
$$
are isomorphisms for such $N$ and all $i\geq 3$. Note that here we
could replace $K_{(v)}$ by the more common completion $K_v$, since
$G_{K_{(v)}}\cong G_v\cong G_{K_v}$.

\medskip

Now let $r\geq 2$. For (a) look at the commutative diagram
$$
\begin{CD} E_r^{p-r,q+r-1}(K) @>{d_r}>>  E^{p,q}_r(K)  @>{d_r}>> E_r^{p+r,q-r+1}(K)  \\
@V{\alpha}VV @V{\beta}VV  @V{\gamma}VV \\
\bigoplus\limits_{ v \mid\infty}E_r^{p-r,q+r-1}(K_{(v)}) @>{d_r}>>
\bigoplus\limits_{ v \mid\infty} E^{p,q}_r(K_{(v)})  @>{d_r}>>  \bigoplus\limits_{ v \mid\infty}E_r^{p+r,q-r+1}(K_{(v)})
\end{CD}
$$
coming from the map of spectral sequences. We may assume $p\geq
1$(since $E^{0,d+1}_r=0$), and hence $(p+r, q-r+1)\neq (2,d)$.
Then $\beta$ is surjective and $\gamma$ is an isomorphism, by
induction assumption (for (a) and (d)). By taking homology of both
rows, we obtain a surjection
$E^{p,q}_{r+1}(K)\twoheadrightarrow\mathop{\textstyle\bigoplus}\limits_{ v \mid\infty}E^{p,q}_{r+1}(K_{(v)})$ as wanted for (a).

\medskip

For (d) we look at the same diagram where now we may assume that
$p\geq 2, (p,q)\neq (2,d)\neq (p+r, q-r+1)$, that $\beta$ and
$\gamma$ are bijective, and that $\alpha$ is surjective (by
induction assumption for (a), (b) and (d)). Hence we get the
isomorphism
$$
E^{p,q}_{r+1}(K)\stackrel{\sim}{\longrightarrow}\textstyle\bigoplus\limits_{ v \mid\infty}E^{p,q}_{r+1}(K_{(v)}) .
$$

\medskip

For (b) and (c) consider the \emph{exact} commutative diagram
$$
\xymatrix{ 0\ar[r] & E^{2,d}_{r+1}(K)\ar[r]\ar[d]^{\raisebox{0,5cm}{$\scriptstyle\beta'$}} &
E^{2,d}_{r}(K)\ar[r]^{d_r}\ar[d]^{\raisebox{0,5cm}{$\scriptstyle\beta$}} &
E^{2+r,d-r+1}_{r}(K)\ar[d]^{\raisebox{0,5cm}{$\scriptstyle\gamma$}}_{\raisebox{0,5cm}{$\scriptstyle\wr$}}\\
0 \ar[r] & \bigoplus\limits_{ v \in S'}E^{2,d}_{r+1}(K_{(v)})\ar[r] &
\bigoplus\limits_{ v \in S'}E^{2,d}_{r}(K_{(v)}) \ar[r]^{\hspace{-5mm}\partial} &
\bigoplus\limits_{ v \mid\infty}E^{2+r,d-r+1}_{r}(K_{(v)})}
$$
for any set of places $S'\supset\{v\!\mid\!\infty\}$, in which $\partial = \mathop{\textstyle\bigoplus}\limits_{v\in S'}d_r(K_{(v)})$(note that
$E^{p,q}_r(K_{(v)})=0$ for $p>2$ and $v\nmid \infty)$. The map
$\gamma$ is an isomorphism by induction assumption (for (d)). Hence
for $S'=\{v\!\mid\!\infty\}$ the surjectivity of $\beta$ implies
the one for $\beta'$, i.e., we get (b) for $r+1$. If $S'$ is the
set of all places, we see that clearly $\ker(\beta')=\ker(\beta)$, and
that coker$(\beta') =$ coker$(\beta)$, since im$(\partial\circ \beta)=$
im$(\partial)$ by induction assumption for (b). Thus we get (c) for
$r$ from (b) and (d) for $r$.

\vspace{1cm}

\noindent We use the lemma to complete the proof of proposition 1.2. From what we
have shown, we have $E^{0,d+2}_\infty = E^{0,d+1}_\infty =0$ for
$K$ and all $K_{(v)}$, and isomorphisms
$$
E^{p,q}_{\infty}(K)\stackrel{\sim}{\longrightarrow}\textstyle\bigoplus\limits_v E^{p,q}_{\infty}(K_{(v)}) .
$$
for all $(p,q)$ with $p+q=d+2, p\geq 3$ (note that $E^{p,q}_2 (K_{(v)})=0$ for
$p\geq 3$ and $v\nmid\infty$). Hence kernel and cokernel of
$$
H^{d+2}(F, M)\rightarrow
\textstyle\bigoplus\limits_
v H^{d+2}(F_{(v)}, M)
$$
can be identified with kernel and cokernel of
$$
E^{2,d}_\infty(K)\rightarrow
\textstyle\bigoplus\limits_v
E^{2,d}_\infty(K_{(v)}) ,
$$
respectively. But these coincide with kernel and cokernel of
$$
E^{2,d}_2(K)=H^2(K,H^d(F\overline K, M))\rightarrow
\textstyle\bigoplus\limits_
v H^2(K_{(v)}, H^d(F\overline K, M))=
\textstyle\bigoplus\limits_
vE^{2,d}_2(K_{(v)})
$$
respectively, by (c) of the lemma. Finally, for any finite $\ell$-primary
$G_K$-module $N$, Poitou-Tate duality gives an exact sequence
$$
H^2(K,N) \longrightarrow \mathop{\textstyle\bigoplus}\limits_{v} H^2(K_{(v)},N) \longrightarrow H^0(K,N^\ast)^\vee \longrightarrow 0\,,
$$
where $N^\ast$ denotes the finite $G_K$-module $Hom(N,\mu)$ where $\mu$ is the group
of roots of unity in $\overline{K}$ and $M^\vee$ is the Pontrjagin dual of a finite
$G_K$-module. But then we have canonical identifications
$$
\begin{array}{rl}
      & H^0(K,Hom(N,\mu))^\vee =  Hom_{G_K}(N,\qzl(1))^\vee \\
\cong & Hom_{G_K}(N(-1),\qzl)^\vee = Hom(N(-1)_{G_K},\qzl)^\vee = (N(-1)_{G_K})^{\vee\vee} \cong N(-1)_{G_K}\,.
\end{array}
$$
This shows the last isomorphism of 1.2

\vspace{0,5cm}

\noindent\textbf{Remarks 1.4} (a) Proposition 1.2 extends to the
case where $F$ is a function field over $K$, but $K$ is not
necessarily separably closed in $F$, by replacing $F_{(v)}$ with
$F\otimes_KK_{(v)}$ and $F\overline{K}$ with $F\otimes_K\overline
K$. The cohomology groups of these rings have to be interpreted as
the \'etale cohomology groups of the associated affine schemes;
with this the proof carries over verbatim. In more down-to-earth
(but more tedious) terms, we may note that $(F\otimes_K\overline K)_{red}
\cong \prod_\sigma(F\otimes_{\tilde{K},\,\sigma}\overline{K})$,
where $\tilde{K}$ is the separable closure of $K$ in $F$ and
$\sigma$ runs over the $K$-embeddings of $\tilde{K}$ into $\overline K$.
Similarly, $F\otimes_KK_{(v)} \cong
\prod_\sigma(\prod_{w}F\otimes_{\tilde{K}}\sigma(\tilde{K})_{(\sigma w)})$,
where $w$ runs over the places of $\tilde{K}$ above $v$, $\sigma w$ is the corresponding
place of $\sigma(\tilde{K})$ above $v$, and
$\sigma(\tilde{K})_{(\sigma w)}$ is the Henselization of $\sigma(\tilde{K})$ at $\sigma w$.
The \'etale cohomology groups referred to above can thus be identified with sums of Galois
cohomology groups of the fields introduced above, and the claim
also follows by applying Proposition 1.2 to $F/\tilde{K}$.

(b) A consequence of Proposition 1.2 is that the restriction map
$$
f':H^{d+2}(F,M)\rightarrow \prod\limits_vH^{d+2}(F_v, M)
$$
has image in the direct sum
$\mathop{\textstyle\bigoplus}\limits_
v\subset\prod_v $ as well, since it factors through the map $f$ in 1.2.
Moreover, as we shall see in \S2, the maps $H^{d+2}(F_{(v)},
M)\rightarrow H^{d+2}(F_v, M)$ are injective, so that
$\ker(f')=\ker(f)$. For $d>0$, however, $H^{d+2}(F_v, \mathbb
Q_\ell/\mathbb Z_\ell(d+1))$ is much bigger than
$H^{d+2}(F_{(v)},\mathbb Q_\ell/\mathbb Z_\ell(d+1))$, and
proposition 1.2 does not extend to the completions. In particular,
everywhere in [Ja4] the completions $K_v$ should be replaced by
the Henselizations $K_{(v)}$ (In loc. cit., Proof of Th. 1' and later,
the notations $FK_v$ and $F\overline K_v$
are problematic; they should be interpreted as $F_v$ and $F_v\overline{K_v}$.
Even then Gal$(\overline {FK_v}/F\overline{K_v})$
= Gal$(\overline{F_v}/F_v\overline{K_v})$ is not isomorphic to
Gal$(\overline F/F\overline K)$, but much bigger, as was kindly
pointed out to me by J.-L. Colliot-Th\'{e}l\`{e}ne and J.-P.
Serre). The comparison of coker$(f')$ and coker$(f)$ is more subtle,
see \S4.

\vspace{1cm}

By Proposition 1.2, the restriction map
$$
H^{d+2}(F, \mathbb Q_\ell/\mathbb Z_\ell(d+1))\rightarrow
\textstyle\bigoplus\limits_
v H^{d+2}(F_{(v)},\mathbb Q_\ell/\mathbb Z_\ell(d+1))
$$
has the same kernel and cokernel as the restriction map
$$
\beta_N:H^2(K, N)\rightarrow
\textstyle\bigoplus\limits_
v H^2(K_{(v)},
N)\cong\textstyle\bigoplus\limits_
v H^2(K_v, N)
$$
for the $G_K$-module $N=H^d(F\overline K, \mathbb Q_\ell/\mathbb Z_\ell(d+1))$.
Here we have used the isomorphism
$G_{K_v}\stackrel{\sim}{\longrightarrow} G_{K_{(v)}}$ to rewrite
the latter map in terms of the more familiar completions $K_v$.
Recall that $F=K(V)$, the function field of a geometrically
irreducible variety $V$ of dimension $d$ over $K$. From this we
obtain
$$
H^d(F\overline K, \mathbb Q_\ell/\mathbb
Z_\ell(d+1))=\lim\limits_{\overrightarrow{U\subset V}} H^d_{\mbox{\scriptsize \'{e}t}}(U\times_K\overline
K,\, \mathbb Q_\ell/\mathbb Z_\ell(d+1))
$$
where the limit is over all affine open subvarieties $U$ of $V$. In fact
\'{e}tale cohomology commutes with this limit ([Mi2] III 1.16 ), so that the
right hand side is the \'etale cohomology group $H^d_{\mbox{\scriptsize
\'{e}t}}(\mbox{Spec}(\overline K(V)), \mathbb Q_\ell/\mathbb
Z_\ell(d+1))$, which can be identified with the Galois cohomology
group on the left hand side. Since
$$
H^2(K,\lim\limits_{\overrightarrow{\phantom{aaaa}}} N_i)=
\lim\limits_{\overrightarrow{\phantom{aaaa}}} H^2(K, N_i)
$$
for a direct limit of
$G_K$-modules $N_i$, and since the same holds for $\bigoplus_v\,H^2(K_v,-)$, it
thus suffices to study the maps
$$
\beta_B:H^2(K,
B)\rightarrow\textstyle\bigoplus\limits_
v H^2(K_{v}, B)
$$
for $B=H^d(U\times_K\overline K, \mathbb Q_\ell/\mathbb
Z_\ell(d+1))$, where $U\subseteq V$ runs through all open
subvarieties of $V$, or through a cofinal set of them. For this we
shall use the following Hasse principle, which generalizes [Ja2]
Theorem 3.

\vspace{0,5cm}

\noindent\textbf{Theorem 1.5} Let $K$ be a global field,
and let $\ell\neq\mbox{char}(K)$ be a prime number.

\begin{itemize}

\item[(a)] Let $A$ be a discrete $G_K$-module which is isomorphic
to $(\mathbb Q_\ell/\mathbb Z_\ell)^m$ for some $m$ as an abelian
group, and mixed of weights $\neq -2$ as a Galois module. Then the
restriction map induces isomorphisms
$$
\beta_A:H^2(K,
A)\stackrel{\sim}{\longrightarrow}\textstyle\bigoplus\limits_
v H^2(K_{v},
A)=\textstyle\bigoplus\limits_{ v\in S
\mbox{ \scriptsize or } v\mid\ell}H^2(K_{v}, A) ,
$$
where $S$ is a finite set of bad places for $A$.

\item[(b)] Let $T$ be a finitely generated free $\mathbb
Z_\ell$-module with continuous action of $G_K$ making $T$ mixed of
weights $\neq 0$. Then for any finite set $S'$ of places of $K$
the restriction map in continuous cohomology
$$
\alpha_T: H^1(K,
T)\rightarrow\prod\limits_{v\notin S'}H^1 (K_v, T)
$$
is injective.

\end{itemize}

\medskip

Before we prove this, let us explain the notion of a mixed
$G_K$-representation and a bad place $v$ for it. A priori, this is
defined for a $\mathbb Q_\ell$-representation $V$ of $G_K$ (i.e.,
a finite-dimensional $\mathbb Q_\ell$-vector space with a
continuous action of $G_K$) - see {De2] (1.2) and (3.4.10), and 1.6
below. We extend it to a module like $A$ above or more generally,
to a discrete $\ell$-primary torsion $G_K$-module of cofinite type
(resp. to a finitely generated $\mathbb Z_\ell$-module $T$ with
continuous action of $G_K$), by calling $A$ (resp. $T$) pure of
weight $w$ or mixed, if this holds for the $\mathbb
Q_\ell$-representation $T_\ell A\otimes_{\mathbb Z_\ell}\mathbb
Q_\ell$ (resp. $T\otimes_{\mathbb Z_\ell}\mathbb Q_\ell)$, where
$T_\ell A=\raisebox{-6pt}{$\stackrel{\textstyle\lim}{\scriptstyle\longleftarrow}$}\phantom{}_n \,A[\ell^n]$ is the Tate module of
$A$. In the same way we define the bad places for $A$ (resp. $T$)
to be those of the associated $\mathbb Q_\ell$-representations. It
remains to recall

\vspace{0,5cm}

\noindent\textbf{Definition 1.6} (a) A $\mathbb
Q_\ell$-representation $V$ of $G_K$ is pure of weight $w\in
\mathbb Z$, if there is a finite set $S\supset\{v\!\mid\!\infty\}$
of places of $K$ such that

\begin{itemize}

\item[(i)] $V$ is unramified outside $S\cup \{v\!\mid\!\ell\}$,
i.e., for $v\notin S, v\nmid\ell$, the inertia group $I_v$ at $v$
acts trivially on $V$,

\item[(ii)] for every place $v\notin S, v\nmid \ell$, the
eigenvalues $\alpha$ of the geometric Frobenius $Fr_v$ at $v$
acting on $V$ are pure of weight $w$, i.e., algebraic numbers with
$$
\mid \!\iota\alpha\!\mid = (Nv)^{\frac{w}{2}}
$$
for every embedding $\iota: {\mathbb Q}(\alpha)\hookrightarrow
\mathbb C$, where $Nv$ is the cardinality of the finite residue
field of $v$.

\end{itemize}

\medskip

\noindent Every such set $S$ will be called a set of bad places
for $V$; the places not in $S$ are called good.

\medskip

\noindent (b) $V$ is called mixed, if it has a filtration
$0=V_0\subset V_1\subset \ldots \subset V_n=V$ by
subrepresentations such that every quotient $V_i/V_{i-1}$ is pure
of some weight $w_i$. The weights and bad places of $V$ are those
present in some non-trivial quotient $V_i/V_{i-1}$.

\vspace{0,5cm}

\noindent\textbf{Remarks and examples 1.7} For a field $L$ denote by $\overline{L}$
its separable closure, and let $G_L = Gal(\overline{L}/L)$ be its absolute Galois group.

\noindent (a) If $v$ is a place
of $K$, any extension $w$ of $v$ to $\overline K$ determines a
decomposition group $G_w\subset G_K$ and an inertia group
$I_w\subset G_w$. The arithmetic Frobenius $\varphi_w$ is a
well-defined element in $G_w/I_w$; under the canonical isomorphism
$G_w/I_w\stackrel{\sim}{\longrightarrow}\mbox{Gal}(k(w)/k(v))$ it
corresponds to the automorphism $x\mapsto x^{Nv}$ of $k(w)$. The
geometric Frobenius $Fr_w$ is the inverse of $\varphi_w$. If $I_w$
acts trivially on $V$, then the action of $Fr_w$ on $V$ is
well-defined. If we do not fix a choice of $w$, everything is
well-defined up to conjugacy in $G_K$, and we use the notation
$G_v, I_v$, and $Fr_v$. Thus ``$I_v$ acts trivially'' means that
one and hence any $I_w$ , for $w\!\mid\! v$, acts trivially, and
then the eigenvalues of $Fr_v$ are well-defined, since they depend
only on the conjugacy class.

\medskip\noindent (b) If $V$ is pure of weight $w$, then the same
holds for every $\mathbb Q_\ell$-$G_K$-subquotient. If $V'$ is
pure of weight $w'$, then $V\otimes_{\mathbb Q_\ell} V'$ is pure
of weight $w+w'$.

\medskip\noindent (c) The representation $\mathbb Q_\ell(1)$ is
unramified outside $S=\{v\!\mid\!\infty\cdot \ell\}$, and for
$v\notin S, \varphi_v$ acts on $\mathbb Q_\ell(1)$ by
multiplication with $Nv$. Therefore $\mathbb Q_\ell(1)$ is pure
of weight $-2$, and $\mathbb Q_\ell(i)$ is pure of weight $-2i$.

\medskip

\noindent (d) Let $A$ or $T$ or $V$ be $G_K$-representations as above
Definition 1.6, which are mixed of weights $\neq 0$. Then $V^{G_K} = 0 = V_{G_K}$
and $T^{G_K}=0$ and $A_{G_K}=0$. The first statement is easily reduced to the
pure case, where it follows from the fact that the eigenvalues of $Fr_v$ as in 1.6 (ii)
are different from 1. The other claims follow from the injection
$T \hookrightarrow T\otimes_{\z_\ell} \q_\ell$ and the surjection
$T_\ell A \otimes_{\z_\ell} \q_\ell \twoheadrightarrow A$.

\medskip

\noindent (e) If $X$ is a smooth and proper variety over $K$, then
the $i$-th \'{e}tale cohomology group
$H_{\mbox{\scriptsize\'{e}t}}^i(\overline X, \mathbb Q_\ell)$ of
$\overline X=X\times_K\overline K$ is pure of weight $i$ by the
smooth and proper base change theorems and by Deligne's proof of
the Weil conjectures over finite fields (cf., e.g., [Ja1] proof
of Lemma 3). The set $S$ can be taken to be the set of places
where $X$ has bad reduction, i.e., such that for $v\notin S,X$ has
good reduction at $v$, viz., a smooth projective model $\mathcal
X_v$ over $\mathcal O_v$, the ring of integers in $K_v$, with
$\mathcal X_v\times_{\mathcal O_v}K_v=X_v$.

\medskip

\noindent (f) For later purposes we note that the whole theory above has a generalization
to an arbitrary finitely generated field $K$ (see \cite{De2} (3.4.10)). A $\q_\ell$-representation
$V$ of $G_K$ (for $\ell\neq \ch(K))$ is called pure of weight $w$, if there is a normal scheme $T$ of finite type over
$\z$ with fraction field $K$ such that $V$ comes from a $\q_\ell$-representation of the algebraic
fundamental group $\pi(T,\Spec(\overline{K}))$ via the natural epimorphism $G_K \rightarrow \pi(T,\Spec(\overline{K}))$
(i.e., from a smooth $\q_\ell$-sheaf on $T$) such that for any closed
point $t\in T$ with residue field $k(t)$ of characteristic $\neq \ell$ the eigenvalues of
the geometric Frobenius $Fr_t$ are pure of weight $w$ in the sense of 1.6 (i) (replace $Nv$
by $Nt$, the cardinality of the residue field $k(t)$ of $t$, which is finite. The geometric Frobenius
$Fr_t$ is the image of the geometric Frobenius under the homomorphism
$G_{k(t)}=\pi(\Spec(k(t)),\overline{k(t)}) \rightarrow \pi(T,\Spec(\overline{K}))$
which is well-defined up to conjugation. The other notions (mixed representations, the notions for $T$ and
$A$) extend literally, as well as the properties (b) to (e) above. In (e) one takes $T$ such that $X/K$
extends to a smooth proper model $\pi: \mathcal X \rightarrow T$, and uses the base change isomorphism
$$
H^i(\overline{X},\q_\ell) \cong H^i(\mathcal X_t\times_{k(t)}\overline{k(t)},\q_\ell)
$$
where $\mathcal X_t = \mathcal X\times_Tk(t)$ is the fiber of $\pi$ over $t\in T$.

\medskip

\noindent (g) Moreover we note that there is even an analogue for a finitely generated field
$K$ and $\ell = p = \ch(K) > 0$. First we note that the notions of pure and mixed representations
still make sense, and that properties (a), (b) and (d) also hold in this situation, while (c) does
not have any counterpart. On the other hand, one has the following analogue of (e).
For a scheme $Z$ of finite type over a perfect field $L$ and $m \in \mathbb N$ let
$$
H^i(Z,\z/p^m\z(j)) := H^{i-j}(Z,W_m\Omega^j_{X,log})
$$
be the \'etale cohomology of the logarithmic part $W_m\Omega^j_{X,log}$ of the de Rham-Witt
sheaf $W_m\Omega^j_X$ (see \cite{Il1} I 5.7 and compare (0.2)). Moreover let
$$
H^i(Z,\q_p(j)) = H^i(Z,\z_p(j))\otimes_{\z_p}\q_p \quad\mbox{  where  }\quad H^i(Z,\z_p(j)) = \mathop{lim}\limits_{\leftarrow\; m} H^i(Z,\z/p^m(j))\;,
$$
with the inverse limit taken with respect to the natural epimorphisms $W_{m+1}\Omega^j_{X,log} \rightarrow W_m\Omega^j_{X,log}$.
Then for $X$ smooth and proper over a finite field $k$ of characteristic $p$, the $\q_p$-$G_k$-representation
$H^i(\overline{X},\q_p(j))$ is finite-dimensional, and it follows from the work of Deligne \cite{De1}, Katz-Messing \cite{KM}
and Milne \cite{Mi3} that it is pure of weight $i-2j$, cf. \cite{Ja6} section 3.
If $X$ is smooth and proper over a finitely generated field $K$ of characteristic $p$
and $\pi: \mathcal X \rightarrow T$ is a smooth proper model as in (f) (so that $T$ is of finite type over $\mathbb F_p$),
then Gros and Suwa (\cite{GrS} Th\'eoreme 2.1) established base change isomorphisms
$$
H^i(X\times_K\overline{K},\q_p(j)) \cong H^i(\mathcal X_t\times_{k(t)}\overline{k(t)},\q_p(j))
$$
for all closed points $t$ in a non-empty open $U \subset T$, where $\overline{K}$ now stands for
an algebraic closure of $K$. These isomorphisms are compatible with the actions
of the absolute Galois groups $G_K$ (on the left) and $G_{k(t)}$ (on the right), so that the representation
$H^i(\overline{X},\q_p(j)) = H^i(X\times_K\overline{K},\q_p(j))$ is pure of weight $i-2j$ in exactly the same sense
as for the $\ell$-adic case in (f).
Here we regard $G_K$ as the Galois group $Gal(\overline{K}/K^{in})$, where $K^{in}$ is the maximal inseparable
extension of $K$ in $\overline{K}$.

\vspace{0,5cm}

\noindent\textbf{Proof of theorem 1.5:} Part (a) is implied by
(b). In fact, $A$ is mixed of weights $\neq -2$ if and only if its
Kummer dual $T=\mbox{Hom}(A, \mu)$ (where $\mu$ is the Galois
module of roots of unity in $\overline K$) is mixed of weights
$\neq 0$, and the kernels of $\beta_A$ and $\alpha_T$ for
$S'=\emptyset$ are dual to each other by the theorem of
Tate-Poitou (and passing to the limits over the finite modules
$A[\ell^n]$ and $T/\ell^nT=\mbox{Hom}(A[\ell^n], \mu)$,
respectively). Moreover, by Tate-Poitou the cokernel of $\alpha_A$
is isomorphic to $H^0(K, T)^\vee \cong A(-1)_{G_k}$, and this is
zero by the hypothesis on the weights. Finally, by local Tate
duality, $H^2(K_v, A)$ is dual to $H^0(K_v, T)$, and for good
places $v\nmid \ell$ this is zero if $T$ is mixed of weights $\neq
0$

\medskip

Part (b) generalizes [Ja2] Theorem 3 a), which covers the case of
a pure $T$. The generalization follows by induction: let
$$
0\rightarrow T'\rightarrow T\rightarrow T''\rightarrow 0
$$
be an exact sequence of $\mathbb Z_\ell$-$G_K$-modules as in (b),
and let $S'$ be a finite set of primes. Then there is a
commutative diagram with exact rows
$$
\xymatrix{\prod\limits_{v\notin S'}H^0 (K_v, T'')\ar[r] &
\prod\limits_{v\notin S'}H^1 (K_v, T')\ar[r] & \prod\limits_{v\notin S'}H^1 (K_v, T)\ar[r] &
\prod\limits_{v\notin S'}H^1(K_v, T'')\\
H^0(K, T'')\ar[u]\ar[r] & H^1(K, T')\ar[u]_{\raisebox{-0,7cm}{$\scriptstyle\beta_{T'}$}}\ar[r] &
H^1(K, T)\ar[u]_{\raisebox{-0,7cm}{$\scriptstyle\beta_{T}$}}\ar[r] & H^1(K,
T'')\ar[u]_{\raisebox{-0,7cm}{$\scriptstyle\beta_{T''}$}} }
$$

%$$
%\begin{CD}
%\prod\limits_{v\notin S'}H^0 (K_v, T'') @>>> \prod\limits_{v\notin S'}H^1 (K_v, T') @>>>
%\prod\limits_{v\notin S'}H^1 (K_v, T) @>>> \prod\limits_{v\notin S'}H^1(K_v, T'')\\
%@AAA   @AA{\beta_{T'}}A   @AA{\beta_T}A   @AA{\beta_{T''}}A  \\
%H^0(K, T'')  @>>>  H^1(K, T')  @>>>  H^1(K, T)  @>>>  H^1(K,T'')
%\end{CD}
%$$

If $\beta_{T''}$ is injective and $H^0(K_v, T'')=0$ for all
$v\notin S'$( which is the case for $T''$ pure of weight $\neq 0$
and $S'$ containing all bad places for $T''$ and all
$v\!\mid\!\ell$, by loc. cit.), then $\beta_T$ is injective if and
only if $\beta_{T'}$ is. Since we may always enlarge the set $S'$,
the proof proceeds by induction on the length of a filtration with
pure quotients, which exists on $T\otimes_{\mathbb Z_\ell}\mathbb Q_\ell$, by definition, and hence on $T$ by pull-back.

\newpage %\vspace{15mm}

\centerline{\noindent\textbf{\S 2 Injectivity of the global-local map for coefficients invertible in $K$}}

\vspace{1,0cm}

Let $K$ be a global field, let $\ell\neq$ char$(K)$ be a prime,
and let $U$ be a smooth, quasi-projective, geometrically
irreducible variety of dimension $d$ over $K$. Following the
strategy of section 1, we study the $G_K$-module $H^d(\overline
U,\, \mathbb Q_\ell/\mathbb Z_\ell)$. Assume the following condition,
which holds for number fields by Hironaka's resolution of
singularities in characteristic zero [Hi].

\begin{itemize}
\item[]
{\bf RS2($U$):} There is a {\it good compactification} for $U$, i.e., a
smooth projective variety $X$ over $K$ containing $U$ as an open subvariety such that
$Y=X\smallsetminus U$, with its reduced closed subscheme
structure, is a divisor with simple normal crossings.
\end{itemize}

Recall that $Y$ is said to have simple normal crossings if its
irreducible components $Y_1,\ldots, Y_N$ are smooth projective
subvarieties $Y_i \subset X$ such that for all
$1\leq i_1<\ldots< i_\nu\leq N$, the
$\nu$-fold intersection $Y_{i_1,\ldots,i_\nu}:=Y_{i_1}\cap\ldots
\cap Y_{i_\nu}$ is empty or smooth projective of pure dimension
$d-\nu$, so the same is true for the disjoint union
$$
Y^{[\nu]}:=\mathop{\coprod}_{1\leq i_1<\ldots< i_\nu\leq N}
Y_{i_1,\ldots,i_\nu}\qquad (1\leq\nu\leq d),
$$
and for $Y^{[0]}:=Y_\emptyset:= X.$

\medskip

This geometric situation gives rise to a spectral sequence
$$
E^{p,q}_2=H^p(\overline{Y^{[q]}},\mathbb Q_\ell(-q))\Rightarrow
H^{p+q}(\overline U,\mathbb Q_\ell)\,,\leqno(2.1)
$$
see, e.g., [Ja2] 3.20. It is called the weight spectral sequence, because it
induces the weight filtration on the $\ell$-adic representation
$H^n(\overline U,\mathbb Q_\ell)$. In fact, $E^{p,q}_2$ is pure of
weight $p+2q$. Therefore the same is true for the
$E^{p,q}_\infty$-terms, and if $\tilde W_q$ denotes the canonical
ascending filtration on the limit term $H^n(\overline U,\mathbb
Q_\ell)$ for which $\tilde W_q/ \tilde W_{q-1}=E^{n-q,q}_\infty$,
then its $n$-fold shift $W_\cdot:=\tilde W_\cdot[-n]$ (i.e., $W_i=\tilde
W_{i-n})$ is the unique weight filtration : $W_i/W_{i-1}\cong
E_\infty^{2n-i,i-n}$ is pure of weight $i$. Moreover, for $r>3$
the differentials
$$
d_r^{p,q}:E_r^{p,q}\longrightarrow E^{p+r,q-r+1}_r
$$
are morphisms between Galois $\mathbb Q_\ell$-representation of
different weights (viz., $p+2q$ and $p+2q-r+2$) and hence vanish,
so that $E_\infty^{p,q}=E_3^{p,q}.$

\medskip

Note that $E_2^{p,q}=0$ for $p<0$ or $q<0$. Hence the weights
occurring in $H^n(\overline U,\mathbb Q_\ell)$ lie in
$\{n,\ldots,2n\},W_{2n-1}$ is mixed of weights $w\leq 2n-1$, and
$$
\begin{array}{ccl}
H^n(\overline U,\mathbb Q_\ell)/W_{2n-1} & = & W_{2n}/W_{2n-1}
=E^{0,n}_3 \\& =&  \mbox{ker}(H^0( \overline{Y^{[n]}},\mathbb
Q_\ell(-n))\stackrel{d_2^{0,n}}{\longrightarrow}
H^2(\overline{Y^{[n-1]}},~\mathbb Q_\ell(-n+1)).
\end{array}
$$
In particular, the Galois action on $(W_{2n}/W_{2n-1})(n)$ factors
through a finite quotient, since this is the case for $H^0
(\overline{Y^{[n]}},\mathbb Q_\ell).$

\medskip

We want to say something similar for $H^n(\overline U,\mathbb
Q_\ell/\mathbb Z_\ell)$, at least for $n=d\;(=\dim~U)$. If $U $ is
affine, then we have an exact sequence
$$
\ldots\rightarrow H^d(\overline U,\mathbb Z_\ell)\rightarrow
H^d(\overline U,\mathbb Q_\ell) \rightarrow H^d(\overline
U,\mathbb Q_\ell/\mathbb Z_\ell)\rightarrow 0,
$$
since $H^{d+1}(\overline U,\mathbb Z_\ell)=0$ by weak Lefschetz
[Mi2] VI 7.2. Thus $B_1=H^d(\overline U,\mathbb Q_\ell/\mathbb
Z_\ell)$ is divisible, and there is an exact sequence
$$
0\longrightarrow A_1\longrightarrow B_1\longrightarrow
C_1\longrightarrow 0
$$
in which $A_1=\mbox{im}(W_{2d-1} H^d(\overline U,\mathbb
Q_\ell)\longrightarrow H^d(\overline U,\mathbb Q_\ell/\mathbb
Z_\ell))$ is divisible and of weights $w\in\{d,\ldots ,2d-1\}$,
and in which $C_1$ is a quotient of $H^d(\overline U,\mathbb
Q_\ell)/W_{2d-1}$, divisible and pure of weight $2d$. We need to
know $C_1$ precisely, not just up to isogeny,
and this requires more arguments - note that in the $\mathbb
Q_\ell/\mathbb Z_\ell$-analogue of (2.1) the differentials
$d_1^{p,q}$ will not in general vanish for $r\geq 3$.

\medskip

For a better control of this spectral sequence we replace $U$ by a
smaller variety, as follows. By the Bertini theorem, there is a
hyperplane $H$ in the ambient projective space whose intersection
with $X$ and all $Y_{i_1,\ldots,i_\nu}$ is transversal, i.e.,
gives smooth divisors in these (in particular, the intersection
with $Y_{i_1,\ldots,i_d}$ is empty for all $d$-tuples
$(i_1,\ldots,i_d))$. This means that $\tilde Y=
\mathop{\cup}\limits_{i=1}^{N+1} Y_i$, with $Y_{N+1}:= H\cap X$,
is again a divisor with strict normal crossings on $X$. As
explained in section 1, it is possible for our purposes to replace
$U$ by the open subscheme $U^0=X\smallsetminus \tilde Y=
U\smallsetminus (H\cap U)$, because such subschemes form a cofinal
subset in the set of all opens $U\subseteq V,\, F=K(V)$. Now we
have the following description for $B_0:=H^d
(\overline{U^0},\mathbb Q_\ell/\mathbb Z_\ell)$.

\vspace{0,5cm}

\noindent\textbf{Proposition 2.2} There is an exact sequence
$$
0\rightarrow A_0\rightarrow H^d(\overline{U^0},\mathbb
Q_\ell/\mathbb Z_\ell) \rightarrow C_0\rightarrow 0
$$
in which $A_0$ is divisible and mixed of weights in $\{d,\ldots
,2d-1\}$, and
$$
C_0=I\otimes_\mathbb Z \mathbb Q_\ell/\mathbb Z_\ell (-d)
$$
for a finitely generated free $\mathbb Z$-module $I$ with discrete
action of $G_K$. Moreover, there is an exact sequence
$$0\rightarrow I'\rightarrow I\rightarrow I''\rightarrow 0$$
of $G_K$-modules with
$$
I''=\mathbb Z[\pi_0(\overline{Y^{[d]}})]
$$
$$
I'=\ker (\mathbb Z[\pi_0(\overline{Y^{[d-1]}\cap H)}]
\stackrel{\beta}{\twoheadrightarrow} \mathbb
Z[\pi_0(\overline{Y^{[d-1]}})])
$$
where $Y^{[d-1]}\cap H:= \mathop{\coprod}\limits_{1\leq i_1<\ldots
<i_{d-1}\leq N} Y_{i_1,\ldots, i_{d-1}}\cap H$, and where $\beta$
is induced by the inclusions $Y_{i_1,\ldots,i_\nu}\cap
H\hookrightarrow Y_{i_1,\ldots,i_\nu}$.

\vspace{0,5cm}

\noindent\textbf{Proof} For $1\leq i_1 <\ldots < i_\nu \leq N$
define
$$
Y^0_{i_1,\ldots,i_\nu}:= Y_{i_1,\ldots,i_\nu}\smallsetminus
(Y_{i_1,\ldots,i_\nu}\cap H)
$$
by removing the smooth hyperplane section with $H$, and let
$Y^{0[\nu]}\subseteq Y^{[\nu]}$ be the disjoint union of these
open subvarieties for fixed $\nu~(\mbox{with}~Y^0{^{[0]}} :=X^0:=
X\smallsetminus (X\cap H))$. Then
$Y^0=\mathop{\cup}\limits^N_{i=1}Y^0_i$ is a divisor with (strict)
normal crossing on $X^0$ with $U^0=X^0\smallsetminus Y^0$, and
hence there is a spectral sequence
$$
E^{p,q}_2=H^p(\overline{{Y^0}{^{[q]}}},\mathbb Q_\ell/ \mathbb
Z_\ell(-q))\Rightarrow H^{p+q}(\overline{U^0}, \mathbb Q_\ell/
\mathbb Z_\ell)\leqno(2.3)
$$
by the same arguments as for (2.1) (the properness is not needed
in the proof).

\medskip

But now the $Y^0_{i_1,\ldots,i_q}$ are affine varieties, as
complements of hyperplane sections, and of dimension $d-q$, so
that
$$
H^p(\overline{{Y^0}{^{[q]}}},\mathbb Q_\ell/\mathbb Z_\ell(-q))=
\left\{\begin{array}{ll} 0 & \mbox{for } p>d-q ,\\
\mbox{divisible} & \mbox{for } p=d-q\, ,
\end{array}\right.
$$
by weak Lefschetz. Moreover, by the Gysin sequences
$$
\ldots\rightarrow H^p(\overline{Y_{\underline{i}}},\mathbb
Q_\ell)\rightarrow H^p({\overline{Y^0_{\underline i}}},\mathbb
Q_\ell) \rightarrow H^{p-1}(\overline{Y_{\underline i}\cap
H},\mathbb Q_\ell(-1)) \rightarrow\ldots\; ,
$$
$H^p(\overline{{Y^0}{^{[q]}}},\mathbb Q_\ell)$ is mixed with
weights $p$ and $p+1$, since
$H^p(\overline{Y_{\underline{i}}},\mathbb Q_\ell)$ and
$H^{p-1}(\overline{Y_{\underline i}\cap H},\mathbb Q_\ell(-1))$
are pure of weights $p$ and $p+1$, respectively. Hence the
spectral sequence (2.3) is much simpler than (2.1) and has the
following $E_2$-layer:

\begin{center}
\myfig
{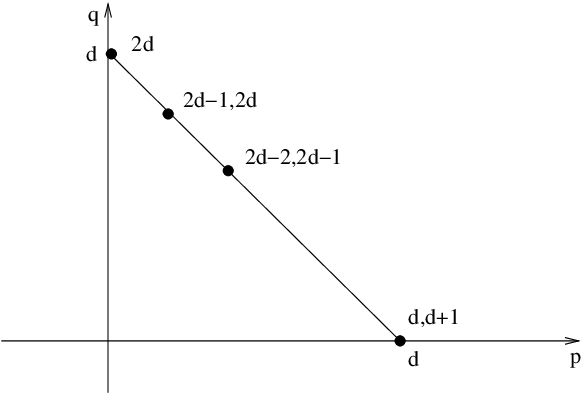}
{0.6}
%\bildunterschrift{}
\end{center}

%\begin{figure}[h]
%\begin{center}
%\scalebox{0.7}
%{\input{grafik2.eepic}}
%\end{center}
%\end{figure}

%\begin{center}
%\vskip1cm
%\myfig{grafik2.eps}{0.6}
%\end{center}

The terms vanish for $p+q>d$, and on the line $p+q=d$ the
$E_2^{p,q}$-terms - and hence also the $E^{p,q}_\infty$-terms
which are quotients - are divisible and mixed of the indicated
weights. Note that $H^0(\overline{Y^{0[d]}},\, \mathbb
Q_\ell/\mathbb Z_\ell(-d)) = H^0(\overline{Y^{[d]}},\qzl(-d)))$ is pure of weight $2d$.

\medskip

Let $F^\cdot$ be the descending filtration on $B_0=H^d(
\overline{U^0},\mathbb Q_\ell/\mathbb Z_\ell)$ for which
$F^\nu/F^{\nu-1}=E^{\nu,d-\nu}_\infty$. Then we see that $F^2$ is
divisible and mixed of weights $\leq 2d-1$. Next,
$$
F^1/F^2\cong E_2^{1,d-1}=H^1(\overline{{Y^0}{^{[d-1]}}}, \mathbb
Q_\ell/\mathbb Z_\ell(-d+1))
$$
is the cohomology of a (usually non-connected) smooth affine
curve, and by the Gysin sequence
$$
\begin{array}{ccl}
0 & \rightarrow  & H^1(\overline{Y^{[d-1]}},\mathbb Q_\ell/\mathbb
Z_\ell)\rightarrow H^1 (\overline{{Y^0}{^{[d-1]}}},\mathbb
Q_\ell/\mathbb Z_\ell)\\ & \rightarrow &
H^0(\overline{Y^{[d-1]}\cap H},\mathbb Q_\ell/\mathbb
Z_\ell(-1))\rightarrow H^2 (\overline{{Y}{^{[d-1]}}},\mathbb
Q_\ell/\mathbb Z_\ell) \rightarrow 0
\end{array}
$$
there is an exact sequence
$$
0\rightarrow A'\rightarrow F^1/F^2\rightarrow C'\rightarrow 0
$$
where $A'$ is divisible of weight $2d-1$ and where
$C'=I'\otimes\mathbb Q_\ell/\mathbb Z_\ell$, with $I'$ defined by
the exact sequence
$$
0\rightarrow I'\rightarrow\mathbb Z[\pi_0(\overline{Y^{[d-1]}\cap
H})] \rightarrow\mathbb Z[\pi_0(\overline{Y^{[d-1]}})]\rightarrow
0.
$$
Finally,
$$
C'':=F^0/F^1\cong H^0(\overline{{Y^0}{^{[d]}}}, \mathbb
Q_\ell/\mathbb Z_\ell (-d))
$$
is the cohomology of ${Y^0}{^{[d]}}=Y^{[d]}$ which is a union of
points, and $C''=I''\otimes \mathbb Q_\ell/\mathbb Z_\ell$ for
$$
I''=\mathbb Z[\pi_0(\overline{Y^{[d]}})].
$$
Let $A_0$ be the preimage of $A'$ in $F^1$, and $C_0=B_0/A_0$.
Then we have exact sequences
$$
\begin{array}{ccccccccc}
0&\rightarrow & A_0 &\rightarrow&B_0 & \rightarrow&C_0 & \rightarrow&0\\
0&\rightarrow & F^2 &\rightarrow&A_0 & \rightarrow&A'  & \rightarrow&0\\
0&\rightarrow & C'  &\rightarrow&C_0 & \rightarrow&C'' &
\rightarrow&0
\end{array}
$$
Hence $A_0$ is divisible and mixed of weights $\leq 2d-1$, and
$C_0$ is divisible of weight $2d$. This determines $A_0$ and $C_0$
uniquely (there is no non-trivial $G_K$-morphism between such
modules), and so the spectral sequence (2.1) for
$\mathop{U}\limits^{0}=X\smallsetminus\tilde Y$ instead of
$U=X\smallsetminus Y$ shows that $C_0$ is a quotient of
$$
\ker(H^0(\overline{\tilde Y^{[d]}},\mathbb Q_\ell (-d))\rightarrow
H^2(\overline{\tilde Y^{[d-1]}},\mathbb Q_\ell (-d+1))).
$$
Hence the action of $G_K$ on $C_0(d)$ factors through a finite
quotient $G$. This in turn shows that the extension
$$
0\rightarrow C'\rightarrow C_0\rightarrow C''\rightarrow 0
$$
comes from an extension
$$
0\rightarrow I'\rightarrow I\rightarrow I''\rightarrow 0
$$
of $G$-modules by tensoring with $\mathbb Q_\ell/\mathbb Z_\ell
(d)$. In fact, applying a Tate twist is an exact functor on
$\mathbb Z_\ell$-$G_k$-modules, and one has isomorphisms (where
the tensor products are over $\mathbb Z$)
$$
\begin{array}{rcl}
\mbox{Ext}^1_G(I'',I')\otimes_\mathbb Z\mathbb Z_\ell &
\stackrel{\sim}{\longrightarrow} & \mbox{Ext}^1_{\mathbb
Z_\ell[G]} (I''\otimes \mathbb Z_\ell, I'\otimes \mathbb Z_\ell)\\
& \stackrel{\sim}{\longrightarrow} & \mbox{Ext}^1_G (I''\otimes
\mathbb Q_\ell/\mathbb Z_\ell, I'\otimes \mathbb Q_\ell/\mathbb
Z_\ell),
\end{array}
$$
since $\mathbb Z_\ell$ is flat over $\mathbb Z$, and since the
functor $T\mapsto T\otimes_{\mathbb Z_\ell}\mathbb Q_\ell/\mathbb
Z_\ell$ is an equivalence between $\mathbb Z_\ell$-lattices and
divisible $\ell$-torsion modules of cofinite type (with action of
$G$) preserving exact sequences. Finally,
$$
\mbox{Ext}^1_G(I'',I')\rightarrow \mbox{Ext}^1_G(I'',
I')\raisebox{-8pt}{$\stackrel{\textstyle\otimes}{\scriptstyle
\mathbb Z}$}~\mathbb Z_\ell
$$
is surjective for a finite group $G$.

\medskip

We are now ready to prove

\vspace{0,5cm}

\noindent\textbf{Theorem 2.4} The restriction map
$$
\beta_B:H^2(K,B)\rightarrow\textstyle\bigoplus\limits_v H^2(K_v,B)
$$
is injective for $B=H^d(\overline{U^0}, \mathbb Q_\ell/\mathbb
Z_\ell (d+1)).$

\vspace{0,5cm}

\noindent\textbf{Proof} We follow the method of [Ja4]. By
applying the $(d+1)$-fold Tate twist to the sequence $0\rightarrow
A_0\rightarrow B_0\rightarrow C_0\rightarrow 0$ of Proposition
2.2, we get an exact sequence
$$
0\rightarrow A\rightarrow B \rightarrow C\rightarrow 0\, .
$$
It induces a commutative diagram with exact rows
$$
\xymatrix{ \ldots \ar[r] &
\raisebox{-6pt}{$\stackrel{\bigoplus}{\scriptstyle v\in S}$}
H^1(K_v, C)\ar[r]&
\raisebox{-6pt}{$\stackrel{\bigoplus}{\scriptstyle v}$}
H^2(K_v,A)\ar[r] &
\raisebox{-6pt}{$\stackrel{\bigoplus}{\scriptstyle v}$}
H^2(K_v, B)\\
\ldots\ar[r] & H^1(K,
C)\ar[r]\ar[u]_{\alpha_{C,S}\hspace{1cm}\displaystyle (\ast)} &
H^2(K,A)\ar[r]\ar[u]_{\beta_A}^\wr & H^2(K,B)\ar[u]_{\beta_B}}
$$

\medskip

$$
\xymatrix{ \ar[r] &
\raisebox{-6pt}{$\stackrel{\bigoplus}{\scriptstyle v}$} H^2(K_v,
C)\ar[r] & \raisebox{-6pt}{$\stackrel{\bigoplus}{\scriptstyle
v\infty}$}
H^3(K_v,A)\ar[r] & \ldots\\
\ar[r] & H^2(K, C)\ar[r]\ar[u]_{\beta_C}&
H^3(K,A)\ar[r]\ar[u]^\wr_{\gamma_A}&\ldots}
$$
for a suitable finite set $S$ of places of $K$. In fact, if
$S_{bad}$ is a set of bad places for $A$, then for any $S\supset
S_{bad}\cup\{v\mid\ell\},\; H^2(K_v,A)=0$ for $v\notin S $, and
thus $(\ast)$ is commutative. By Theorem 1.5, $\beta_A$ is an
isomorphism, since $A$ is divisible and mixed of weights $\leq
-3$, and by Tate duality, $\gamma_A$ is an isomorphism (for all
torsion modules $A$). To show the injectivity of $\beta_B$ by the
5-lemma, it therefore suffices to show that $C$ satisfies
$$
\begin{array}{ccrccl}
(H) & (i) & \alpha_{C,S}:H^1(K,C) & \rightarrow &
\textstyle\bigoplus\limits_{v\in S}
H^1(K_v,C)&  \mbox{ is surjective for all finite } S.\\
& (ii) & \beta_C:H^2(K,C) & \rightarrow &
\textstyle\bigoplus\limits_v H^2(K_v,C) &  \mbox{ is injective.}
\end{array}
$$
Let $I, I'$ and $I''$ be as in Proposition 2.2, so that $C=I
\otimes_\mathbb Z~\mathbb Q_\ell/\mathbb Z_\ell(1)$. We have exact
sequences
$$
\begin{array}{ccccccccl}
0 & \rightarrow & I' & \rightarrow & I & \rightarrow & I'' &
\rightarrow & 0\\
0 & \rightarrow & I' & \rightarrow & I_2 & \rightarrow & I_3 &
\rightarrow & 0\; ,
\end{array}
$$
in which $I'', \, I_2$ and $I_3$ are permutation modules, i.e., of
the form $\mathbb Z[M]$ for a $G_K$-set $M$. Thus $(H)$ holds for
$C$ by repeated application (first to $I'', I_2$ and $I_3$, then
to $I'$, and finally to $I$) of the following result

\vspace{0,5cm}

\noindent\textbf{Theorem 2.5} Let $I_1$, $I_2$ and $I_3$ be
finitely generated free $\mathbb Z$-modules with discrete
$G_K$-action, and let $C_i=I_i \otimes_\mathbb Z\mathbb
Q_\ell/\mathbb Z_\ell (1)$ for $i= 1,2,3$. Assume that $I_3$ is a
permutation module.

\begin{itemize}

\item[(a)] Property $(H)$ holds for $C_3$.

\item[(b)] If $0  \rightarrow  I_1  \rightarrow  I_2  \rightarrow
 I_3  \rightarrow  0$ is an exact sequence, then $(H)$ holds for
 $C_1$ if and only if it holds for $C_2$.

\end{itemize}

\medskip

The following observation will help to prove part (b).

\vspace{0,5cm}

\noindent\textbf{Lemma 2.6} Let $I$ be a finitely generated free
$\mathbb Z$-module with discrete $G_K$-action, and let $T$ be the
torus over $K$ with cocharacter module $X_\ast(T)=I$. Then
property $(H)(i)$ (resp. $(H)(ii))$ holds for $C=I \otimes_\mathbb
Z~\mathbb Q_\ell/\mathbb Z_\ell (1)$ if and only if $T$ satisfies
$$
\begin{array}{rlcll}
(H'_\ell) (i): & \!\alpha_{T,S,\ell}:H^1(K,T)\{\ell\} &
\rightarrow  &\!\textstyle\bigoplus\limits_{v\in S}
H^1(K_v,T)\{\ell\}   \mbox{ is surjective for all finite } S.\\
(\mbox{resp. }(H'_\ell) (ii):  & \!\beta_{T,\ell}:H^2(K,T)\{\ell\}
& \rightarrow & \textstyle\bigoplus\limits_v H^2(K_v,T)\{\ell\}
\mbox{ is injective.})
\end{array}
$$

\vspace{0,5cm}

\noindent\textbf{Proof} Recall that $T(\overline K) =I
\otimes_\mathbb Z\overline {K^\times}$ and $H^i(K,T)=
H^i(K,T(\overline K))$ by definition. Since $\ell\neq$ char$(K),
~T(\overline K)$ is $\ell$-divisible, and the Kummer sequences
$$
0\rightarrow I\otimes_{\mathbb Z}\mu_{\ell^n}\rightarrow T(\overline
K)\stackrel{\ell^n}{\longrightarrow} T(\overline K)\rightarrow 0
\leqno(2.7)
$$
identify $C$ with $T(\overline K)\{\ell\}$, the $\ell$-primary
torsion subgroup of $T(\overline K)$. Similar results hold for the
fields $K_v$, and the cohomology sequences associated to (2.7) for
all $n$ give rise to a commutative diagram with exact rows
$$
\xymatrix{\raisebox{-8pt}{$\stackrel{\bigoplus}{\scriptstyle v\in
S}$}~T(K_v)\otimes \mathbb Q_\ell/\mathbb Z_\ell\ar[r]&
\raisebox{-8pt}{$\stackrel{\bigoplus}{\scriptstyle v\in
S}$}~H^1(K_v,C)\ar[r]&\raisebox{-8pt}{$\stackrel{\bigoplus}
{\scriptstyle v\in S}$}~H^1(K_v,T)\{\ell\} \ar[r]&0\\
T(K)\otimes \mathbb Q_\ell/Z_\ell\ar[r]\ar[u]_{\omega_{T,S}} &
H^1(K,C)\ar[r]\ar[u]_{\alpha_{C,S}} &
H^1(K,T)\{\ell\}\ar[r]\ar[u]_{\alpha_{T,S,\ell}} & 0}
$$
and to a commutative diagram with horizontal isomorphisms
$$
\xymatrix{ \raisebox{-6pt}{$\stackrel{\bigoplus}{\scriptstyle
v}$}~H^2(K_v,C)\ar[r]^{\hspace{-0,5cm}\sim} &
\raisebox{-6pt}{$\stackrel{\bigoplus}{\scriptstyle v}$}~H^2(K_v,
T)\{\ell\}\\
H^2(K,C)\ar[r]^{\hspace{-0,5cm}\sim}\ar[u]_{\beta_C} &
H^2(K,T)\{\ell\}\; .\ar[u]_{\beta_{T,\ell}}}
$$
Here the vertical maps are induced by the various restriction
maps, and we used that $T(K)=H^0(K,T (\overline K))$ and $T(K_v)=
H^0(K_v,T(\overline K_v))$ for a separable closure $\overline
{K_v}$ of $K_v$. Note that $H^i(K,T)$ and $H^i(K_v,T)$ are torsion
groups for $i\geq 1$.

\medskip

Now the map $\omega_{T,S}$ is surjective for any torus $T$ and any
finite set of places $S$ ([Ja4] Lemma 2). This proves the lemma.

\vspace{0,5cm}

\noindent\textbf{Proof of theorem 2.5} Let $I_3$ be a permutation
module. Then $I_3$ is a direct sum of modules of the form
$I_0=\mbox{Ind}^K_{K'}(\mathbb Z)=\mathbb Z[G_K/G_{K'}]$ for some
finite separable extension $K'$ of $K$. Let $T_0$ be the torus
with cocharacter module $I_0$. Then
$$
H^i(K,T_0)\cong H^i(K',\mathbb G_m)=H^i(K',\overline
K^\times)\leqno{(2.8)}
$$
by Shapiro's lemma, and similarly
$$
H^i(K_v, T_0)\cong
\raisebox{-6pt}{$\stackrel{\bigoplus}{\scriptstyle w\mid
v}$}H^i(K'_\omega,\mathbb G_m)\, ;\leqno{(2.9)}
$$
where $w$ runs through the places of $K'$ above $v$. Thus
$$
H^1(K,T_0)=0=H^1(K_v,T_0)
$$
by Hilbert's theorem 90, and
$\beta_{T_0}: H^2(K,T_0) \rightarrow\raisebox{-6pt}{$\stackrel{\bigoplus}{\scriptstyle
v}$}H^2(K_v,T_0)$ is injective by the classical theorem of
Brauer-Hasse-Noether for $K'$. This shows property $(H'_\ell)$ for
the torus $T_3$ with cocharacter module $I_3$, for all primes
$\ell$, and hence part (a) of theorem 2.5.

\medskip

For part (b), let $T_i$ be the torus with cocharacter module
$I_i~(i=1,2,3)$. Then we have an exact sequence
$$
0\rightarrow T_1 \rightarrow T_2 \rightarrow T_3 \rightarrow 0
$$
with $H^1(K,T_3)=0=H^1(K_v,T_3)$ by assumption and the above. This
gives exact commutative diagrams
$$
\xymatrix{ \raisebox{-8pt}{$\stackrel{\bigoplus}{\scriptstyle v\in
S}$}T_3(K_v)\ar[r]^{\hspace{-0,5cm}\delta} &
\raisebox{-8pt}{$\stackrel{\bigoplus}{\scriptstyle v\in
S}$}H^1(K_v,T_1)\ar[r] &
\raisebox{-8pt}{$\stackrel{\bigoplus}{\scriptstyle v\in
S}$}H^1(K_v,T_2)\ar[r] & 0\\
T_3(K)\ar[r]\ar[u]_\omega &
H^1(K,T_1)\ar[r]\ar[u]_{\alpha_{T_1,S}} &
H^1(K,T_2)\ar[r]\ar[u]_{\alpha_{T_2,S}} & 0}
$$
and
$$
\xymatrix{0\ar[r] &
\raisebox{-6pt}{$\stackrel{\bigoplus}{\scriptstyle
v}$}H^2(K_v,T_1)\ar[r] &
\raisebox{-6pt}{$\stackrel{\bigoplus}{\scriptstyle
v}$}H^2(K_v,T_2)\ar[r] &
\raisebox{-6pt}{$\stackrel{\bigoplus}{\scriptstyle
v}$}H^2(K_v,T_3)\\
0\ar[r] & H^2(K,T_1)\ar[r]\ar[u]_{\beta_{T_1}} &
H^2(K_v,T_2)\ar[r]\ar[u]_{\beta_{T_2}} &
H^2(K_v,T_3)\ar[u]_{\beta_{T_3}}& .}
$$
Since $\beta_{T_3}$ is injective by assumption, one has an
isomorphism
$\ker\beta_{T_1}\stackrel{\sim}{\longrightarrow}\ker\beta_{T_2}$.
On the other hand, the groups $H^1(K_v,T_1)$ have finite exponent
$n$ (by Hilbert's theorem 90 we can take $n=[K':K]$, if $K'/K$ is
a finite Galois extension splitting $T_1$). Hence $\delta$ factors
through $\raisebox{-8pt}{$\stackrel{\bigoplus}{\scriptstyle v\in
S}$}T_3(K_v)/n$. But $\omega\otimes \zn$ is
surjective for every $n$: Indeed,
$$
K^\times/(K^\times)^n\rightarrow
\raisebox{-8pt}{$\stackrel{\bigoplus}{\scriptstyle v\in
S}$}K^\times_v/(K^\times_v)^n
$$
is surjective for all $n$ by weak approximation for $K$, and the
same for all finite extensions $K'$ of $K$ gives the result for
$T_3$ (cf. (2.8) and (2.9) for $i=0$). This gives an isomorphism
coker $\alpha_{T_1,S}\stackrel{\sim}{\longrightarrow}$ coker
$\alpha_{T_2,S}$ and hence (b).

\medskip

This completes the proof of Theorem 2.4, and we can now show the following
theorem which is a variant of Theorems 0.2, in which the fields $F_v$ are replaced by the fields $F_{(v)}$.

\vspace{0,5cm}

\noindent\textbf{Theorem 2.10} Let $F$ be a function field in $d$
variables over $K$, such that $K$ is separably closed in $F$,
and let $\ell$ be a prime invertible in $K$. Then the restriction map
$$
 H^{d+2}(F,\mathbb Q_\ell/\mathbb Z_\ell(d+1))\longrightarrow
\textstyle\bigoplus\limits_v H^{d+2}(F_{(v)},\mathbb Q_\ell/\mathbb Z_\ell(d+1))
$$
is injective.

\vspace{0,5cm}

\noi\textbf{Proof} Let $F = K(V)$ for a geometrically integral variety $V$ of dimension
$d$ over $K$.
By Proposition 1.2 it is equivalent to show the injectivity of
$$
\alpha: \; H^2(K,H^d(F\overline{K},\qzl(d+1))) \longrightarrow \textstyle\bigoplus\limits_v
H^2(K_{(v)},H^d(F\overline{K},\qzl(d+1)))\,.
$$
Let $x$ be an element in the kernel of $\alpha$. By the limit property recalled above
Theorem 1.5, there is an open affine $U\subset V$ such that $x$ is the image of an element
$y$ lying in the kernel of
$$
H^2(K,H^d(\overline{U},\qzl(d+1))) \longrightarrow \textstyle\bigoplus\limits_v H^2(K_{(v)},H^d(\overline{U},\qzl(d+1)))\,.
$$
If $K$ is a number field, then we may assume that $U$ is smooth over $K$, and there is a good compactification
$U \subset X$ as in property RS2$(U)$ at the beginning of this section. Thus the claim
follows immediately by restricting to the subset $U^0$ constructed before Proposition 2.2,
and applying Theorem 2.4.

If $K$ has positive characteristic, we use the following result of Gabber (see the notes \cite{Il2} by Illusie,
Theorem 1.3), which refines the theorem of de Jong on alterations.

\v\noindent\textbf{Theorem 2.11}(Gabber, see \cite{Il2} and \cite{ILO})
If $X$ is separated and integral of finite type over a field $L$ and $\ell$ is a prime which is invertible in $L$,
and $Y \subset X$ is a proper closed subscheme, then there exists a finite extension $L'/L$ of degree prime to $\ell$
and a connected smooth quasi-projective variety $X'$ over $L'$ together with a proper surjective $L$-morphism
$\pi: X' \rightarrow X$ such that the extension of function fields $L'(X')/L(X)$ is finite of degree prime to $\ell$,
and such that $Y' = \pi^{-1}(Y)$, with the reduced subscheme structure, is a divisor with strict normal crossings on $X'$.

\v
We apply this to a compactification $U \subset X$ for our affine variety with a proper integral variety $X$ over $K$,
and the closed subset $Y = X - U$. Let $\pi: X' \rightarrow X$ and $Y' = \pi^{-1}(Y)$ be as in $(G)$, so that
$X'$ is smooth projective, without loss of generality geometrically irreducible, and $Y'$ is a simple normal crossings divisor.
Let $U' = X'-Y'$, and let $(U')^0 \subset U'$ be constructed as the complement of a well-chosen hyperplane section like before
Proposition 2.2. Then the image $y'$ of $y$ under the restriction map for $U' \rightarrow U$ lies in the kernel of
$$
H^2(K',H^d(\overline{(U')^0},\qzl(d+1))) \longrightarrow \textstyle\bigoplus\limits_w H^2(K'_{(w)},H^d(\overline{(U')^0},\qzl(d+1)))\,,
$$
where $w$ runs over all places of $K'$. Thus $y'=0$ by Theorem 2.4. By restricting to $F'\overline{K'}$ and applying
once more Proposition 1.2 this implies that the image of $x$ under
$$
H^{d+2}(F,\mathbb Q_\ell/\mathbb Z_\ell(d+1))\longrightarrow H^{d+2}(F',\mathbb Q_\ell/\mathbb Z_\ell(d+1))
$$
is zero. It remains to remark that this restriction map is injective, because the degree $[F':F]$ is prime
to $\ell$. In fact, we can decompose the extension $F'/F$ as $F'/F_i/F$, where $F_i$ is the maximal
inseparable extension inside $F'/F$. Then the restriction from $F$ to $F_i$ is an isomorphism, and the
restriction $Res$ from $F_i$ to $F'$ is injective, since, for the corestriction $Cor$ from $F'$ to $F_i$
we have $Cor Res =$ multiplication by $[F':F_i]$, which is prime to $\ell$.

\bigskip
To have the same result with $F_v$ in place of $F_{(v)}$, and thus obtain Theorem 0.2, it
suffices to show:

\vspace{0,5cm}

\noi \textbf{Theorem 2.12} For any $n\in\mathbb N$ and all $i,j\in\mathbb Z$,
the restriction map
$$
H^i(F_{(v)},\zn(j))\rightarrow H^i(F_v,\zn(j))
$$
is injective.

\v
\noi This is related to a more precise rigidity result (for $n$ invertible in $K$)
on the Kato complexes recalled in Theorem 4.11, which we shall also need in the
following sections. However, as was pointed out to me by J.-L.
Colliot-Th\'{e}l\`{e}ne, the injectivity above follows by a
simple argument, and in the following general version:

\vspace{0,5cm} \noi \textbf{Theorem 2.13} Let $K/k$ be a field
extension satisfying  the following property:

\medskip
\noi {\bf (SD)} If a variety $Y$ over $k$ has a
$K$-rational point, then it also has a $k$-rational point.

\medskip\noi Let $F$ be a set-valued contravariant functor
on the category of all $k$-schemes such that

\medskip\noi {\bf (FP)} For any inductive system $(A_i)$ of $k$-algebras and $A=
\raisebox{-6pt}{$\stackrel{\textstyle\lim}{\scriptstyle\longrightarrow}$}\phantom{}_i
A_i$, the natural map
$\raisebox{-6pt}{$\stackrel{\textstyle\lim}{\scriptstyle\longrightarrow}$}\phantom{}_i
\, F(A_i)\cong F(A)$ is an isomorphism. (Here we write $F(B):=F(\,\mbox{Spec } B)$ for a $k$-algebra
$B$.)

\medskip\noi Let $V$ be a geometrically integral variety
over $k$, and write $k(V)$ (resp. $K(V))$ for the function field
of $V$ (resp. $V\times_k K)$. Then the map
$$
F(k(V))\rightarrow F(K(V))
$$
is injective.

\v
\noi \textbf{Proof} The field $K$ can be written as
the union of its subfields $K_i$ which are finitely generated (as
fields) over $k$. Every $K_i$ can of course be written as the
fraction field of a finitely generated $k$-algebra $A_i$.

Now let $\alpha \in F(k(V))$, and assume that $\alpha$ vanishes in
$F(K(V))$. By (FP), there is an $i$ such that $\alpha$ already
vanishes in $F(K_i(V))$. Moreover, there is a non-empty affine
open $V'\subseteq V$ and a $\beta \in F(V')$ mapping to $\alpha$
in $F(k(V))$. Finally there is a non-empty affine open $U\subseteq
Z_i\times_k V'$, where $Z_i=\,\mbox{Spec } A_i$, such that $\beta$
vanishes under the composite map $F(V')\rightarrow
F(Z_i\times_kV')\rightarrow F(U)$.

Now it follows from Chevalley's theorem that the image of $U$
under the projection $p:Z_i\times_kV'\rightarrow Z_i$ contains a
non-empty affine open $U'$ ($p$ maps constructible set to
constructible sets, and is dominant). Now $U'$ has a $K$-point
Spec$(K)\rightarrow$ Spec$(K_i)\hookrightarrow U'$. Hence, by
property (SD), $U'$ has a $k$-rational point $Q$. Then $W =
p^{-1}(Q)\cap U$ is open and non-empty in
$p^{-1}(Q)=Q\times_kV'\cong V'$. By functoriality, $\beta$
vanishes in $F(W)$, and thus $\alpha$ in $F(k(V))$.

\vspace{0,5cm} \noi \textbf{Proof of Theorem 2.12} We may apply 2.13
to the extension $K_v/K_{(v)}$ and the functor $F(X)=
H^i_{\mbox{\scriptsize\'{e}t}}(X,M)$ for any fixed discrete
$G_{K_{(v)}}$-module $M$ (regarded as \'{e}tale sheaf by
pull-back) to get the injectivity of
$$
H^i(F_{(v)},M)\rightarrow H^i(F_v, M)\; .
$$
In fact, property (SD) (for ``strongly dense'') is known to hold
in this case (cf. \cite{Gre} Theorem 1), and the commuting with limits as
in (FP) (for ``finitely presented'') is a standard property of
\'{e}tale cohomology (cf. [Mi2] III 1.16).

\vspace{0,5cm}

\vspace{15mm}

\centerline{\noindent\textbf{\S 3 A crucial exact sequence, and a Hasse principle for unramified
cohomology}}

\vspace{1,0cm}

To investigate the cokernel of $\beta_B$ (notations as in section
2), we could follow the method of [Ja4] and show that it is
isomorphic to coker$(\beta_C)$. By describing the edge morphisms
in the spectral sequence (2.3) we could prove the crucial Theorem 3.1 below
for global fields. Instead, we prefer to argue more directly, which
allows to treat arbitrary finitely generated fields and use 3.1 also
for the remaining sections.

\vspace{0,5cm}

We shall make repeated use of the following. Let
$i:Y\hookrightarrow X$ be a closed immersion of smooth varieties
over a field $L$, of pure codimension $c$. Then, for every integer
$n$ invertible in $L$ and every integer $r$, one has a long exact
Gysin sequence
$$
\ldots \rightarrow H^{\nu-1}(U,\zn(r))\stackrel{\delta}{\rightarrow}H^{\nu-2c}(Y,\zn(r-c))\stackrel{i_\ast}{\rightarrow}H^\nu(X,\zn(r))\stackrel{j^\ast}{\rightarrow} H^\nu(U,\zn(r))\rightarrow\ldots
$$
where $U=X\smallsetminus Y$ is the open complement of $Y$ and
$j:U\hookrightarrow X$ is the open immersion. We call $i_\ast$ and
$\delta$ the Gysin map and the residue map for $i:Y\hookrightarrow
X$, respectively. If $i':Y'\hookrightarrow Y$ is another closed
immersion, with $Y'$ smooth and of pure codimension $c'$ in $Y$,
then the diagram of Gysin sequences
$$
\xymatrix{ H^{\nu-1}(U,\zn(r))\ar[r]^{\hspace{-0,3cm}\delta} & H^{\nu-2c}(Y,\zn(r-c))\ar[r]^{\hspace{0,6 cm}i_\ast}& H^\nu(X,\zn(r))\ar[r]^{\hspace{1,3cm}j^\ast}& \\
H^{\nu-1}(U',\zn(r))\ar[r]^{\hspace{-1cm}\delta}\ar[u]_{j'^\ast} &
H^{\nu-2(c+c')}(Y',\zn(r-c-c'))\ar[r]^{\hspace{1,3cm}(i\circ
i')_\ast}\ar[u]_{i'_\ast} & H^\nu(X,\zn(r))\ar[r]^{\hspace{1,5cm}(j\circ j')^\ast}\ar@{=}[u] & }
$$
is commutative, where $j':U\hookrightarrow U'$ is the open immersion.
In fact, the first sequence comes from the long exact relative sequence involving $H^\ast_Y(X,\z/n\z(r)))$,
together with canonical Gysin isomorphisms
$$
H^{\nu-2c}(Y,\z/n\z(r-c)) \mathop{\rightarrow}\limits^{\sim} H^\nu_Y(X,\z/n\z(r)))\,.
$$

If $L$ is a perfect field of characteristic $p>0$ and $n=p^m$, then the one has still Gysin morphisms $i_\ast$ with
the transitivity property, by work of Gros \cite{Gro1}, but the remaining properties are not in general true
anymore, except for the following special case: If $X$ is smooth of pure dimension $d$,
then one has canonical Gysin isomorphisms
$$
H^{\nu-2c}(Y,\z/p^m\z(d-c)) \mathop{\rightarrow}\limits^{\sim} H^\nu_Y(X,\z/p^m\z(d)))\,,
$$
see \cite{Su} Cor. 2.6., and gets an exact Gysin sequence as above for $r=d$.

\medskip\noindent
With these preparations we can now prove a crucial exact sequence for a specialization map
which is not only used for Theorem 3.8 below, giving a Hasse principle for unramified cohomology,
but is also essential in the proofs of the Theorems 0.5, 0.6 and 0.7

\medskip

\noindent\textbf{Theorem 3.1}  Let $K$ be a finitely generated field with algebraic closure $\overline{K}$,
and let $X$ be a smooth, proper, irreducible variety of dimension $d$ over $K$.
Let $Y=\bigcup^r_{i=1} Y_i$, with $r\geq 1$, be a union of smooth
irreducible divisors on $X$ intersecting transversally such that
$X \smallsetminus Y_1$ is affine (this holds, e.g., if $X$ is projective and $Y_1$ is a
smooth hyperplane section), and let $U=X \smallsetminus Y$. Then,
for any prime $\ell$, and with the notations of
the beginning of section 2, the sequence
$$
0\rightarrow H^d(\overline U,\mathbb Q_\ell/\mathbb Z_\ell
(d))_{G_K} \stackrel{e}{\rightarrow} H^0(\overline{Y^{[d]}},
\mathbb Q_\ell/\mathbb Z_\ell)_{G_K} \stackrel{d_2}{\rightarrow}
H^2(\overline{Y^{[d-1]}},\mathbb Q_\ell/\mathbb Z_\ell(1) )_{G_K}
$$
is exact, where we write $\overline{X}=X\times_K\overline{K}$, and similarly for the other varieties,
and where we regard $G_K$ as $Gal(\overline{K}/K^{per})$ for the perfect hull $K^{per}$
of $K$ in $\overline{K}$, which is the maximal inseparable extension of $K$ inside $\overline{K}$ and
is a perfect field. Moreover $e$ and $d_2$ are defined as follows. The
specialization map $e$ is induced by the compositions
$$
H^d(\overline U,\mathbb Q_\ell/\mathbb Z_\ell
(d))\stackrel{\delta}{\rightarrow} H^{d-1}(\overline{Y_{i_d}
\smallsetminus (\textstyle\bigcup\limits_{i\not= i_d}
Y_i)},\mathbb Q_\ell/\mathbb Z_\ell (d-1)) \leqno(3.2)
$$
$$
\stackrel{\delta}{\rightarrow} H^{d-2}(\overline{ Y_{i_{d-1},i_d}
\smallsetminus (\textstyle\bigcup\limits_{i\not= i_{d-1},i_d}
Y_i)}, \mathbb Q_\ell/\mathbb Z_\ell (d-2))\longrightarrow \ldots
$$
$$
\ldots\stackrel{\delta}{\rightarrow}H^1(\overline
{Y_{i_2,\ldots,i_d} \smallsetminus
(\textstyle\bigcup\limits_{i\not= i_2,\ldots,i_d} Y_i)}, \mathbb
Q_\ell/\mathbb Z_\ell (1)) \stackrel{\delta}{\rightarrow}
H^0(\overline{Y_{i_1,\ldots, i_d}},\mathbb Q_\ell/\mathbb Z_\ell),
$$
where each $\delta$ is the connecting morphism in the obvious
Gysin sequence. On the other hand
$d_2=\sum\limits^d_{\mu=1}(-1)^\mu\,\delta_\mu$, where
$\delta_\mu$ is induced by the Gysin map associated to the
inclusions
$$
Y_{i_1,\ldots,i_d}\hookrightarrow Y_{i_1,\ldots,\hat
i_\nu,\ldots,i_d}
$$
(and $\hat i_\nu$ means omission of $i_\nu$, as usual).

\vspace{0,5cm}

\noindent\textbf{Proof} We note that here the absolute Galois group $G_K$ of $K$ can be regarded
as the Galois group $Gal(\overline{K}/K^{per})$, where $K^{per}\subset \overline{K}$ is the
perfect hull of $K$ (the maximal inseparable extension of $K$ in $\overline{K}$).
For $\ell$ invertible in $K$, we could replace the algebraic closure of $K$ by its separable closure,
and, by a standard property of \'etale cohomology, we get isomorphic groups above,
which are the ones used in section 2. For $\ell =\ch(K)$
however, we need $\overline{K}$ to be the algebraic closure.

Write $H^i(\overline Z,j)$ instead of $H^i(\overline Z,\mathbb Q_\ell/\mathbb Z_\ell(j))$, for short, and note
that $U$ is affine, because $X\setminus Y_1$ is affine and $U \hookrightarrow X$
is an affine morphism, because $Y$ is defined by a locally principal ideal.
Hence $H^d(\overline{U},d)$ is divisible, since $H^{d+1}(\overline{U},\z/\ell\z(d))=0$
by weak Lefschetz, which also holds for $\ell = \ch(K)$, see \cite{Su} Lemma 2.1.
We now proceed by induction on $r$, the number of components of $Y$. If
$r=1$, then the Gysin sequence
$$
\ldots \rightarrow H^d(\overline X,d)\rightarrow H^d(\overline U,
d)\rightarrow H^{d-1} (\overline Y,d-1)\rightarrow\ldots
$$
shows that $H^d(\overline U,d)$ is mixed with weights  $-d$ and
$-d+1$, see 1.7 (e) to (f). Hence, using 1.7 (d), we can only have $H^d(\overline U,d)_{G_K}\not= 0$
and $Y^{[d]}\not= \emptyset$ for $d=1$. In this case we have an
exact sequence
$$
0\rightarrow H^1(\overline X,1)\rightarrow H^1(\overline U,1)
\stackrel{\delta}{\rightarrow} H^0(\overline Y,0)\rightarrow
H^2(\overline X,1)\rightarrow 0
$$
Without loss of generality, we may assume that $X$ is
geometrically irreducible over $K$ (otherwise this is the case
over a finite extension $K'$ of $K$, and everything reduces to
this situation, since we have induced modules). Letting
$C=\mbox{im}(\delta)$, we have
$$
H^1(\overline U,1)_{G_K} \stackrel{\sim }{\rightarrow}C_{G_K}
$$
since $H^1(\overline X, 1)_{G_K}=0~(H^1(\overline X,1)$ is
divisible and of weight $-1$), and there is an exact sequence
$$
0\rightarrow C\rightarrow \mbox{Ind}^K_{K(x)}(\mathbb
Q_\ell/\mathbb Z_\ell) \rightarrow \mathbb Q_\ell/\mathbb
Z_\ell\rightarrow 0
$$
where $K(x)$ is the residue field of the unique point $x\in Y$,
which is a separable extension of $K$, by assumption.
But this sequence stays exact after taking cofixed modules: the
action of $G_K$ factors through a finite quotient $G$, and
$H_1(G,\mathbb Q_\ell/\mathbb Z_\ell)=0$ (this group is dual to
$H^1(G,\mathbb Z_\ell)=0)$. Putting things together, we have an
exact sequence
$$
0\rightarrow H^1(\overline U,1)_{G_K} \stackrel{e}{\rightarrow}
H^0(\overline Y,0)_{G_K}  \stackrel{d_2}{\rightarrow}H^2(\overline
X,1)_{G_K}\rightarrow 0\,.
$$

\medskip

Now let $r>1$. Then $Z=\bigcup_{i=1}^{r-1} Y_i $ is a divisor
with normal crossings on $X$ that fulfills all the assumptions of
the theorem, and the same is true for $Z_r = Y_r\cap
Z=\bigcup_{i=1}^{r-1} (Y_r\cap Y_i)$ on $Y_r$.

\medskip

We claim that we obtain a commutative diagram

\smallskip

$$
\xymatrix{ & 0\ar[d] & 0\ar[d] & 0\ar[d] & \\
& H^d(\overline{X\smallsetminus Z},
d)_{G_K}\ar[r]\ar[d]^{e\hspace{1,5cm} \textstyle(1)} &
H^d(\overline{U},
d)_{G_K}\ar[r]^{\hspace{-1,5cm}\delta}\ar[d]^{e\textstyle\hspace{2cm}
(2)} & H^{d-1}(Y_r\smallsetminus (Y_r\cap
Z),d-1)_{G_K}\ar[r]\ar[d]^e &
0\\
0\ar[r] & H^0(\overline{Z^{[d]}},
0)_{G_K}\ar[r]\ar[d]^{d_2\hspace{1,4cm} \textstyle(3)}  &
H^0(\overline{Y^{[d]}},
0)_{G_K}\ar[r]\ar[d]^{d_2\textstyle\hspace{1,9cm} (4)} &
H^0((\overline{Y_r\cap Z)^{[d-1]}}, 0)_{G_K}\ar[r]\ar[d]^{d_2} & 0
\\
0\ar[r] & H^2(\overline{Z^{[d-1]}},1)_{G_K}\ar[r] &
H^2(\overline{Y^{[d-1]}},1)_{G_K}\ar[r] & H^2((\overline{Y_r\cap
Z)^{[d-2]}},1)_{G_K}\ar[r] & 0 }\leqno{(3.3)}
$$
with exact rows: The first row comes from the Gysin sequence for
$(X \smallsetminus Z, Y_r \smallsetminus Z_r)$
$$
\ldots\rightarrow H^d(\overline{X \smallsetminus Z},d)\rightarrow
H^d(\overline U,d)  \stackrel{\delta}{\rightarrow}
H^{d-1}(\overline{Y_r \smallsetminus Z_r},d-1)\rightarrow 0
$$
in which $H^{d+1}(\overline {X \smallsetminus Z},d)=0$ by weak
Lefschetz. Next note that
$$
\begin{array}{rccccl}
&& Y^{[\nu]}& = & \coprod\limits_{1\leq i_1<\ldots<i_\nu\leq r} &
Y_{i_1,\ldots,i_\nu}
\\\\
&&  Z^{[\nu]} & = & \coprod\limits_{1\leq i_1<\ldots<i_\nu\leq r-1} &
Y_{i_1,\ldots,i_\nu} \\\\
Y_r\cap Z^{[\nu-1]} & = & (Y_r\cap Z)^{[\nu-1]} & = & \coprod\limits_{1\leq i_1<\ldots <i_{\nu
-1}\leq r-1} & Y_r\cap Y_{i_1,\ldots,i_{\nu -1}} \\\\
&&& = &
\coprod\limits_{1\leq i_1< \ldots <i_\nu = r} &
Y_{i_1,\ldots,i_\nu}
\end{array}
$$
so that $Y^{[\nu]} = Z^{[\nu]}\;\coprod\; (Y_r \cap Z)^{[\nu-1]}$. Hence one has commutative diagrams
\smallskip
$$
\xymatrix{ 0 \ar[r]& H^i(\overline{Z^{[d]}},j)
\ar[r]\ar[d]^{d_2\hspace{1,5cm}\textstyle(3')} &
H^i(\overline{Y^{[d]}},j)
\ar[r]\ar[d]^{d_2\hspace{1,5cm}\textstyle(4')} &
H^i((\overline{Y_r\cap
Z)^{[d-1]}},j)\ar[r] \ar[d]^{d_2}& 0\\
0\ar[r]& H^{i+2}(\overline{Z^{[d-1]}},j+1)\ar[r] &
H^{i+2}(\overline{Y^{[d-1]}},j+1)\ar[r] &
H^{i+2}(\overline{(Y_r\cap Z)^{[d-2]}},j+1)\ar[r] & 0}\leqno(3.4)
$$
with canonically split exact rows, where the two left maps $d_2$ are $d_2=
\sum\limits^\nu_{\mu=1}(-1)^\mu \delta_\mu,$ with $\delta_\mu$
being induced by the inclusions
$$
Y_{i_1,\ldots,i_\nu}\hookrightarrow Y_{i_1,\ldots,\widehat{i_\mu},
\ldots,i_\nu}
$$
and where the right hand $d_2$ is similarly defined as $d_2=
\sum\limits^{\nu-1}_{\mu=1}(-1)^\mu \delta_\mu,$ with $\delta_\mu$
being induced by the inclusions
$$
Y_{i_1,\ldots,i_{\nu-1}}\cap Y_r\hookrightarrow
Y_{i_1,\ldots,\widehat{i_\mu},\ldots,i_{\nu-1}}\cap Y_r\, .
$$
In fact, the commuting of (3') is trivial, and the square (4')
commutes since it commutes with $\delta_\mu, 1\leq\mu\leq\nu -1$
in place of $d_2$, whereas $\delta_\nu$ vanishes after projection
onto $(Y_r\cap Z)^{(\nu -1)}$ (the last component of
$(i_1\ldots,\hat i_\nu)$ cannot be $r$). This implies the
commutativity of (3) and (4), and the exactness of the two
involved rows.

\medskip

The commutativity of (2) is clear: For $1\leq i_1<\ldots<i_d=r$
the specialization map (3.2) is the composition
$$
H^d(\overline U,d) \stackrel{\delta}{\rightarrow}
H^{d-1}(\overline{Y_r\smallsetminus Z},d-1)
\stackrel{\delta}{\rightarrow} H^{d-2}
(\overline{Y_{i_{\nu-1}}\cap
Y_r\smallsetminus(\raisebox{-8pt}{$\stackrel{\bigcup}{\scriptstyle
i\not= i_{\nu-1,r}}$} (Y_i\cap Y_r})),d-2)
$$
$$
\ldots \stackrel{\delta}{\rightarrow}
H^1(\overline{Y_{i_2},\ldots,r
\smallsetminus(\raisebox{-8pt}{$\stackrel{\bigcup}{\scriptstyle
i\not= i_2,\ldots,r}$}(Y_i\cap Y_r)},1)
\stackrel{\delta}{\rightarrow} H^0(\overline{Y_{i_1,\ldots,r}},0).
$$
The commutativity of (1) is implied by the commutativity of
$$
\xymatrix{ e: & \!\!\!\!H^d(\overline U,d)
\ar[r]^{\hspace{-1,3cm}\delta} & H^{d-1}(
\overline{Y_{i_d}\smallsetminus
\raisebox{-8pt}{$\stackrel{\displaystyle\cup}{\scriptstyle i\not=
i_d}$} Y_i}, d-1)\ar[r] & \ldots & \ldots  \ar[r] &
H^0(\overline{Y_{i_1,\ldots,i_d}},0)\\
e: &
\!\!\!\!H^d(\overline{X-Z},d)\ar[r]^{\hspace{-1,3cm}\delta}\ar[u]
& H^{d-1}( \overline{Y_{i_d}\smallsetminus
\raisebox{-8pt}{$\stackrel{\displaystyle\cup}{\scriptstyle i\not=
i_{d,r}}$}Y_i}, d-1)\ar[r]\ar[u] & \ldots & \ar[u] \ldots
\ar[r]& H^0(\overline{Y_{i_1,\ldots,i_d}},0)\ar@{=}[u] }
$$
for $1\leq i_1 <\ldots< i_d <r$, where the vertical maps are the
restriction maps for the open immersions obtained by deleting
$Y_r$ everywhere (note that $Y_{i_1,\ldots,i_d}\cap Y_r=\emptyset$
for $i_d\not= r$). This commutativity follows from the
compatibility of the corresponding Gysin sequences with
restriction to open subschemes.

\medskip

Given the diagram (3.3), we can carry out the induction step: It is easy
to check that the middle column is a complex, and by
induction the left and right column are exact. Hence the middle column is
exact, by a straightforward diagram chase.

\medskip

We give a first application to function fields. Recall the following definition
[CT2] 2.1.8 and 4.1.1.

\vspace{0,5cm}

\noindent\textbf{Definition 3.5} Let $k$ be a field and let $F$ be
a function field over $k$. For an integer $n$ invertible in $k$,
the unramified cohomology
$H^i_{nr}(F/k,\zn(j))\subseteq H^i(F,\zn(j))$ is defined as the subset of elements lying in
the image of
$$
H^i_{\mbox{\scriptsize \'{e}t}}(\mbox{Spec}\,A,\zn(j))\rightarrow H^i(F,\zn(j))
$$
for all discrete valuation rings $A\subseteq F$ containing $k$.

\vspace{0,5cm}

If $\lambda$ is a discrete valuation of $F$ which is trivial on
$k$, and if $A_\lambda$ and $k(\lambda)$ are the associated
valuation ring and residue field, respectively, then one has an
exact Gysin sequence
$$
\ldots H^i_\et(\Spec A_\lambda,\zn(j)) \rightarrow H^i(F,\zn(j))
\mathop{\rightarrow}\limits^{\delta_\lambda}
H^{i-1}(k(\lambda),\zn(j-1)) \rightarrow \ldots ...\,,
$$
since purity is known to hold in this situation. We call the map
$\delta_\lambda$ the residue map for $\lambda$. This shows:

\vspace{0,5cm}

\noindent\textbf{Lemma 3.6} One has
$$
H^i_{nr}(F/k,\zn(j)) = \ker(H^i(F,\zn(j))
\rightarrow \prod_\lambda H^{i-1}(k(\lambda),\zn(j-1)))\,,
$$
where the sum is over all discrete valuations $\lambda$ of $F/k$,
and the components of the map are the residue maps
$\delta_\lambda$.

\vspace{0,5cm}

We will need the following fact (cf. [CT2] 4.1.1).

\vspace{0,5cm}

\noindent\textbf{Proposition 3.7} Let $X$ be a smooth proper
variety over $k$, and let $F=k(X)$ be its function field. Then
$$
H^i_{nr}(F/k,\z/n\z(j))=\ker(H^i(F,\z/n\z(j))\stackrel{\delta_X}{\rightarrow}\raisebox{-12pt}{$\stackrel{\textstyle\bigoplus}{\scriptstyle
x\in X^1}$}\,H^{i-1}(k(x),\z/n\z(j-1)))
$$
where $X^i=\{x\in X\!\mid\!\dim\mathcal O_{X,x}=i\}$ for
$i\geq 0,\, k(x)$ is the residue field of $x\in X$, and $\delta$ is
the map from the Bloch-Ogus complexes for \'{e}tale cohomology
[BO]. In particular,
$$
H^i_{nr}(F/k,\zn(j))\cong H^0_{\mbox{\scriptsize
Zar}}(X,\mathcal H^i_n(j))
$$
where $\mathcal H^i_n(j)$ is the Zariski sheaf on $X$ associated
to the presheaf $U\mapsto H^i_{\mbox{\scriptsize
\'{e}t}}(U,\zn(j))$.

\vspace{0,5cm}\noi
\textbf{Proof} Since we need a variant below, we recall the beautiful argument.
First note that, by definition of the Bloch-Ogus sequence, the
components of $\delta_X$ are the residue maps $\delta_{X,x} :=
\delta_{\lambda(x)}$, where $\lambda(x)$ is the discrete valuation
associated to $x$ (so that $A_{\lambda(x)} = \mathcal O_{X,x}$ and
$k(\lambda(x)) = k(x)\,)$. This shows that the kernel of 3.6 is contained in the
kernel of 3.7. Conversely, let $A\subset F$ be a discrete valuation ring. Then
by properness of $X$ we have a factorization $\Spec(F) \rightarrow \Spec(A) \rightarrow X$,
and hence a factorization
$$
\Spec(F) \rightarrow \Spec(A) \rightarrow \Spec(\mathcal O_{X,x}) \rightarrow X\,,
$$
where $x\in X$ is the image of the closed point of $\Spec(A)$. By the results of
Bloch and Ogus \cite{BO}, since $X$ is smooth, the sequence
$$
H^i(\Spec(\mathcal O_{X,x}),\z/n\z(j)) \stackrel{j^\ast}\rightarrow H^i(\Spec(F),\zn(j))
\rightarrow \mathop{\textstyle\bigoplus}\limits_{\scriptstyle x}\,H^{i-1}(k(x),\z/n\z(j-1))
$$
is exact, where $x$ runs over the codimension 1 points of $\Spec(\mathcal O_{X,x})$.
Therefore any element in the kernel of 3.7 lies in the image of $j^\ast$, hence in the
image of $H^i(\mbox{Spec}\,A,\zn(j))\rightarrow H^i(F,\zn(j))$, by the above factorization.
Since $A$ was arbitrary, the element lies in the unramified cohomology.

\v\noi

The second main result of the present section is now:

\vspace{0,5cm}

\noindent\textbf{Theorem 3.8} Let $K$ be a global field, let $n \in \mathbb N$
be invertible in $K$, and let
$F$ be a function field in $d$ variables over $K$, $d>0$, such
that $K$ is separably closed in $F$. For every place $v$ of
$K$ let $K_{(v)}$ be the Henselization of $K$ at $v$, and let
$F_{(v)}= FK_{(v)}$ be the corresponding function field over
$K_{(v)}$. Then the restriction maps induce an isomorphism
$$
\begin{CD}
H^{d+2}_{nr}(F/K,\zn(d+1)) @>{\sim}>> \mathop{\bigoplus}\limits_v\, H^{d+2}_{nr}(F_{(v)}/K_{(v)},\zn(d+1))\, .
\end{CD}
$$

\vspace{0,5cm}

\noindent\textbf{Proof} It suffices to consider the case $n = \ell^m$, where
$\ell$ is a prime invertible in $K$. Moreover it suffices to show that the map
$$
\begin{CD}
H^{d+2}_{nr}(F/K,\mathbb Q_\ell/\mathbb Z_\ell(d+1)) @>>> \mathop{\bigoplus}\limits_v \, H^{d+2}_{nr}(F_{(v)}/K_{(v)},\mathbb Q_\ell/\mathbb Z_\ell(d+1))
\end{CD}\leqno{(3.9)}
$$
is an isomorphism. In fact, if this holds, the bijectivity for $n=\ell^m$
follows from the commutative diagram with exact columns
$$
\begin{CD}
H^{d+2}_{nr}(F/K,\mathbb Q_\ell/\mathbb Z_\ell(d+1)) @>>> \mathop{\bigoplus}\limits_v \, H^{d+2}_{nr}(F_{(v)}/K_{(v)},\mathbb Q_\ell/\mathbb Z_\ell(d+1)) \\
@AA{\ell^m}A @AA{\ell^m}A   \\
H^{d+2}_{nr}(F/K,\mathbb Q_\ell/\mathbb Z_\ell(d+1)) @>>> \mathop{\bigoplus}\limits_v \, H^{d+2}_{nr}(F_{(v)}/K_{(v)},\mathbb Q_\ell/\mathbb Z_\ell(d+1)) \\
@AAA   @AAA  \\
H^{d+2}_{nr}(F/K,\zlm(d+1)) @>>>  \mathop{\bigoplus}\limits_v \, H^{d+2}_{nr}(F_{(v)}/K_{(v)},\zlm(d+1)) \\
@AAA   @AAA \\
0   @.   0
\end{CD}
$$
The exactness of the columns follows from Lemma 3.7 and the exactness of
$$
0 \rightarrow H^{i+1}(L,\zlm(i)) \longrightarrow H^{i+1}(L,\qzl(i)) \stackrel{\ell^m}\longrightarrow H^{i+1}(L,\qzl(i))
$$
for any field $L$ and any natural number $i$, which in turn follows from the theorem of Rost and Voevodsky, i.e.,
the proof of the Bloch-Kato conjecture $BK(L,i,\ell)$, see the introduction.

\medskip
We know already from Theorem 2.10 that (3.9) is injective; therefore it suffices to show
the surjectivity in (3.9).

\medskip
{\bf Case 3.8.1.} First assume that there is a smooth projective variety $X$ over $K$ with function field
$K(X)=F$. This is certainly the case if $K$ is a number field. In fact, there is a
geometrically irreducible variety $U$ over $K$ with $K(U)=F$, and after possibly shrinking $U$
we may assume that $U$ is smooth. Then, by resolution of singularities (more precisely
by property RS2$(U)$ from the beginning of section 2) it can be
embedded in a smooth projective variety $X$ over $K$ as an open
subvariety. Then, abbreviating $H^i(?,j)$ for $H^i(?,\mathbb Q_\ell/\mathbb Z_\ell(j))$, Proposition 3.7
gives a commutative diagram with exact rows

\smallskip\begin{footnotesize}
$$
\begin{CD}
0   @>>>  \bigoplus\limits_v H^{d+2}_{nr}(F_{(v)}/K_{(v)},d+1)   @>>>  \bigoplus\limits_v H^{d+2}(F_{(v)},d+1)
 @>{\bigoplus\limits_v \delta_{X_{(v)}}}>>  \bigoplus\limits_v \bigoplus\limits_{y\in X^1_{(v)}} H^{d+1}(K_{(v)}(y),d)\\
 @.   @AA{\beta'}A   @AA{\beta(F/K)}A  @AA{\beta''}A  \\
 0  @>>>  H^{d+2}_{nr}(F/K,d+1)  @>>>  H^{d+2}(F,d+1)  @>{\delta_X}>>  \bigoplus\limits_{x\in X^1} H^{d+1}(K(x),d)
\end{CD} \leqno{(3.10)}
$$
\end{footnotesize}

\noi
where $X_{(v)}=X\times_K K_{(v)}$ and in which $\beta(F/K)$ is the
restriction map, $\beta ''$ is induced by the restrictions for the
field extensions $K_{(v)}/K(x)$ for $y$ lying above $x$, and
$\beta '$ is the induced map. Note that $F_{(v)} \cong
F\otimes_KK_{(v)}$ is the function field of $X_{(v)}$ over
$K_{(v)}$. The commutativity of the right square is easily checked
(contravariance of Gysin sequences for pro-\'etale maps). Now, for
$x\in X^1$ and a place $v$ of $K$, every $y\in X_{(v)}$ lying
above $x$ is again of codimension $1$, since $X_{(v)}\rightarrow
X$ is integral. Hence
$$
\raisebox{-23pt}{$\stackrel{\coprod}{\stackrel{y\in
X^1_{(v)}}{\scriptstyle y\mid x}}$} \mbox{Spec}(K_{(v)}(y))=
\raisebox{-20pt}{$\stackrel{\coprod}{\stackrel{y\in
X_{(v)}}{\scriptstyle y\mid x}}$}\mbox{Spec}(K_{(v)}(y))=
X_{(v)}\times_X K(x)\; ,\leqno{(3.11)}
$$
the fibre of the pro-\'{e}tale morphism $X_{(v)}\rightarrow X$
over $x$. This is again isomorphic to
$$
\begin{array}{rcl}
(X\times_K K_{(v)})\times_X  K(x) & \cong & \mbox{Spec} (
K(x)\otimes_K  K_{(v)})\\
& \cong & \mbox{Spec}( K(x)\otimes_{K\{x\}} (K\{x\}\otimes_K
K_{(v)}))\\
& \cong & \coprod\limits_{w \mid v}\mbox{Spec}(
K(x)\otimes_{K\{x\}}K\{x\}_{(w)})
\end{array}\leqno{(3.12)}
$$
where $K\{x\}$ is the separable closure of $K$ in $K(x)$ (which is
a finite extension of $K$) and where $w$ runs over the places $w$
of $K\{x\}$ above $v$. This shows that $\beta''$ can be identified
with the map
$$
\mathop{\textstyle\bigoplus}\limits_{x\in X^1} \beta
(K(x)/K\{x\}):\,\mathop{\textstyle\bigoplus}\limits_{x\in
X^1}\,H^{d+1}(K(x),d) \longrightarrow
\mathop{\textstyle\bigoplus}\limits_{x\in
X^1}\mathop{\textstyle\bigoplus}\limits_{w\in
P(K\{x\})}\,H^{d+1}(K(x)_{(w)},d)\leqno{(3.13)}
$$
where $P(K\{x\})$ is the set of places of the global field
$K\{x\}$. Hence $\beta''$ is injective as well as $\beta(F/K)$, by
Theorem 2.10. (Note that $K(x)$, for $x\in X^1$, is a function
field in $d-1$ variables over $K\{x\})$. By diagram (3.10) it
now suffices to show that the following map is injective.
$$
\mbox{coker }\beta(F/K)\longrightarrow \mbox{coker } \beta'' =
\mathop{\textstyle\bigoplus}\limits_{x \in X^1}\mbox{coker
}\beta(K(x)/K\{x\})\, . \leqno{(3.14)}
$$

\vspace{0,5cm}

\noindent\textbf{Lemma 3.15} The map (3.14) can be identified with
the map of cofixed modules
$$
H^d(F\overline K,d)_{G_K} \longrightarrow
\mathop{\textstyle\bigoplus}\limits_{x\in
X^1}\,H^{d-1}(K(x)\otimes_K\overline K,d-1)_{G_K} \cong
(\,\mathop{\textstyle\bigoplus}\limits_{y\in
\overline{X}^1}\,H^{d-1}(\overline{K}(y),d-1)\,)_{G_K}
$$
induced by the residue map $\delta_{\overline{X}}$ for
$\overline{X} = X\times_K\overline K$.

\vspace{0,5cm}

\noindent\textbf{Proof} By (3.11) and (3.12), the map $\beta''$
can also be identified with the map
$$
\mathop{\textstyle\bigoplus}\limits_{x\in X^1}\, [\;\beta
(K(x)/K):\,H^{d+1}(K(x),d) \longrightarrow
\mathop{\textstyle\bigoplus}\limits_{v\in
P(K)}\,H^{d+1}(K(x)\otimes_KK_{(v)},d)\;]\,.
$$
\noindent
Therefore, the map (3.14) can be identified with the map
 $\mbox{ coker }\beta_1 \rightarrow \mbox{ coker } \beta_2\;$ induced by the commutative diagram
$$
\xymatrix{\raisebox{-6pt}{$\stackrel{\textstyle\bigoplus}{\scriptstyle
v}$}~H^2(K_v,H^d(F\overline K,d+1))\ar[r] &
\bigoplus\limits_{x\in X^1}\bigoplus\limits_{
v}H^2(K_v,H^{d-1}(K(x)\otimes_K\overline K,d))\\
H^2(K,H^d(F\overline K,d+1))\ar[r]\ar[u]_{\beta_1} &
\bigoplus\limits_{x\in
X^1}H^2(K,H^{d-1}(K(x)\otimes_K\overline K,d))\ar[u]_{\beta_2}
\,, }
$$
where the vertical maps are the obvious restriction maps, and the horizontal maps are induced by the residue maps
$$
H^d(F\overline K,d+1) \longrightarrow
\mathop{\textstyle\bigoplus}\limits_{x\in X^1}\,
H^{d-1}(K(x)\otimes_K\overline K,d) \cong
\mathop{\textstyle\bigoplus}\limits_{y\in
\overline{X}^1}\,H^{d-1}(\overline{K}(y),d)
$$
for $\overline{X}$. This follows from Proposition 1.2, Remark 1.4
(a) and the fact that the Hochschild-Serre spectral sequence is
compatible with the connecting morphisms for Gysin sequences. The
latter statement follows from the fact that the Hochschild-Serre
spectral sequence for \'etale (hyper)cohomology of complexes is
functorial with respect to morphisms in the derived category, and
that the Gysin isomorphisms are compatible with pro-\'etale base
change.

\noindent Finally, for all discrete torsion $\mathbb Z_\ell$-$G_K$-modules
$M$ there are canonical isomorphisms
$$
\coker[\;\beta_M:\,H^2(K,M) \rightarrow
\textstyle\mathop{\textstyle\bigoplus}\limits_{v}\,H^2(K_v,M)\;]\; \mathop{\longrightarrow}\limits^{\sim} \; M(-1)_{G_K}\,,
\leqno{(3.16)}
$$
which are functorial in $M$, see (the proof of) Proposition 1.2. This proves Lemma 3.15.

\vspace{0,5cm}

We are now ready to prove Theorem 3.8 (in our situation, assuming the existence
of $X$ with $F=K(X)$ as above). By Lemma 3.15 it suffices to show the following
more general theorem, which will also be used in the later sections.

\v\noi
\textbf{Theorem 3.17} Let $K$ be a finitely generated field with perfect hull $L$ and algebraic closure $\overline{K}$,
let $X$ be a smooth projective irreducible variety of dimension $d$ over $K$, and
let $\ell$ be a prime. Assume that $\ell$ is invertible in $K$ or that condition RS1$(U)$
(see the beginning of section 2) holds for any open $U \subset X\times_K L$. Then the map
$$
(\mathop{\textstyle\bigoplus}_{y \in \overline{X}^0} H^d(\overline{K}(y),\qzl(d)))_{G_K} \longrightarrow
(\mathop{\textstyle\bigoplus}_{x \in \overline{X}^1}\;H^{d-1}(\overline{K}(x),\qzl(d-1)))_{G_K}
$$
induced by the Bloch-Ogus complex for $\overline{X}=X\times_K\overline{K}$ (via taking coinvariants under $G_K$)
is injective.

\v\noi
\textbf{Proof} We note that here we regard the absolute Galois group $G_K$ of $K$ as $Gal(\overline{K}/L)$,
and we may replace $K$ by $L$ and call this $K$ again. Moreover the above map can be identified with a map
$$
H^d(K(X)\otimes_K\overline{K},\qzl(d))_{G_K} \longrightarrow
\mathop{\textstyle\bigoplus}_{x\in X^1}\;H^{d-1}(K(x)\otimes_K\overline{K},\qzl(d-1))_{G_K}\,.
$$
Let $a$ be an element in the kernel of the above map. Then there is an open $U\subset X$ such that
$a$ is the image of an element $a_U\in H^d(\overline{U},\qzl(d))_{G_K}$, where $\overline{U}=U\times_K\overline{K}$.

\medskip
{\bf Case 3.17.1.} First assume that RS1$(U)$ holds and that $K$ is infinite. Then there exists an open embedding
$U \subset X'$ into a smooth projective variety $X'$ such that the complement $Y= X'\smallsetminus U$ is a
divisor with simple normal crossings, say with smooth components $Y_i$
$(i=1,\ldots,r)$. By possibly applying Bertini's theorem as in section 2
(before Proposition 2.2) and removing a suitable hyperplane
section (which does not matter for our purposes), we may assume
that $X' \smallsetminus Y_1$ is affine, i.e., that $U\subseteq
X' \supseteq Y$ satisfies the assumptions of Theorem 3.1.

\medskip
Next we note that the kernel of the map in 3.17 only depends on $F=K(X)$ and not on the smooth
projective model $X$ of $F$. In the case of a global field $K$ and $\ell$ invertible in $K$
this is clear from Lemma 3.15, the diagram (3.10) and Proposition 3.7 for $F$ and the $F_{(v)}$.
In general, the argument is the same as in the proof of Proposition 3.7, noting the following two facts.
By properness, for any $x' \in (X')^1$ the discrete valuation ring $\cO_{X',x'}$ dominates a local ring
$\cO_{X,y}$ of $X$. Moreover, for this regular ring and any finite Galois extension $M/K$ the
Bloch-Ogus sequence of $Gal(M/K)$-modules
%\begin{footnotesize}
$$
H^d(\Spec(\mathcal O_{X,y}\otimes_KM),\Lambda(d)) \stackrel{j^\ast}\rightarrow H^d(K(X)\otimes_KM),\Lambda(d))
\rightarrow \mathop{\textstyle\bigoplus}\limits_{\scriptstyle x \in X^1}\,H^{d-1}(K(x)\otimes_KM,\Lambda(d-1))
$$
%\end{footnotesize}
\noi
is exact for $\Lambda=\qzl$, by purity for the semi-local ring $\cO_{X,y}\otimes_KM$. In addition, it stays exact after taking
coinvariants under $Gal(M/K)$, because the Gersten resolution is universally exact (see \cite{CTHK} Cor. 6.2.4
together with Ex. 7.3 (1) (for $\ell \neq \ch(K)$) and Ex. 7.4 (3) (for $\ell = \ch(K)$ and the Tate twist $d$).
By passing to the inductive limit we get the corresponding result for $\overline{K}$ and $G_K$ in place of
$M$ and $Gal(M/K)$. The same holds for the discrete valuation ring $\mathcal O_{X',x'}$.
As in the proof of 3.7 we get that the kernel of 3.17 for $X$ lies in the kernel of 3.17 for $X'$.
Interchanging the roles of $X$ and $X'$ we get the wanted equality.

\medskip
Therefore we may replace $X$ above by $X'$ and call it $X$ again.
Now we claim that $a_U$ lies in the kernel of the map
$$
e:H^d(\overline U,\qzl(d))_{G_K}\rightarrow H^0(\overline{Y^{[d]}},\qzl(0))_{G_K}
$$
introduced in Theorem 3.1. Since the assumptions of 3.1 are
fulfilled for $U$, we then conclude that $a_U$ is zero
and hence $a$ is zero as wanted.

\medskip
With the notations of
(3.2), the claimed vanishing of $e(a_U)$ follows from
the following commutative diagram, for each $(i_1,\ldots, i_d)$
and each $y\in Y_{i_1,\ldots , i_d}$.
%, where we write $H^-(-,j)$ for $H^-(-,\qzl(j))$, as before.

\smallskip

$$
\begin{CD}
H^d(\overline U,d)_{G_K}  @>>>  H^{d-1}(\overline{Y_{i_d} \!\smallsetminus\! (\bigcup\limits_{i\neq i_d}Y_i)}, d-1)_{G_K}
@>>>  H^{d-2}(\overline{Y_{i_{d-1}, i_d}\!\smallsetminus\!(\!\!\!\bigcup\limits_{i\neq i_{d-1},i_d}\!\!\!Y_i}) ,d-2)_{G_K} \\
@VVV      @VVV    @VVV  \\
H^d(F\overline K,d)_{G_K}  @>>>   H^{d-1}(K(y_{i_d})\otimes_K\overline K,d-1)_{G_K} @>>>
H^{d-2}(K(y_{i_{d-1},i_d})\otimes_K\overline K,d-2)_{G_K}
\end{CD}\leqno{(3.18)}
$$
$$
\xymatrix{ \ldots\ar[r]  & H^0(\overline {\{y\}},0)_{G_K}\ar@{=}[d]   \\
\ldots\ar[r]  & H^0(K(y_{i_1,\ldots,i_d})\otimes_K\overline
K,0)_{G_K} }
$$
in which $y_{\underline i}$ is the generic point of the component
$Y_{\underline i}^y$ of $Y_{\underline i}$ in which $y$ lies, for
any $\underline i=(i_r,\ldots,i_d)$, so that $K(y_{\underline i})$
is the function field of $Y_{\underline i}^y$. In fact, the maps
in the bottom line are all induced by the residue maps, by
definition, and the image of $\overline{a_U}$ under the left
vertical map is $\overline a$. As we have noted, the image of
$\overline a$ in $H^{d-1}(K(y_{i_d})\otimes_K\overline
K,d-1)_{G_K}$ vanishes, for every choice of $(i_1,\ldots,i_d),\,
1\leq i_1<\ldots <i_d\leq r$ (note that $y_{i_d}\in \tilde X^1)$.
Therefore the image of $\overline{a_U}$ in
$H^0(\overline{\{y\}},0)_{G_K}$ vanishes for every $y\in
\overline{Y^{[d]}}$ as claimed.

\v
{\bf Case 3.17.2.} For the case of a finite field $K$ we note the following. First of all we have
canonical functorial isomorphisms $M_{G_K}\cong H^1(K,M)$ for all discrete $G_K$-modules $M$.
Therefore the map in 3.17 can be identified with the map
$$
H^{d+1}(K(X),\qzl(d)) \rightarrow \mathop{\textstyle\bigoplus}\limits_{x\in X^1} H^d(k(x),\qzl(d-1))\,,
$$
and it follows directly from Proposition 3.7 that the kernel of this map is
independent of $X$ and in fact equal to the unramified cohomology $H^{d+1}_{nr}(K(X)/K,\qzl(d))$.
If $a$, $U$ and $a_U$ are as above, and we have a good compactification $U \subset X \supset Y$ as above,
we may not have a suitable hyperplane section defined over $K$, but we get one after
taking a base extension to a field extension $K'/K$ of degree prime to $\ell$. Then we conclude
that $a$ maps to zero under the restriction $Res: H^{d+1}(K(X),\qzl(d)) \rightarrow H^{d+1}(K'(X),\qzl(d))$.
But this map is injective, by the existence of the corestriction $Cor$ in the other direction with
$Cor Res =$ multiplication with $[K':K]$ which is prime to $\ell$.

\v
{\bf Case 3.17.3.} Finally consider the case that $\ch(K)=p >0$ and $\ell\neq p$, and that we have no
good compactification of $U$, where $a$, $U$ and $a_U$ are as above.
By the weaker resolution of singularities of Gabber-Illusie
(property (G) in the proof of Theorem 2.10) we get a diagram
$$
\begin{array}{ccccc}
U' & \subset & X' & \supset & Y' \\
\downarrow  &&  \downarrow\rlap{$\pi$}  &&  \downarrow  \\
U & \subset & X & \supset & Y
\end{array}\eqno{(3.19)}
$$
where $X'$ is an geometrically irreducible smooth and projective variety over a finite extension $K'$ of $K$
with $\ell$ not dividing $[K':K]$, $\pi$ is a proper surjective morphism which is generically finite of degree prime to $\ell$,
$U' = \pi^{-1}(U)$, and $Y'=\pi^{-1}(Y)$ is a divisor with strict normal crossings on $X'$.
Since $\ell\neq p$ the $\qzl$-cohomology does not change under radicial maps, and we may pass to the perfect hull
of $K$ in $K'$ and thus assume that $K'/K$ is separable. Then $X$ and $X'$ are smooth projective over $K$.

For any smooth variety $V$ over $K$, let $\overline{\mathcal H}^i(d)_{G_K}$ be the Zariski sheaf on $V$ associated to
the presheaf $U \mapsto H^i(\overline{U},\qzl(d))_{G_K}$ for $U \subseteq V$ open, where $\overline{U}=U\times_K\overline{K}$.
Then the Bloch-Ogus theory and the universal exactness used above show that the kernel of the map in Theorem 3.17
is canonically isomorphic to $H^0(X,\overline{\mathcal H}^d(d)_{G_K})$, and the pull-back maps for \'etale cohomology
induce a natural pull-back map
$$
\pi^\ast: \; H^0(X,\overline{\mathcal H}^d(d)_{G_K}) \rightarrow H^0(X',\overline{\mathcal H}^d(d)_{G_K})\,.
$$
We claim that this map is injective. In fact, this maps embeds into the restriction map for
the function fields, which can be factored as
$$
H^d(K(X)\otimes_K\overline{K})_{G_K} \rightarrow H^d(K'(X)\otimes_K\overline{K})_{G_K} \rightarrow H^d(K'(X')\otimes_{K'}\overline{K})_{G_{K'}}\,.\leqno{(3.20)}
$$
where we have omitted the coefficients $\qzl(d)$. Both restriction maps are injective, because the degrees
$[K'(X):K(X)]$ and $[K'(X'):K'(X)]$ are prime to $\ell$ (see the corestriction argument at the end of the
proof of Theorem 2.10, which also works for the modules of coinvariants).
Thus the restriction map is injective and $\pi^\ast$ is injective as well.

\medskip
Now the element $a \in H^0(X,\overline{\mathcal H}^i(d)_{G_K})$ maps to an element $a'\in H^0(X',\overline{\mathcal H}^i(d)_{G_K})$,
which is the image represented by the image $a_{U'}$ of $a_U$ under the restriction map
$H^d(\overline{U},\qzl(d))_{G_K} \rightarrow H^d(\overline{U'},\qzl(d))_{G_K'}$ which is defined in
the same way as for (3.19). By the arguments given for the previous case of the proof we get that
$a'=0$, and so $a=0$ by the injectivity of (3.20). This finishes the proof of Theorem 3.17.

\v
{\bf Case 3.8.2} With similar arguments we can now also complete the proof of Theorem 3.8, in the case where
the function field $F$ over $K$ does not have any smooth projective model, but the prime $\ell$ is
invertible in $K$. Let $U$ be an affine integral geometrically irreducible variety of dimension $d$
over $K$ with function field $K(U)=F$. Then we have an open embedding $U \subset X$ into a projective integral variety,
and we get again a diagram as in (3.19) (the smoothness of $U$ or $X$ was not needed).
Let $F'=K'(X')$. It follows from the definition of unramified cohomology that the morphism
$F/K \rightarrow F'/K'$, i.e., the commutative diagram
$$
\begin{array}{ccc}
F & \rightarrow & F' \\
\uparrow &  & \uparrow \\
K  & \rightarrow & K'
\end{array}
$$
induces a restriction map $H_{nr}^{d+1}(F/K,\qzl(d))\longrightarrow H_{nr}^{d+1}(F'/K',\qzl(d))$.
(Any discrete valuation of $F'$ over $K'$ induces by restriction a discrete valuation of $F$ over $K$.)
The same holds for the morphism $F_{(v)}/K_{(v)} \rightarrow F'_{(w)}/K'_{(w)}$ for a place
$v$ of $K$ and a place $w$ of $K'$ above $v$. Moreover, the commutative diagrams
$$
\begin{array}{ccc}
F_{(v)}/K_{(v)} & \rightarrow & F'_{(w)}/K'_{(w)} \\
\uparrow & & \uparrow \\
F/K  &  \rightarrow  & F'/K'
\end{array}
$$
induces a commutative diagram
$$
\begin{CD}
\textstyle\bigoplus_v\;H_{nr}^{d+1}(F_{(v)}/K_{(v)},\qzl(d)) @>>>  \textstyle\bigoplus_w\; H_{nr}^{d+1}(F'_{(w)}/K'_{(w)},\qzl(d)) \\
@AA{\beta'(F/K)}A @AA{\beta'(F'/K')}A  \\
H_{nr}^{d+1}(F/K,\qzl(d))  @>>>   H_{nr}^{d+1}(F'/K',\qzl(d))\,.
\end{CD}
$$
By the first part of the proof (having the existence of the smooth projective model $X'$ of $F'$),
the cokernel of $\beta'(F'/K')$ is zero. Now we claim that the map
$\coker\; \beta'(F/K) \rightarrow \coker\;\beta'(F'/K')$ induced by the diagram is injective;
then we have $\coker\;\beta'(F/K)=0$ as wanted. First of all, the above diagram is obtained
from the following commutative diagram of restriction maps
$$
\begin{CD}
\textstyle\bigoplus_v\;H^{d+1}(F_{(v)}/K_{(v)},\qzl(d)) @>>>  \textstyle\bigoplus_w\; H^{d+1}(F'_{(w)}/K'_{(w)},\qzl(d)) \\
@AA{\beta(F/K)}A @AA{\beta(F'/K')}A  \\
H^{d+1}(F/K,\qzl(d))  @>>>   H^{d+1}(F'/K',\qzl(d))
\end{CD}
$$
by passing to the unramified subgroups, so that we have a commutative diagram
$$
\begin{CD}
\coker\;\beta(F/K) @>{r}>> \coker\;\beta(F'/K') \\
@AA{i}A   @AA{i'}A  \\
\coker\;\beta'(F/K) @>{r'}>> \coker\;\beta'(F'/K')\,.
\end{CD}
$$
We claim that the maps $i$ and $r$ are injective; then we obtain the injectivity of $r'$.

\v
The injectivity of $i$ follows from the commutative diagram with exact rows
$$
\begin{CD}
0 @>>> \mathop{\textstyle\bigoplus}\limits_v\;H^{d+2}_{nr}(F_{(v)}/K_{(v)},d+1)
@>>> \mathop{\textstyle\bigoplus}\limits_v\;H^{d+2}(F_{(v)},d+1)
@>{\oplus_v\delta_v}>> \mathop{\textstyle\bigoplus}\limits_v\mathop{\textstyle\bigoplus}\limits_{\mu\in P_v}\;H^{d+1}(K_{(v)}(\mu),d)\\
@.  @AA{\beta'(F/K)}A @AA{\beta(F/K)}A @AA{\beta''}A \\
0 @>>> H^{d+2}_{nr}(F/K,d+1) @>>> H^{d+2}(F,d+1) @>{\delta}>> \mathop{\textstyle\bigoplus}\limits_{\lambda\in P}\;H^{d+1}(K(\lambda),d)\,.
\end{CD}
$$
Here we have omitted the coefficients $\qzl$, $P$ (resp. $P_v$) is the set of discrete valuations of $F/K$ (resp. $F_{(v)}/K_{(v)}$),
and the components of $\delta$ (resp. $\delta_v$) are the residue maps for the valuations $\lambda \in P$ (resp. $\mu\in P_v$).
The restriction map $\beta''$ is defined as follows. If the valuation $\mu$ of $F_{v)}/K_{(v)}$ lies over the valuation $\lambda$
of $F/K$, then the corresponding component is induced by the inclusion of the corresponding valuation rings; otherwise the
component is zero. The map $\beta''$ is injective, by similar arguments as in the beginning of the proof of 3.8: If $A_\lambda$
is the valuation ring of $\lambda$ and $M/K$ is a finite separable field extension, then $A_\lambda\otimes_KM$ is a regular
semi-local ring of dimension 1, and hence the integral closure of $A_\lambda$ in $F\otimes_KM$. Hence the extensions of the
valuation $\lambda$ to $F\otimes_KM$ correspond to the fiber above the closed point of $A_\lambda$, i.e., to the points of
$\Spec(K(\lambda)\otimes_KM)$. Therefore the extensions of $\lambda$ to $F_{(v)}$ correspond to the points of
$Spec(K(\lambda)\otimes_KK_{(v)})\cong \oplus_w\,\Spec(K\{\lambda\}_{(w)})$, where $K\{\lambda\}$ is the separable
closure of $K$ in $K(\lambda)$ (which is a finite extension), and $w$ runs over all places of $K\{\lambda\}$ lying above $v$.
Thus the restriction of $\beta''$ to the component for $\lambda$ can be identified with the map
$\beta(K(\lambda)/K\{\lambda\})$, which is injective by Theorem 2.10. (Note that $K(\lambda)$ is a geometrically irreducible
function field in $d-1$ variables over $K\{\lambda\}$.) The injectivity of $\beta''$ now implies by a diagram chase that
$i: \coker\;\beta'(F/K) \rightarrow \coker\;\beta(F/K)$ is injective.

\v
Now we consider the injectivity of $r$. Since the $\ell$-adic cohomology does not change under
radicial/inseparable extensions, we may assume that $F'/F$ and $K'/K$ are separable.
Then we get a commutative diagram with exact rows
$$
\begin{CD}
H^{d+2}(F',d+1) @>>> \mathop{\textstyle\bigoplus}\limits_v\;H^{d+2}(F'\otimes_KK_{(v)},d+1) @>>> \coker\;\beta(F'/K') @>>> o \\
@AA{Res}A   @AA{Res}A   @AA{r}A  @. \\
H^{d+2}(F,d+1) @>>> \mathop{\textstyle\bigoplus}\limits_v\;H^{d+2}(F\otimes_KK_{(v)},d+1) @>>> \coker\;\beta(F/K) @>>> 0 \\
\end{CD}
$$
induced by the restriction for $F'/F$. The cokernel in the upper row can indeed
be identified with $\coker\;\beta(F'/K')$, compare Remark 1.4 (a), and then the right hand map
can be identified with $r$ as indicated. Now the finite \'etale map $\pi: \Spec(F')\rightarrow \Spec(F)$
also induces compatible downward maps $\pi_\ast$ in the left square, such that $\pi_\ast$ is the usual corestriction Cor
on the left and such that $\pi_\ast Res$ is the multiplication with $[F':F]$ in both cases. Since
$F'/F$ and the extensions $K_{(v)}/K$ are separable, these properties follow from obvious calculations
in Galois cohomology which are left to the readers. We obtain an induced map $\pi_\ast$ on the right with
$\pi_\ast r$ being the multiplication with $[F':F]$. Since this degree is prime to $\ell$,
and we consider $\qzl$-coefficients, we obtain the injectivity of $r$ as claimed.

\v\noindent
\textbf{Remark 3.18} In the considerations of this section we have preferred to work with
the explicit descriptions (via Gysin and specialization maps) of the
maps $e$ and $d_2$ in Theorem 3.1, but we note that they coincide with
the corresponding edge morphisms and differentials of the weight spectral sequence (2.3) for
$U \subset X \supset Y$, up to signs.

\vspace{15mm}

\centerline{\noindent\textbf{\S 4 A Hasse principle for
Bloch-Ogus-Kato complexes}}

\vspace{1,0cm}

Let $X$ be an excellent scheme, let $n\geq 1$ be an integer,
and let $r,s\in\mathbb Z$. Under some conditions on $X$, $n$ and $(r,s)$,
there are homological complexes of Gersten-Bloch-Ogus-Kato type

\medskip\noindent
$C^{r,s}(X,\zn):\\
\phantom{C^{r,s}}\ldots \rightarrow \bigoplus\limits_{x\in X_i}H^{r+i}(k(x),\zn(s+i))\stackrel{\partial}{\rightarrow} \bigoplus\limits_{x\in X_{i-1}} H^{r+i-1}(k(x),\zn(s+i-1))\rightarrow\ldots\\
\phantom{C^{r,s}\ldots \rightarrow \bigoplus\limits_{x\in
X_i}H^{r+i}}\ldots \rightarrow  \bigoplus\limits_{x\in X_0}
H^r(k(x),\zn(s)) $

\medskip\noindent
where the term for $X_i=\{x\in X\mid\dim (x)=i\}$ is placed in
degree $i$.

\medskip\noindent
If $X$ is separated of finite type over a field $L$, $n$ is invertible in $L$
and $(r,s)$ are arbitrary, these complexes where defined by Bloch and Ogus [BO],
by using the \'{e}tale homology for such schemes, defined by
$$
H_a(X,\zn(b)):= H^{-a}(X,Rf^!\zn(-b)),\leqno{(4.1)}
$$
where $Rf^!$ is the functor on constructible \'{e}tale $\zn$-schemes defined in [SGA 4] XVIII, for the structural morphism
$f:X\rightarrow$ Spec $L$. In fact, Bloch and Ogus constructed a
niveau  spectral sequence
$$
E^1_{p,q} (X,\zn (b)) = \textstyle\bigoplus\limits_{x\in
X_p}H_{p+q}(k(x),\zn(b)) \Rightarrow H_{p+q}(X, \zn(b))
$$
where, by definition, $H_a(k(x),\zn(b)) =
\lim\limits_{\raisebox{2pt}{$\longrightarrow$}}  H_a(U, \zn(b))$ for $x\in X$, where the limit is over all open subschemes
$U\subseteq \overline{\{x\}}$ of the Zariski closure of $x$. By
purity, there is an isomorphism $ H_a(U, \zn(b)) \cong
H^{2p-a} (U, \zn (p-b)) $ for $U$ irreducible and smooth
of dimension $p$ over $L$. Thus one has a canonical isomorphism
$$
E^1_{p,q} (X,\zn (b)) \cong
\textstyle\bigoplus\limits_{x\in X_p}H^{p-q}(k(x),\zn(p-b))\, .\leqno{(4.2)}
$$
This is clear for a perfect field $L$, because then
$\overline{\{x\}}$ is generically smooth. So the limit can be
carried out over the smooth $U\subseteq \overline{\{x\}}$, and
$$
\lim\limits_{\raisebox{2pt}{$\longrightarrow$}} \quad
H^{2p-a}(U,\zn(p-b)) = H^{2p-a}(k(x),\zn(p-b)),
$$\vspace{-0,9cm}
$$
\hspace{-9,3cm}\scriptstyle U\subseteq \overline{\{x\}}
$$
since $\lim\limits_{\raisebox{2pt}{$\longleftarrow$}} U=$ Spec
$k(x)$, and since \'{e}tale cohomology commutes with this limit.
For a general field $L$, we may pass to the separable hull,
because of invariance of \'{e}tale cohomology with respect to base
change with radical morphisms.

\medskip
Using the identification (4.2), one may define
$$
C^{r,s}(X,\zn) = E^1_{\cdot,-r}(X,\zn(-s))\, .
$$
With this definition, one obtains the following description of the
differential
$$
\partial: \textstyle\bigoplus\limits_{x\in X_i} H^{r+i}(k(x),\zn(s+i)) \rightarrow \textstyle\bigoplus\limits_{x\in X_{i-1}} H^{r+i-1}(k(x),\zn(s+i-1))\, .
$$
We may assume that $L$ is perfect. For $y\in X_i$ and $x\in
X_{i-1}$ let $\partial_{y,x} = \partial^X_{y,x}$ be the
$(y,x)$-component of $\partial$. If $x\notin \overline{\{y\}}$,
then $\partial_{y,x}=0$. If $x$ is a smooth point of
$\overline{\{y\}}$, then there is an open smooth neighbourhood
$x\in U\subseteq \overline{\{y\}}$. Moreover, any $\alpha \in
H^{r+i}(k(x),\zn(s+i))$ lies in the image of $
H^{r+i}(V,\zn(s+i))\rightarrow H^{r+i}(k(x),\zn(s+i)) $ for some open $V\subseteq U$. Moreover, by making $U$
(and $V$) smaller we may assume that $Z=U\setminus V$ is
irreducible and smooth as well, and that $x$ is the generic point
of $Z$. Then one has a commutative diagram
$$
\xymatrix{H^{r+i}(k(y),\zn(s+i)\ar[r]^{\hspace{-0,7cm}\partial_{y,x}} & H^{r+i-1}(k(x),\zn(s+i-1))\\
H^{r+i}(V,\zn(s+i)\ar[r]^{\hspace{-0,7cm}\partial }\ar[u]
& H^{r+i-1}(Z,\zn(s+i-1))\ar[u]}\leqno{(4.3)}
$$
where the vertical maps come from passing to the generic points,
and $\partial$ is the connecting morphism for the Gysin sequence
for $(U,Z)$. This determines $\partial_{y,x}(\alpha)$. If
$x\in\overline{\{y\}}$, but is not necessarily a smooth point of
$Y=\overline{\{y\}}$, let $\tilde Y\rightarrow Y$ be the
normalization of $Y$. Any point $x'\in\tilde Y$ above $x$ has
codimension $1$ and thus is a regular point in $\tilde Y$. Since
the niveau spectral sequence is covariant with respect to proper
morphisms, there is a commutative diagram
$$
\xymatrix{H^{r+i}(k(y),\zn(s+i))\ar[rr]^{\hspace{-0,7cm}\bigoplus\limits_{x'\mid x}\partial_{y,x'}^{\tilde Y}} & &
\textstyle\bigoplus\limits_{x'\mid x}H^{r+i-1}(k(x'),\zn(s+i-1))\ar[d]^{\pi_\ast}\\
H^{r+i}(k(y),\zn(s+i))\ar[rr]^{\hspace{-0,7cm}\partial_{y,x}}\ar@{=}[u] & &
H^{r+i-1}(k(x),\zn(s+i-1))}\leqno{(4.4)}
$$
where $\pi_\ast$ is induced by $\pi:\tilde Y \rightarrow
Y\hookrightarrow X$. One can check that
$\pi_\ast((\alpha_{x'}))=\sum\limits_{x'\mid x}\mbox{Cor}_{x'\mid
x}(\alpha_{x'})$, where $\mbox{Cor}_{x'\mid x}:
H^\mu(k(x'),\zn(\nu))\rightarrow H^\mu(k(x),\zn(\nu))$ is the corestriction for the finite extension
$k(x')/k(x)$ (this also makes sense if this extension has some
inseparable part). Since $\partial^{\tilde Y}_{y,x'}$ can be
treated as before, this determines $\partial_{y,x}$.

\vspace{0,5 cm} For a function field $L$ of transcendence degree $d$ over a perfect field $k$ of
characteristic $p>0$, a separated scheme $X$ of finite type over $L$ and $n$ a power of $p$ it was shown in
\cite{JS} and \cite{JSS} 2.11.3 that the theory of Bloch and Ogus
can be literally extended to this situation for the case $b=-d$, where the cohomology groups $H^i(X,\zn(j))$
and $H^i(k(x),\zn(j))$ are defined as in (0.2).

\vspace{0,5cm}
For a general excellent scheme $X$, and arbitrary $n$,
the complexes $C^{r,s}(X,\zn)$ were defined by Kato (and named
$C^{r,s}_n(X)$, cf. [Ka]), in a more direct way, by using the
Galois cohomology of discrete valuation fields, assuming the following condition:

\begin{itemize}
\item[($\ast$)]
If $r = s+1$ and $p$ is a prime dividing $n$, then for any $x\in X_0$ with
$\ch(k(x))=p$ one has $[k(x):k(x)^p]\leq s$.
\end{itemize}

It is shown in [JSS] that both definitions agree (up to well-defined signs)
for varieties over fields in the cases discussed above.

\vspace{0,5cm} Now let $K$ be a global field, and let $X$ be a
variety over $K$. For every place $v$ of $K$ let $X_v= X\times_K
K_v$. Then condition $(\ast)$ holds for $X$ and the $X_v$ for $(r,s) = (2,1)$
and arbitrary $n$. Moreover,
one has natural restriction maps $\alpha_v: C^{r,s}(X,\zn)\rightarrow C^{r,s}(X_v,\zn)$,
and Kato stated the following conjecture (see Conjecture 2 in the introduction).

\vspace{0,5cm}

\noindent\textbf{Conjecture 4.5} Let $X$ be connected, smooth and
proper. Then the $\alpha_v$ induce isomorphisms
$$
H_a(C^{2,1}(X,\zn))\stackrel{\sim}{\rightarrow}
\textstyle\bigoplus\limits_vH_a(C^{2,1}(X_v,\zn)) \mbox{
for all } a\neq 0,
$$
and an exact sequence
$$
0\rightarrow H_0(C^{2,1}(X,\zn))\rightarrow
\textstyle\bigoplus\limits_vH_0(C^{2,1}(X_v,\zn))\rightarrow \zn\rightarrow 0\, .
$$

\vspace{0,5cm}

\textbf{Remark 4.6} (a) For $X=$ Spec $(L)$, $L$ any finite extension
of $K$, the cohomology groups vanish for $a\neq 0$, and the
sequence for $a=0$ becomes the exact sequence
$$
0\rightarrow H^2(L,\zn(1))\rightarrow
\textstyle\bigoplus\limits_{w\in P(L)} H^2(L_w,\zn
(1))\rightarrow \zn\rightarrow 0\, ,
$$
which is the $n$-torsion of the classical exact sequence
$$
0\rightarrow Br(L)\rightarrow \textstyle\bigoplus\limits_{w\in
p(L)} Br (L_w)\stackrel{\sum\limits_w \mbox{\scriptsize
inv}_w}{\longrightarrow} \mathbb Q/\mathbb Z\rightarrow 0
$$
for the Brauer groups (where $\mbox{\scriptsize inv}_w:
Br(L_w)\stackrel{\sim}{\rightarrow}\mathbb Q/\mathbb Z$ is the
'invariant' map). Thus Kato's conjecture is a generalization of
this famous sequence to higher dimensional varieties.\\
(b) As we will see below, the $\alpha_v$ induce a map
$$
\alpha_{X,n}: C^{2,1}(X,\zn) \longrightarrow \textstyle\bigoplus\limits_v \, C^{2,1}(X_v,\zn)\,.\leqno{(4.7)}
$$
Let $C'(X,\zn)$ be its cokernel. Then conjecture 4.5 is implied to the following two statements:

\smallskip\noi
(i) $\alpha_{X,n}$ is injective.\\
(ii) $H_0(C'(X,\zn)) = \zn$, and $H_a(C'(X,\zn)) = 0$ for $a > 0$.

\smallskip\noi
Conversely conjecture 4.5 implies (i) and (ii) by the known case (a) and induction on dimension,
provided the occurring function fields have smooth and proper models
over the perfect hull of $K$ (which holds over number fields).
\vspace{0,5cm}

We prove the above conjecture for $n$ invertible, in the following form (cf. Remark 4.6 (b)).

\vspace{0,5cm}

\noindent\textbf{Theorem 4.8} Let $K$ be a global field, let $n \in \mathbb N$ be invertible in $K$,
and let $X$ be a connected smooth proper variety over $K$.

\smallskip
\noi(a) The map $\alpha_{X,n}: C^{2,1}(X,\zn) \longrightarrow \mathop{\textstyle\bigoplus}_v\;C^{2,1}(X_v,\zn)$
is well-defined and injective.

\smallskip
\noi(b) Let $C'(X,\zn)$ be the cokernel of $\alpha_{X,n}$. Then
$$
H_a(C'(X,\z/n\z))=\left\{\begin{array}{ccc} 0 & , & a\neq 0,\\
\z/n\z & , & a=0\, .\end{array}\right.
$$

\vspace{0,5cm}\noi

\noi\textbf{Proof of Theorem 4.8 (a):} First note that the restriction map $\alpha_v$ factors as
$$
\alpha_v: C^{2,1}(X,\zn)\mathop{\longrightarrow}\limits^{\beta_v}
C^{2,1}(X_{(v)},\zn)\rightarrow
C^{2,1}(X_{v},\zn)
$$
where $X_{(v)} = X\times_K K_{(v)}$. These maps of complexes have
components
$$
\textstyle\bigoplus\limits_{x\in X_i}H^{i+2}(k(x),\zn(i+1))
\rightarrow\hspace{-2mm}\bigoplus\limits_{x\in(X_{(v)})_i}\hspace{-2mm}H^{i+2}(k(x),\zn(i+1))\rightarrow
\hspace{-2mm}\bigoplus\limits_{x\in(X_{v})_i}\hspace{-2mm}H^{i+2}(k(x),\zn(i)),
$$
which in turn can be written as the sum, over all $x\in X_i$, of
maps
$$
H^{i+2}(k(x),\zn(i+1))\rightarrow
\textstyle
\hspace{-5mm}\bigoplus\limits_{\begin{array}{c}\scriptstyle
x'\in(X_{(v)})_i\vspace{-0,2cm}\\ \scriptstyle x'\mid x
\end{array}}\hspace{-5mm}H^{i+2}(k(x'),\zn(i+1))\rightarrow
\hspace{-5mm}\bigoplus\limits_{\begin{array}{c}\scriptstyle
x''\in(X_{v})_i\vspace{-0,2cm}\\ \scriptstyle x''\mid x
\end{array}}\hspace{-5mm}H^{i+2}(k(x''),\zn(i+1))
$$
By the same reasoning as in the proof of Theorem 3.8, the first
map can be identified with
$$
H^{i+2}(k(x),\zn(i+1))\rightarrow\textstyle\bigoplus\limits_{w\mid v}
H^{i+2}(k(x) K\{x\}_{(w)}, \zn(i+1))
$$
where $K\{w\}$ is the separable closure of $K$ in $k(x)$, which is
a finite extension of $K$, and where $w$ runs over all places of
$K\{x\}$ above $v$. Note that $k(x)$ is a function field of
transcendence degree $i$ over $K\{x\}$. Therefore the restriction
maps above induce an injective map into the direct sum
$$
H^{i+2}(k(x),\zn(i+1))\rightarrow\textstyle\bigoplus\limits_{w\in P(K\{x\})}
H^{i+2}(k(x)_{(w)}, \zn(i+1))
$$
by Proposition 1.2 and Theorem 2.10 (since the latter implies Theorem 0.1, see the introduction).
This shows that we get maps
$$
\alpha_{X,n}:\; C^{2,1}(X,\zn)\mathop{\longrightarrow}\limits^{\beta_{X,n}}\textstyle\bigoplus\limits_v
C^{2,1}(X_{(v)},\zn)\rightarrow\bigoplus\limits_v C^{2,1}(X_v,\zn)\leqno{(4.9)}
$$
of which the first one is injective. The claim of Theorem 4.8 (a)
therefore follows from the next claim.

\vspace{0,5cm}

\noindent \textbf{Proposition 4.10} Let $K$ be a global field, and
let $v$ be a place of $K$. For every variety $V$ over $K_{(v)}$,
every integer $n$ invertible in $K$ and all $r, s\in\mathbb Z$,
the natural map
$$
C^{r,s}(V,\zn)\rightarrow
C^{r,s}(V\times_{K_{(v)}}K_v,\zn)
$$
is injective.

\vspace{0,5cm}\noi\textbf{Proof} In degree $i$, this map is the sum
over all $x\in V_i$ of restriction maps
$$
H^{r+i}(k(x),\zn(s+i))\rightarrow
\textstyle\bigoplus\limits_{\begin{array}{c}\scriptstyle x'\in
\tilde V_i\vspace{-0,2cm}\\ \scriptstyle x'\mid x
\end{array}}H^{r+i}(k(x'),\zn(s+i))
$$
where $\tilde V=V\times_{K_{(v)}} K_v$. For $x\in V_i$, $k(x)$ is
the function field of the integral subscheme (of dimension
$i)\,Z=\overline{\{x\}}\subseteq V$. Because $K_v/K_{(v)}$ is
separable, and $K_{(v)}$ is algebraically closed in $K_v,\;\tilde
Z= Z\times_{K_{(v)}}K_v\hookrightarrow\tilde V$ is a closed
integral subscheme whose generic point $\tilde x$ is in $\tilde
V_i$ and lies above $x$. Let $L$ be the algebraic closure of
$K_{(v)}$ in $k(x)$. Then $\tilde L=L\otimes_{K_{(v)}}K_v$ is a
field, $Z$ is geometrically integral over $L$ with function field
$L(Z)=k(x)$, and $\tilde Z=Z\times_L\tilde L$ with function field
$\tilde L(\tilde Z)= k(\tilde x)$. Moreover, $L$ is henselian,
with completion $\tilde L$. Thus it follows from Theorem 2.11 that
the natural map
$$
H^{r+i}(k(x),\zn(s+i))\rightarrow H^{r+i}(k(\tilde
x),\zn(s+i))
$$
is injective for all $r,s,i\in\mathbb Z$ and all $n\in\mathbb N$
invertible in $K_{(v)}$. This proves Proposition 4.10 and thus Theorem 4.8 (a).

\v

Now we start the proof of Theorem 4.8 (b). The following
rigidity result is shown in [Ja5].

\vspace{0,5cm}

\noindent \textbf{Theorem 4.11} Let $K$ be a global field, and let
$v$ be a place of $K$. For every variety $V$ over $K_{(v)}$, every
integer $n$ invertible in $K$ and all $r,s\in\mathbb Z$, the
natural morphism of complexes
$$
C^{r,s}(V,\zn)\rightarrow
C^{r,s}(V\times_{K_{(v)}}K_v,\zn)
$$
is a quasi-isomorphism, i.e., induces isomorphisms in the
homology.

\vspace{0,5cm}

In view of this result, and of the factorization (4.9), it
suffices to prove Theorem 4.8 (b) after replacing $X_v$ by
$X_{(v)}$ for each $v$. In fact, by Theorem 4.11 we have a canonical quasi-isomorphism
$$
\Cb(X,\mathbb Q_\ell/\mathbb Z) \mathop{\longrightarrow}\limits^{quis}
C'(X,\zn)\,,\leqno(4.12)
$$
where the complex $\Cb(X,\mathbb Q_\ell/\mathbb Z)$
is defined by the exact sequence
$$
0 \rightarrow C^{2,1}(X,\zn)\mathop{\rightarrow}\limits^{\beta_{X,n}}\textstyle\bigoplus\limits_v
C^{2,1}(X_{(v)},\zn)\rightarrow
\Cb(X,\zn)\rightarrow 0\, .\leqno{(4.13)}
$$
So our task is to show $H_0(\Cb(X,\zn))=
\zn$, and $H_a(\Cb(X,\zn))=0$ for $a\neq 0$, if $X$ is connected,
smooth, projective over $K$. Note that all complexes in (4.12) are
concentrated in degrees $0,\ldots, d:=\dim(X)$.

\vspace{0,5cm}\noi
Next we note the following.

\vspace{0,5cm}\noi
\textbf{Lemma 4.14}
For $n\in \mathbb N$ invertible in $K$ the complex
$\overline{C}(X,\zn)$ can be canonically identified with the complex $C^{0,0}(\overline X,\zn)_{G_K}$:
$$
\ldots \rightarrow \mathop{\textstyle\bigoplus}\limits_{x\in
X_r}~H^{r}(K(x)\otimes_K\overline K,\zn(r))_{G_K} \rightarrow
%\mathop{\textstyle\bigoplus}\limits_{x\in X_{r-1}}~H^{r}(K(x)\otimes_K\overline K,r-1)_{G_K}
%\rightarrow
\dots \quad\ldots \rightarrow \mathop{\textstyle\bigoplus}\limits_{x\in
X_0}~H^{0}(K(x)\otimes_K\overline K,\zn(0))_{G_K}\,,
$$
obtained from the Kato complex $C^{0,0}(\overline{X},\zn)$ by taking coinvariants.

\vspace{0,5cm}\noi
\textbf{Proof} This follows easily via the arguments used in the proof of Lemma 3.15,
together with the explicit description of the differentials in this complex in (4.4)
and the covariance of the Hochschild-Serre spectral sequence for corestrictions.

\vspace{0,5cm}\noi
With these tools at hand, we can reduce the proof of Theorem 4.8 (b) to a $\qzl$-version.
Note that it suffices to prove Theorem 4.8 (b) for $n=\ell^m$, for any prime $\ell$
invertible in $K$ and any natural number $m$.
For any  prime $\ell$ and any integers $r, s$ and any scheme $Z$ where it is defined
define the Kato complex $C^{r,s}(Z,\qzl)$ as the direct limit of the complexes
$C^{r,s}(Z,\zln)$ via the transition maps induced by the canonical injections
$\zln \rightarrow \z/\ell^{n+1}$. Then we have in fact:

\vspace{0,5cm}

\noi\textbf{Lemma 4.15} Let $K$ be a global field, let $X$ be a connected smooth proper variety over $K$,
let $\ell$ be any prime. Define the map
$$
\beta_{\ell^\infty}: C^{2,1}(X,\qzl) \longrightarrow \textstyle\bigoplus\limits_v C^{2,1}(X_{(v)},\qzl)
$$
as the inductive limit of the maps $\beta_{X,\ell^m}$ for all $m\in \mathbb N$, and let
$\overline{C}(X,\qzl)$ be its cokernel.

\medskip\noi
(a) The injectivity of $\beta_{X,\ell^\infty}$ is equivalent to the injectivity of the $\beta_{\ell^m}$
for all $m$.

\medskip\noi
(b) To have $H_a(\overline{C}(X,\qzl)) = 0$ for $a\neq 0$ and $H_0(\overline{C}(X,\qzl)=\qzl$ is equivalent to having
$H_a(\overline{C}(X,\zlm))=0$ for all $a\neq 0$ and $H_0(\overline{C}(X,\zlm))=\zlm$ for all $m\in \mathbb N$.

\vspace{0,5cm}\noi
\textbf{Proof} Let $F$ be a function field in $d$ variables over $K$. Then
the theorem of Rost-Voevodsky, more precisely, the validity of the condition BK$(F,d+1,\ell)$ recalled in
the introduction, implies that the sequence
$$
0 \rightarrow H^{d+2}(F,\zln(d+1)) \rightarrow H^{d+2}(F,\mathbb Q/\mathbb Z_\ell(d+1))
\mathop{\rightarrow}\limits^{\ell^n} H^{d+2}(F,\mathbb Q/\mathbb Z_\ell(d+1))
$$
is exact (see the introduction), and the same holds for the fields $F_v$, for all places $v$ of $K$.
By applying this to all residue fields of $X$ and $X_v$, for all $v$, we get a commutative diagram
with exact rows
$$
\begin{CD}
0 @>>> \textstyle\bigoplus_v C^{2,1}(X_{(v)},\zln) @>>> \textstyle\bigoplus_v C^{2,1}(X_{(v)},\qzl) @>{\ell^n}>> \textstyle\bigoplus_v C^{2,1}(X_{(v)},\qzl)\\
@.   @AA{\beta_{X,\ell^n}}A    @AA{\beta_{X,\ell^\infty}}A    @AA{\beta_{X,\ell^\infty}}A  \\
0 @>>>  C^{2,1}(X,\zln)  @>>>  C^{2,1}(X,\qzl) @>{\ell^n}>> C^{2,1}(X,\qzl)\,,
\end{CD}\leqno{(4.16)}
$$
and we deduce the claim in (a).

Now we consider the cokernels of the vertical maps.
First assume that $K$ is a global function field. Then we claim that we even
have an exact sequence
$$
0 \rightarrow H^{d+2}(F,\zln(d+1)) \rightarrow H^{d+2}(F,\mathbb Q/\mathbb Z_\ell(d+1))
\mathop{\rightarrow}\limits^{\ell^n} H^{d+2}(F,\mathbb Q/\mathbb Z_\ell(d+1))\rightarrow 0\,,
$$
similarly for all $F_{(v)}v$. In fact, we have $H^{d+3}(F,\zln(d+1))=0$: If $\ell\neq\ch(F)$,
then $F$ has $\ell$-cohomological dimension $d+2$, and if $\ell=p=\ch(F)$, then we have
$H^{d+3}(F,\zln(d+1)) = H^2(F,W_n\Omega_{F,log})$, but $F$ has $p$-cohomological dimension $1$.
Exactly the same reasoning works for $F_{(v)}$. Writing, for $n$ a positive integer or $n=\infty$,
$$
C_n := \coker[H^{d+2}(F,\zln(d+1)) \longrightarrow \textstyle\bigoplus_v\; H^{d+2}(F_{(v)},\zln(d+1))]\,,
$$
where we set $\z/\ell^\infty\z:= \qzl$, we obtain an exact sequence
$$
0 \rightarrow C_n \rightarrow C_\infty \rightarrow C_\infty \rightarrow 0\,.
$$
Applied to the points of $X$ and the $X_v$ and the morphisms $\beta_{X,\ell^m}$ for $n\in \mathbb N\cup\{\infty\}$,
we now get an exact sequence
$$
0 \rightarrow \overline{C}(X,\zln) \rightarrow \overline{C}(X,\qzl) \mathop{\rightarrow}\limits^{\ell^n} \overline{C}(X,\qzl) \rightarrow 0\,,
$$
and the claim of 4.15 (b) follows in this case.

Now let $K$ be a number field. If $\ell\neq 2$ or if $K$ has no real places, then
$F$ has $\ell$-cohomological dimension $d+2$, and we can argue in the same way.
In general, we can argue in the following way.
In any case, a function field $F$ of transcendence degree $d$ over an
algebraically closed field has $\ell$-cohomological dimension $d$ for $\ell$ invertible in $F$.
It follows that for any variety $X$ over $K$ the sequence
$$
C^{0,0}(\overline{X},\zln) \rightarrow C^{0,0}(\overline{X},\mathbb Q_\ell/\mathbb Z_\ell)
\mathop{\rightarrow}\limits^{\ell^n} C^{0,0}(\overline{X},\mathbb Q_\ell/\mathbb Z_\ell) \rightarrow 0 \leqno{(4.17)}
$$
is exact, where $\overline{X} = X\times_K\overline{K}$ for an algebraic closure $\overline{K}$ of $K$.
Obviously this complex stays exact if we pass to the co-invariants under $G_K$, the absolute Galois group
of $K$.

\smallskip\noi
Now $\beta_{X,\ell^\infty}$ is injective by Theorem 2.10 (Note that $\ell$ is invertible in $K$).
By (4.17) and Lemma 4.14 we get the following commmutative diagram with exact rows and columns
$$
\begin{CD}
@. @. @. 0 \\
@. @. @. @AAA \\
0 @>>> C^{2,1}(X,\mathbb Q_\ell/\mathbb Z_\ell) @>{\beta_{X,\ell^\infty}}>> \textstyle\bigoplus\limits_v
C^{2,1}(X_{(v)},\mathbb Q_\ell/\mathbb Z_\ell) @>>> \Cb(X,\mathbb Q_\ell/\mathbb Z_\ell) @>>> 0 \\
@. @AA{\ell^n}A @AA{\ell^n}A @AA{\ell^n}A \\
0 @>>> C^{2,1}(X,\mathbb Q_\ell/\mathbb Z_\ell) @>{\beta_{X,\ell^\infty}}>> \textstyle\bigoplus\limits_v
C^{2,1}(X_{(v)},\mathbb Q_\ell/\mathbb Z_\ell) @>>> \Cb(X,\mathbb Q_\ell/\mathbb Z_\ell) @>>> 0 \\
@. @AAA @AAA @AA{i}A \\
 @. C^{2,1}(X,\zln) @>{\beta_{X,\ell^n}}>> \textstyle\bigoplus\limits_v
C^{2,1}(X_{(v)},\zln) @>>> \Cb(X,\zln) @>>> 0 \\
@. @AAA @AAA @. \\
@. 0 @. 0
\end{CD}
$$
A simple diagram chase now shows that $\beta_{X,\ell^n}$ and $i$ are injective, which
gives an exact sequence
$$
\begin{CD}
0 @>>> \Cb(X,\zln) @>>>  \Cb(X,\mathbb Q_\ell/\mathbb Z_\ell) @>>> \Cb(X,\mathbb Q_\ell/\mathbb Z_\ell) @>>> 0\,.
\end{CD}
$$
This implies Lemma 4.15 (b).

\v
\noi\textbf{Definition 4.18}
Let $L$ be a perfect field. We say that resolution of singularities holds over $L$, or that ${\text{\bf (RS)}_L}$ holds,
if the following two conditions hold.

%\medskip
%\noindent $\bf (RS)_L:$
\smallskip\noi
${{\text{\bf (RS1)}_L}:}$ For any integral and proper variety $X$ over $L$,
there exists a proper birational morphism $\pi: \tilde{X} \rightarrow X$ such that $\tilde{X}$
is smooth over $L$.

\smallskip
\noi
${{\text{\bf (RS2)}_L}:}$ For any smooth affine variety $U$ over $L$, there is an
open immersion $U \hookrightarrow X$ such that $X$ is projective smooth over $L$ and $Y = X - U$
(with the reduced subscheme structure) is a simple normal crossing divisor on $X$.}

\v\noi
By Hironaka's fundamental results [Hi], resolution of singularities holds over fields $L$ of characteristic zero.

\v\noi
Using the quasi-isomorphism (4.12) as well as Lemmas 4.14 and 4.15,
Theorem 4.8 (b) is obviously implied by the following result.

\v
\noi\textbf{Theorem 4.19} Let $K$ be a finitely generated field with algebraic closure $\overline{K}$
and perfect hull $K^{per}$, let $\ell$ be a prime, and
let $X$ be an irreducible smooth projective variety over $K$.
Assume that $\ell$ is invertible in $K$ or that resolution of singularities over $L=K^{per}$.
Then for the Kato complex $\overline{C}(X,\qzl) := C^{0,0}(\overline{X},\qzl)_{G_K}$ one has
$$
H_a(\overline{C}(X,\qzl)) = \left\{\begin{array}{cl}
\qzl & ,~~a=0\\
0 & ,~~a\neq 0\, .
\end{array}\right.
$$
Here $G_K$, the absolute Galois group of $K$, is regarded as $Gal(\overline{K}/K^{per})$.

\v
\noi
Note that Theorem 4.19 implies Theorem 0.7 in the introduction. By the following lemma,
it also implies Theorem 0.6 in the introduction, concerning Kato's conjecture over finite fields.

\v
\noi\textbf{Lemma 4.20} Let $k$ be a finite field, and let $X$ be any variety over $k$.

\smallskip\noi
(a) One has a canonical isomorphism of complexes
$$
C^{1,0}(X,\zln) \cong C^{0,0}(X\times_k\overline{k},\zln)_{G_k}\,.
$$

\smallskip\noi
(b) The canonical sequence
$$
0 \rightarrow C^{1,0}(X,\zln) \rightarrow C^{1,0}(X,\qzl) \mathop{\rightarrow}\limits^{\ell^n} C^{1,0}(X,\qzl) \rightarrow 0\,.
$$
is exact.

\v

\noi\textbf{Proof} Let $F$ be a function field of transcendence degree $m$ over $k$. Then
one has canonical isomorphisms
$$
H^{m+1}(F,\zln(m)) \cong H^1(k,H^m(F\overline{k},\zln(m))) \cong H^m(F\overline{k},\zln(m))_{G_k}\,,
$$
where $F\overline{k}$ is the function field over $\overline{k}$ deduced
from $F$ (i.e., $F\overline{k} = F\otimes_{\{k\}}\overline{k}$, where
$\{K\}$ is the algebraic closure of $k$ in $F$). In fact, the first isomorphism follows from
the Hochschild-Serre spectral sequence, because $F\overline{k}$ has $\ell$-cohomological dimension $m$,
and the second isomorphism comes from the canonical identification $H^1(k,M) = M_{G_k}$ for
any $G_k$-module $M$, if $k$ is a finite field. By applying this to all fields $k(x)$ for $x\in X$,
we obtain (a).

\v
\noi
(b) follows from the exact sequence
$$
0 \rightarrow H^{m+1}(F,\zln(m)) \rightarrow H^{m+1}(F,\qzl(m)) \mathop{\rightarrow}\limits^{\ell^n} \rightarrow H^{m+1}(F,\qzl(m)) \rightarrow 0\,,
$$
in which the exactness on the left follows from the results of Bloch-Kato-Gabber (for $\ell = p = \ch(k)$)
and Merkurjev-Suslin-Rost-Voevodsky (for $\ell$ invertible in $k$), see the introduction, and the
exactness on the right follows from the cohomological dimension of $F$, compare the proof of Lemma 4.15.

\vspace{0,5cm} The proof of Theorem 4.19 will be given in the next two sections. The idea is to
`localize' the question; but for this we will have to leave the
realm of smooth projective varieties. First recall that the
complexes $C^{r,s}(X,\zn)$ exist for arbitrary varieties $X$
over a field $L$, under the conditions on $X$, $n$ and $(r,s)$ stated at the
beginning of this section.
If $K$ is a global field, the restriction map
$C^{2,1}(X,\mathbb Q_\ell/\mathbb Z_\ell)\rightarrow \Pi_v
C^{2,1}(X_{(v)},\mathbb Q_\ell/\mathbb Z_\ell)$ still has image in
the direct sum and is injective (by Proposition 1.2 and the same argument as
for 4.8 (a)), and we may define $\Cb(X,\mathbb Q_\ell/\mathbb
Z_\ell)$ for arbitrary varieties $X$ over $K$ by exactness of
the sequence
$$
0\rightarrow C^{2,1}(X,\mathbb Q_\ell/\mathbb
Z_\ell)\rightarrow\textstyle\bigoplus\limits_v
C^{2,1}(X_{(v)},\mathbb Q_\ell/\mathbb Z_\ell)\rightarrow
\Cb(X,\mathbb Q_\ell/\mathbb Z_\ell)\rightarrow 0\, .\leqno{(4.21)}
$$
At the same time, by the same arguments as in Lemma 4.14, we have a canonical isomorphism
$$
\overline{C}(X,\mathbb Q_\ell/\mathbb Z_\ell)  \cong  C^{0,0}(\overline{X},\mathbb Q_\ell/\mathbb Z_\ell)_{G_K}
$$
for $\ell$ invertible in $K$.

\vspace{0,5cm} \noindent \textbf{Definition 4.22} Let $L$ be a
field, and let $\mathcal C$ be a category of schemes of finite
type over $L$ such that for each scheme $X$ in $\mathcal C$ also
every closed immersion $i:Y\hookrightarrow X$
and every open immersion $j:U\hookrightarrow X$ is in $\mathcal C$. \\
(a) Let $\mathcal C_\ast$ be the category with the same objects as
$\mathcal C$, but where morphisms are just the proper maps in
$\mathcal C$. A homology theory on $\mathcal C$ is a sequence of
covariant functors
$$
H_a(-):\mathcal C_\ast\rightarrow (\mbox{abelian
groups})\quad\quad (a\in\mathbb Z)
$$
satisfying the following conditions:

\begin{itemize}

\item[(i)] For each open immersion $j:V\hookrightarrow X$ in
$\mathcal C$, there is a map $j^\ast:H_a(X)\rightarrow H_a(V)$,
associated to $j$ in a functorial way.

\item[(ii)] If $i:Y\hookrightarrow X$ is a closed immersion in
$\mathcal C$, with open complement $j:V\hookrightarrow X$, there
is a long exact sequence (called localization sequence)
$$
\ldots \stackrel{\delta}{\longrightarrow}
H_a(Y)\stackrel{i_\ast}{\longrightarrow}
H_a(X)\stackrel{j^\ast}{\longrightarrow}
H_a(V)\stackrel{\delta}{\longrightarrow}H_{a-1}(Y)\longrightarrow\ldots
$$
(The maps $\delta$ are called the connecting morphisms.) This
sequence is functorial with respect to proper maps or open
immersions, in an obvious way.

\end{itemize}

\medskip
\noindent (b) A morphism between homology theories $H$ and $H'$ is
a morphism $\phi:H\rightarrow H'$ of functors on $\mathcal
C_\ast$, which is compatible with the long exact sequences from
(ii).

\vspace{0,5cm} \noindent \textbf{Lemma 4.23} (a) Let $L$ be a
field, and let $r, s,$ and $n\geq 1$ be fixed integers with $n$ invertible in
$L$, or $r\neq s+1$, or $p = \ch(L)\mid n$ and $r=s+1$ and $[L:L^p]\leq p^s$.
There is a natural way to extend the assignments
$$
X \quad \raisebox{1pt}{$\shortmid$}\hspace{-5pt}\rightsquigarrow \quad
H^{r,s}_a(X,\zn):= H_a(C^{r,s}(X,\zn)) \quad\quad  (a\in\mathbb Z)
$$
to a homology theory on the category of all varieties over $L$.

\medskip\noindent
(b) The same holds for the assignment
$$
X \quad \raisebox{1pt}{$\shortmid$}\hspace{-5pt}\rightsquigarrow \quad
\Hb^{r,s}_a(X,\zn):= H_a(\Cb^{r,s}(X,\zn)) \quad \quad (a\in\mathbb Z)\, ,
$$
where $\Cb^{r,s}(X,\zn) := C^{r,s}(\overline{X},\zn)_{G_L}$, with
$\overline{X} = X\times_L\overline{L}$ for a separable closure of $L$.

\vspace{0,5cm} \noindent \textbf{Proof} (a): The Bloch-Ogus-Kato
complexes are covariant with respect to proper morphisms and
contravariant with respect to open immersions. The localization
sequence for a closed immersion $i:Y\hookrightarrow X$ with open
complement $j:U=X\smallsetminus Y\hookrightarrow X$ is obtained by
the short exact sequence of complexes
$$
0\rightarrow C^{r,s}(Y,\zn)\stackrel{i_\ast}{\rightarrow}
C^{r,s}(X,\zn)\stackrel{j^\ast}{\rightarrow}
C^{r,s}(U,\zn)\rightarrow 0
$$
which are componentwise canonically split.(cf. also [JS] Corollary 2.10).

\medskip\noindent
(b): This follows from (a), because the mentioned splitting is equvariant, so that
the sequences stay exact after taking coinvariants.

\vspace{0,5cm} The mentioned localization is now obtained by the
following observation.

\vspace{0,5cm} \noindent \textbf{Lemma 4.24} Let $L$ be a perfect
field, let $\mathcal C$ be a category of schemes of finite type over $L$ as in 4.20,
and let $\varphi: H\rightarrow \widetilde H$ be a morphism of
homology theories on the category ${\mathcal C}_\ast$ of all
schemes in $\mathcal C$ with proper morphisms. For every integral variety $Z$ over $L$ let
$L(Z)$ be its function field. Define
$$
H_a(L(Z)):= \lim\limits_{\longrightarrow} H_a(U)\, ,
$$
where the limit is over all non-empty open subvarieties $U$ of
$Z$, and define $\widetilde H_a(L(Z))$ similarly. Suppose the
following holds for every integral variety $Z$ of dimension $d$ over $L$.
\begin{itemize}

\item[(i)] $H_a(L(Z))=0$ for $a\neq d$,

\item[(ii)] $\widetilde H_a(L(Z))=0$ for $a\neq d$, and

\item[(iii)] the map $\varphi:H_d(L(Z))\rightarrow \widetilde
H_d(L(Z))$ induced by $\varphi$ is an isomorphism.
\end{itemize}

\noindent Then $\varphi$ is an isomorphism of homology theories.

\vspace{0,5cm} Before we give a proof for this, we note the
following.

\vspace{0,5cm} \noindent \textbf{Remark 4.25} The homology theories
of 4.23 (a) clearly satisfy condition 4.24 (i). In fact, setting
$$
C^{r,s}(L(X),\zn) = \lim\limits_{\longrightarrow} C^{r,s}(U,\zn)
$$
for an integral $X$, where the limit is over all non-empty open subvarieties $U$ of
$X$, we trivially have $C^{r,s}_a(L(X),\zn)=0$ for $a\neq \dim X$, because for any
$x\in X$ different from the generic point there is an non-trivial open $U\subset X$
not containing $x$.
Hence 4.24 (i) also hold for the homology theories from 4.23 (b).
The proof of Theorem 4.19 will then achieved as follows. In the next section
we will define a homology theory $H^W_\ast(-,\qzl)$, over any field $K$ of
characteristic 0, or over a perfect field of positive characteristic assuming
suitable resolution of singularities,
which a priori satisfies
$$
H^W_a (X,\qzl)= \left\{\begin{array}{cl}
\qzl & ,~~a=0\\
0 & ,~~a\neq 0\, .
\end{array}
\right. \leqno{(4.26)}
$$
if $X$ is smooth, projective and irreducible. Moreover we will
show that 4.22 (ii) holds for $H^W$. Still under the same assumptions
we will construct a morphism
$$\varphi: \Hb_\ast(-,\mathbb Q_\ell/\mathbb Z_\ell) := \Hb^{0,0}_\ast(-,\mathbb Q_\ell/\mathbb Z_\ell)\rightarrow H^W(-,\qzl)$$
of homology theories which satisfies 4.24 (iii) if $K$ is finitely generated .
Thus, by Lemma 4.24, $\varphi$ is an isomorphism, and hence
(4.26) also holds for $\Hb_\ast(-,\mathbb Q_\ell/\mathbb Z_\ell)$,
which proves 4.19.

\vspace{0,5cm}\noindent\textbf{Proof of Lemma 4.24} For every
homology theory $H$ over $L$ there is a strongly converging niveau
spectral sequence
$$
E^1_{p,q}(X)=\textstyle\bigoplus\limits_{x\in X_p}
H_{p+q}(k(x))\Rightarrow H_{p+q} (X)\leqno{(4.27)}
$$
for every $X$, cf. [BO], and also [JS]. If $\widetilde E^1_{p,q}
\Rightarrow \widetilde H_{p+q}$ is associated to another homology
theory $\widetilde H$, then every morphism $\varphi:H\rightarrow\widetilde
H$ induces a morphism $E\rightarrow\widetilde E$ of these spectral
sequences, compatible with $\varphi$ on the $E^1$-terms and limit
terms. In the situation of 4.24, the conditions (i), (ii) and
(iii) imply that $\varphi$ induces isomorphisms on the
$E^1$-terms, and hence $\varphi$ also gives an isomorphism between
the limit terms, i.e., between $H$ and $\widetilde H$.

\vspace{15mm}

\centerline{\noindent\textbf{\S 5 Weight complexes and weight cohomology}}

\vspace{1,0cm}

Let $k$ be a field. Let $X$ be a smooth, projective variety of
dimension $d$ over $k$, and let $Y=\bigcup\limits^r_{i=1} Y_i$ be
a divisor with simple normal crossings in $X$ -- with a fixed
ordering of the smooth components as indicated.

\vspace{0,5cm} \noindent \textbf{Definition 5.1} Let $F$ be a
covariant functor on the category $\mathcal{SP}_k$ of smooth
projective varieties with values in an abelian category $\mathcal
A$ which is additive in the sense that the natural arrow
$$
F(X_1)\textstyle\scriptstyle \bigoplus \displaystyle F(X_2)
\rightarrow F(X_1 \scriptstyle\coprod\displaystyle X_2)
$$
is an isomorphism in $\mathcal A$, where $X_1\scriptstyle\coprod
\textstyle X_2$ is the sum (disjoint union) of two varieties $X_1,
X_2$ in $\mathcal{SP}_k$. Then define $L^i F(X,Y)$ as the $i$-th
homology of the complex

\medskip\noindent
$$
C.F(X,Y): ~~~~~~~~~ 0\rightarrow F(Y^{[d]})\rightarrow
F(Y^{[d-1]})\rightarrow\ldots \rightarrow F(Y^{[1]})\rightarrow
F(X) \rightarrow 0\, .
$$
Here $F(Y^{[j]})$ is placed in degree $j$, and the differential
$\partial:F(Y^{[j]})\rightarrow F(Y^{[j-1]})$ is
$\sum\limits^j_{\nu= 1}(-1)^\nu\delta_\nu$ where $\delta_\nu$ is
induced by the inclusions
$$
Y_{i_1,\ldots, i_j} \hookrightarrow Y_{i_1,\ldots, \hat i_\nu,
\ldots, i_j}
$$
for $1\leq i_1< \ldots <i_j\leq r$ (and where $Y^{[0]} =X$, as
usual).

\vspace{0,5cm} \noindent \textbf{Remark 5.2.} There is the dual
notion of an additive, contravariant  functor $G$ from $\mathcal{
SP}_k$ to $\mathcal A$, and here we define $R^i G(X,Y)$ to be the
$i$-th cohomology of the complex
$$
C^\cdot G(X,Y): ~~~~~~~~~ G(X)\rightarrow G(Y^{[1]})\rightarrow
\ldots\rightarrow G(Y^{[d-1]})\rightarrow G(Y^{[d]})\, ,
$$
with $G(Y^{[j]})$ is placed in degree $j$.

\vspace{0,5cm} We may apply this to the following functors. Let
$Ab$ be the category of abelian groups.

\vspace{0,5cm} \noindent \textbf{Definition 5.3} For any abelian
group $A$ define the covariant functor $H_0(-,A):\mathcal{
SP}_k\rightarrow Ab$ and the covariant functor $H^0(-,A):\mathcal{
SP}_k\rightarrow Ab$ by
$$
\begin{array}{rcl}
H_0(X,A) & = & \textstyle\bigoplus\limits_{\alpha\in\pi_0(X)}
A = A\otimes_\z\z[\pi_0(X)]\\\\
H^0(X,A) & = & A^{\pi_0(X)} = \mbox{ Map }(\pi_0(X),A)\, .
\end{array}
$$\\
where $\z[M]$ is the free abelian group on a set $M$, and
$\mbox{Map}(M,N)$ is the set of maps between two sets $M$ and $N$.
(Hence if $A$ happens to be a ring, then $H_0(X,A)$ is the free
$A$-module on $\pi_0(X)$, and $H_0(X,A)=$ Hom$_A(H_0(X,A),A)$ is
its $A$-dual.) We write $C^W_\cdot(X,Y;A)$ for $C_\cdot
H_0(-,A)(X,Y)$ and call
$$ H^W_i(X,Y;A) := L^iH_0(-,A)(X,Y)=H_i(C_\cdot(X,Y;A))$$
the weight homology of $(X,Y)$. Similarly define $C^\cdot_W(X,Y;A)=
C^\cdot H_0(-,A)(X,Y)$ and call $H^i_W(X,Y;A)=H^i(C^\cdot_W(X,Y;A))$
the weight cohomology of $(X,Y)$.

\vspace{0,5cm} \noindent \textbf{Proposition 5.4} Let $Y_{r+1}$ be
a smooth divisor on $X$ such that the intersections with the
connected components of $Y^{[j]}$ are transversal for all $j$ and
connected for all $j\leq d-2$. (Note: If $k$ is infinite, then by
the Bertini theorems such a $Y_{r+1}$ exists by taking a suitable
hyperplane section, since $\dim Y^{[j]}= d-j\geq 2$ for $j\leq
d-2$.) Let $Z=\bigcup\limits_{i=1}^{r+1}Y_i$ (which, by the
assumption, is again a divisor with normal crossings on $X$).
Then, for any abelian group $A$
$$
H^W_i(X,Z;A)=0=H^i_W(X,Z;A)\mbox{ for } i\leq d-1\, .
$$

\vspace{0,5cm} \noindent \textbf{Proof} Fix $A$ and omit it in the
notations. We get a commutative diagram
$$
\xymatrix{0\ar[r] & 0 \ar[r]\ar[d]& H_0(Y^{[d-1]}\cap Y_{r+1})
\ar[r]^{\delta_{d-1}}\ar@{->>}[d]^{\psi_{d-1}}  & H_0(Y^{[d-2]}\cap Y_{r+1})\ar[r]\ar[d]^{\psi_{d-2}}_\wr & \ldots   \\
0\ar[r] & H_0( Y^{[d]})\ar[r]^{\partial_d} & H_0( Y^{[d-1]})\ar[r]^{\partial_{d-1}}& H_0( Y^{[d-2]})\ar[r] & \ldots  \\
& \ldots\ar[r]& H_0(Y^{[1]}\cap Y_{r+1})\ar[r]\ar[d]^{\psi_{1}}_\wr & H_0( Y_{r+1})\ar[r]\ar[d]^{\psi_{0}}_\wr & 0\\
&  \ldots\ar[r]& H_0( Y^{[1]})\ar[r]&H_0(X)\ar[r] & 0}
$$
where the bottom line is the complex $C^W_\cdot(X,Y)$ and the top line
is $C^W_\cdot (Y_{r+1}, Y\cap Y_{r+1})$ (note that
$Y\cap Y_{r+1}= \bigcup\limits^{r}_{i=1} (Y_i\cap Y_{r+1})$ is a
divisor with strict normal crossings on the smooth, projective
variety $Y_{r+1}$), and where $\psi_\nu$ is induced by the
inclusion $Y^{[\nu]}\cap Y_{r+1}\hookrightarrow Y^{[\nu]}$.

\vspace{0,5cm} \noindent By the assumption, $\psi_\nu$ is an
isomorphism for $\nu\leq d-2$, and a (non-canonically) split
surjection for $\nu=d-1$. Hence we have isomorphisms
$$
\begin{array}{lcl}
\mbox{isomorphisms}\;  & H^W_i(Y_{r+1}, Y\cap Y_{r+1})\stackrel{\sim}{\longrightarrow}
H^W_i(X,Y) & \mbox{ for } i\leq d-2\,, \quad \mbox{and}\\\\
\mbox{a surjection}\;\quad    &  H^W_{d-1}(Y_{r+1}, Y\cap Y_{r+1})\twoheadrightarrow H^W_{d-1}(X,Y)\,. &
\end{array}
$$
Moreover, let $C_{\cdot\cdot}$ be the associated double complex,
with $H_0(X)$ placed in degree $(0,0)$ and $\psi_\nu$ being
replaced by $(-1)^{\nu}\psi_\nu$. Then the associated total
complex $s(C_{\cdot\cdot})$ is just the complex $C^W_\cdot (X,Z)$.
Hence the result follows, and we have exact sequences
$$
0\rightarrow H^W_d(X,Y) \rightarrow H^W_d(X,Z) \rightarrow H^W_{d-1}(Y_{r+1}, Y\cap
Y_{r+1})\rightarrow H^W_{d-1}(X,Y)\rightarrow 0\, ,
$$
$$
0 \rightarrow \ker(\psi)_{d-1} \rightarrow H^W_d(X,Z) \rightarrow
H^W_0(Y^{[d]}) \rightarrow 0 \, .
$$
The proof for $H^i_W (X,Z)$ is dual. (Note, however, that in general
$H^W_\cdot (X,Y)$ and $H^\cdot_W(X,Y)$ are related by a coefficient
theorem in an non-trivial way).

\vspace{0,5cm} \noindent \textbf{Corollary 5.5} $H^W_d(X,Z;\mathbb
Z)$ is a finitely generated free $\mathbb Z$-module, and we have an isomorphism
$H^W_d(X,Z;A)\cong H_d(X,Z;\,\mathbb Z)\otimes_\mathbb Z A$. The same
holds for $H^d_W(X,Z)$.

\vspace{0,5cm} \noindent \textbf{Proof} The first statement
follows since $\ker(\psi)_{d-1}$ and $H^0(Y^{[d]})$ have this
property for $A=\mathbb Z$, and the second claim follows from 5.4
and the universal coefficient theorem. Similarly for $H^d_W(X,Z)$.

\vspace{0,5cm} \noindent \textbf{Corollary 5.6} $H^W_d(X,Z;\mathbb
Q_\ell/\mathbb Z_\ell)$ is divisible.

\vspace{0,5cm} Now let $U=X\smallsetminus Y$.

\vspace{0,5cm} \noindent \textbf{Proposition 5.7} Let $\overline{U}=U\times_k\overline{k}$ where
$\overline{k}$ is the algebraic closure of $k$. Then there are
canonical homomorphisms
$$
H^d_\et(\overline U,\zn(d))\stackrel{e}{\longrightarrow}
H^W_d(X,Y;\zn)
$$
for all $n\in\mathbb N$. If $k$ is finitely generated, and
if  $X\smallsetminus Y_\nu$ is affine for one $\nu\in\{1,\ldots,
r\}$, then these induces isomorphisms
$$
H^d_\et(\overline U,\mathbb Q_\ell/\mathbb
Z_\ell(d))_{G_k}\stackrel{\sim}{\longrightarrow} H^W_d(X,Y;\mathbb
Q_\ell/\mathbb Z_\ell)
$$
for all primes $\ell$.

\vspace{0,5cm} This is just a reformulation of Theorem 3.1, in which
the construction of $e$ does not depend on the assumption that
some $X\smallsetminus Y_\nu$ is affine.

\vspace{0,5cm} \noindent \textbf{Remark 5.8} In particular, with
the notations and assumptions of 5.4, this applies to $(X,Z)$ and
$U=X\smallsetminus Z$.

\vspace{0,5cm} We want to have these results in a more functorial
setting. This is possible if resolution of singularities holds,
in a suitable form.

\vspace{0,5cm}

\noindent \textbf{Theorem 5.9} Let $k$ be a field, with perfect hull $L=k^{per}$,
let $A$ be an abelian group, and assume that one of the following two conditions is fulfilled.

\medskip\noi (a) Condition ${{\text{\bf (RS1)}_L}}$ from 4.18 holds.

\medskip\noi (b) The group $A$ is a $\z_{(\ell)}$-module, where $\ell$ is a prime invertible in $k$
and $\z_{(\ell)}$ is the localization of $\z$ with respect to the prime ideal $(\ell)$.

\medskip\noi
Then there exists a homology theory
(in the sense of definition 4.22) $\,(H^W_a(-,A), a\in \z)\,$ on the
category $({\mathcal V}_k)_\ast$ of all varieties over $k$ with proper morphisms such that the following holds.

\smallskip\noi (i) For any smooth, projective and connected variety $X$ over $k$ one has
$$
H^W_a(X,A)=\left\{\begin{array}{lcl} 0 & & a\neq 0\\ A & &
a=0\; .\end{array}\right.
$$

\noi
(ii) If $X$ is smooth projective over $k$ and $Y$ is a divisor with simple normal crossings on $X$,
then one has a canonical isomorphism for $U = X\setminus Y$
$$
H^W_a(U,A) \cong H^W_a(X,Y;A)
= H_a(\textstyle\bigoplus\limits_{\pi_0(Y^{[e]})} A \rightarrow \ldots \rightarrow
\textstyle\bigoplus\limits_{\pi_0(Y^{[1]})} A \rightarrow \textstyle\bigoplus\limits_{\pi_0(X)} A)
$$
where the right hand side is defined in Definition 5.3.

\vspace{0,5cm}

\noindent \textbf{Proof} First assume that $k$ is perfect.
We want to show that the covariant functor (cf. 5.3)
$$
F: ~ X \quad \funct \quad H_0(X,A) = \bigoplus\limits_{\pi_0(X)} A
$$
on the category ${\mathcal {SP}}_k$ of all smooth projective varieties over $k$ extends to a homology
theory on all of ${\mathcal V}_k$, and fulfills (ii).
By the method of Gillet and Soul\'e (\cite{GS} proof of 3.1.1) this holds if $(RS1)_k$
holds and if $F$ extends to a contravariant functor on
Chow motives over $k$, i.e., admits an action of algebraic
correspondences modulo rational equivalence. But the latter is clear --
in fact, one has $F(X) = \Hom(CH^0(X),A)$, where $CH^j(X)$ is the
Chow group of algebraic cycles of codimension $j$ on $X$, modulo
rational equivalence. Under condition (b), the claim follows from
Theorem 5.13 below.

\medskip\noi
If $k$ is general, we just define
$$
H^W_a(Z,A) := H^W_a(Z\times_kk^{per})
$$
where the theory on the right is the one existing over $k^{per}$ by our
assumptions and the case of a perfect field. Note that $Z\times_kk^{per}$ is
again connected for connected $Z$.

\vspace{0,5cm}

\noindent\textbf{Theorem 5.10} Let $k$ be a field and assume that resolution of
singularities holds over the perfect hull $L$ of $k$ (see Definition 4.18). Then the homology
theory $H^W_\ast(-,A)$ of Theorem 5.9 has the property 4.22 (ii).

\vspace{0,5cm}

\noindent\textbf{Proof:} By construction we may assume $k$ is perfect.
For every integral variety $Z$ over $k$ we have to show
$$
H^W_a(k(Z),A) := \lim\limits_{\longrightarrow}
H^W_a(V,A) = 0 \mbox{ for } a \neq e\,,\leqno{(5.11)}
$$
where the inductive limit is over all non-empty open subvarieties
$V \subset Z$.

\smallskip
Now assume $(RS)_k$ holds. Then, by perfectness of $k$, for every
non-empty open subvariety $V\subset Z$ there is an non-empty
smooth open subvariety $U \subset V$ and an open embedding $U
\hookrightarrow X$ into a smooth projective variety $X$ such that
the complement $Y = X \setminus U$ is a divisor with strict normal
crossings.
If $k$ is infinite, then, by Bertini's theorem, there exists a
smooth hyperplane section $H$ of $X$ whose intersection with all
connected components of $Y^{[i]}$ is smooth, and connected for $i
\leq d-2$. Writing $Z = Y \cup H$ (which a divisor with strict
normal crossings on $X$) and $U^0 = X \setminus Z \subset U
\subset V$, we get $H^W_a(U^0,A) = H^W_a(X,Z;A) = 0$ for $a
\neq e$ by Property 5.9 (ii) and Proposition 5.4. Since $V$ was arbitrary, we get
property (5.11).

If $k$ is finite, we use a suitable norm argument. By what has
been shown, for each prime $p$ we find such a hyperplane section
after base change to an extension $k'/k$ of degree $[k':k] = p^r$,
a power of $p$ (the maximal pro-$p$-extension of $k$ is an
infinite field). Then the map
$$
H^W_a(V_{k'},A) \rightarrow H^W_a(k'(Z_{k'}),A)
$$
is zero.

Now we note that there is a homology theory $H^W(-,A;k')$ on
all varieties over $k$, defined by $H^W_a(Z,A;k') = H^W_a(Z_{k'},A)$
and the induced structure maps. This is also the
homology theory which is obtained by the method of Theorem 5.9, by
extending the covariant functor
$$
F': \mathcal{SP}_k \longrightarrow Ab~,~ X \;\funct\;
\bigoplus\limits_{\pi_0(X_{k'})} A
$$
to a homology theory on all varieties. There is an obvious morphism of functors
$Tr: F' \rightarrow F$ (trace, or norm), induced by the natural maps $\pi_0(X_{k'}) \rightarrow \pi_0(X)$.
On the other hand there is also a morphism of functors $Res: F \rightarrow F'$ (restriction)
such that $Tr Res = [k':k]$.

This is best seen by noting that for any smooth
projective variety $X/k$ one has a canonical isomorphism
$$
F(X) = (\bigoplus\limits_{\pi_0(X\times_k\overline{k})} A)_{G_k} \stackrel{\sim}\longrightarrow \bigoplus\limits_{\pi_0(X)} A\,,
$$
where $\overline{k}$ is an algebraic closure of $k$.
Similarly, $F'(X) = B_{G_{k'}}$, where
$$
B = \bigoplus\limits_{\pi_0((X\times_k k')\times_{k'}\overline{k})} A = \bigoplus\limits_{\pi_0(X\times_k\overline{k})} A\,.
$$
In these terms, $Tr$ is induced by the natural map $B_{G_{k'}} \rightarrow B_{G_k}$.
Conversely, for any profinite group $G$, any open subgroup $U$ and any discrete $G$-module $C$ we have a
well-defined functorial map
$$
Cor^\vee(C):\; C_G \longrightarrow C_H \quad \quad , \quad \mbox{class of}\quad a \quad \mapsto \quad \mbox{class of}\quad \sum_{\sigma \in G/H} \sigma a \;,
$$
and the composition $Cor^\vee \circ \pi$ is the multiplication by $(G:H)$. Applied to $(G,H,C) = (G_k,G_{k'},B)$ we
get the claim.

By the construction of Gillet and Soul\'e (or by Theorem 5.13 and Remark 5.15 below), $Tr$
extends to morphism of homology theories $Tr: H^W(-,A;k') \rightarrow H^W(-,A)$, and one checks that the induced maps
$H^W_a(Z_{k'},A) \rightarrow H^W_a(Z,A)$ are just the maps
obtained from functoriality for proper morphisms. Similarly, $Res$ extends to a morphism of homology theories
$Res: H(-,A) \rightarrow H^W(-,A;k')$, and one has $Tr Res = [k':k] = p^r$,
because this holds for the restriction to the functors $F$ and
$F'$. The outcome is that the kernel of $Res$ is killed by $p^r$.
From the commutative diagram
$$
\begin{array}{ccc}
H^W_a(V_{k'},A) & \rightarrow & H^W_a(k'(Z_{k'}),A) \\
\uparrow && \uparrow \\
H^W_a(V,A) & \rightarrow & H^W_a(k(Z),A)\,,
\end{array}
$$
we then get that the image of $H^W_a(V,A)$ in $H^W_a(k(Z),A)$ is killed
by $p^r$ for $a\neq e$, because its image in $H^W_a(k'(Z_{k'}),A)$ is
zero. Since this holds for all $V$ (with varying powers of $p$), we
conclude that every element in the group $H^W_a(k(Z),A)$ is
killed by a power of $p$. Since this also holds for any second
prime $q\neq p$, we conclude the  vanishing of $H^W_a(k(Z),A)$.

\vspace{0,5cm}

\noi\textbf{Theorem 5.12} Let $\ell$ be a prime, and let $K$ be a finitely generated field such that
one of the following conditions is fulfilled.

\medskip\noi (A) $(RS)_{K^{per}}$ holds.

\medskip\noi (B) $\ell$ is invertible in $K$.

\medskip\noi
Let $\Hb_\ast(-,\qzl)$ and $H^W_\ast(-,\qzl)$ be the homology theories
on ${\mathcal V}_{K,d}$ defined in 4.21 (b)
and 5.9, respectively. There exists a morphism
$$
\varphi: \Hb_\ast(-,\qzl) \longrightarrow H^W_\ast(-,\qzl)
$$
of homology theories such that the properties (i) - (iii) of Lemma 4.22
are fulfilled for $\Hb$, $H^W$ and $\varphi$. Consequently, $\varphi$ is an isomorphism
of homology theories.

\vspace{0,5cm}

Evidently this theorem implies Theorem 4.18, in view of 5.9 (i). Theorem 5.12 will be proved
in the next section. Here we will start with some preparations.

\vspace{0,5cm}
Since $H^W_\ast(-,\qzl)$ is defined via the method of Gillet and Soul\'e
in \cite{GS} 3.1.1, we need to analyze the constructions in \cite{GS}
more closely. There functors on Chow motives with values in
abelian categories are extended to homology theories on all
varieties. We give a more general version for complexes, in the following form.

\vspace{0,5cm}

\noi\textbf{Theorem 5.13} Let $k$ be a perfect field, let
$\mathcal C_{\geq 0}(\mathcal A)$ be the category of
non-negative homological complexes $\ldots \rightarrow C_2
\rightarrow C_1 \rightarrow C_0$ in an abelian category $\mathcal
A$, and let
$$
C: \mathcal{SP}_k \rightarrow \mathcal C_{\geq 0}(\mathcal A)
$$
be a covariant functor on the category $\mathcal{SP}_k$ of all
smooth projective varieties over $k$. Assume that the associated
functors
$$
H^C_a: \mathcal{SP}_k \rightarrow \mathcal A ~,~ X \funct H_a(C(X))
$$
extend to contravariant functors on the category $\mathcal{CHM}^{eff}(k,d)$
of effective Chow motives (i.e., $\z$-linear motives modulo
rational equivalence) generated by all smooth projective varieties over $k$.
Assume further that one of the following conditions holds.

\medskip\noi
(a) Condition $(RS)_k$ holds.

\smallskip\noi
(b) $\mathcal A$ is a $\z_{(\ell)}$-linear category, where $\ell$ is a prime invertible in $k$.

\medskip\noi
Then there is a natural way to
extend the above functors $H^C_a$ to a homology theory on the
category ${\mathcal V}_k$ of all varieties over $k$.

\v\noi
The proof will treat both cases in a parallel way and use the following tools and lemmas.

\v\noi
\textbf{Definition 5.15} (a) A morphism $f: X \rightarrow Y$ of schemes is called an envelope,
if for any point $y\in Y$ there is a point $x\in X$ such that $f(x)=y$ and $k(y)\rightarrow k(x)$ is an isomorphism.

\smallskip\noi
(b) Let $\ell$ be a prime. A morphism $f: X \rightarrow Y$ of schemes is called an $\ell'$-envelope
if for any point $y\in Y$ there is a point $x\in X$ such that $f(x)=y$ and $k(x)/k(y)$ is a finite
extension of degree prime to $\ell$.

\v\noi
Our definitions differ from those in \cite{GS} and \cite{GS2} in that we don't assume at this point that $f$ is proper.
Note that any envelope is an $\ell'$-envelope for all primes $\ell$. We treat both cases in a parallel
way by allowing also $\ell=\infty$, in which case an $\infty'$-envelope is just an envelope.

\v\noi
\textbf{Definition 5.16}
Call an envelope or $\ell'$-envelope $f: X \rightarrow Y$ small, if for any generic point
$\eta \in Y$ the fibre $X\times_Y\eta$ is a finite $k(\eta)$-scheme.

\v\noi
\textbf{Definition 5.17}
Let $X.$ be a simplicial variety over $k$, and let $Y$ be a variety over $k$. A morphism $X. \rightarrow Y$
is called a (small) $\ell'$-hyperenvelope, if for all $n\geq 0$ the natural map $X_n \rightarrow (cosk^Y_{n-1}(sk_{n-1}(X.)))_n$
is a (small) $\ell'$-envelope.

\v\noi
\textbf{Lemma 5.18}
(a) If $(RS)_k$ holds, then any variety $Y$ in $\cV_k$ possesses a smooth small hyperenvelope $X. \rightarrow Y$,
i.e., one where all $X_i$ are smooth. If $Y$ is proper, then there exists a smooth projective small
hyperenvelope, i.e., one where the $X_i$ are projective in addition.

(b) If $\ell$ is a prime which is invertible in $k$, then any variety $Y$ in $\cV_{k,d}$ possesses a smooth small
$\ell'$-hyperenvelope $X. \rightarrow Y$, and a smooth projective small $\ell'$-hyperenvelope if $Y$ is proper.

\v\noi\textbf{Proof:}
Under $(RS)_k$, any $Y$ possesses a small envelope $X \rightarrow  Y$ with smooth $X$. If $\ell$ is invertible
in $k$, then by the result of Gabber-Illusie (property (G) in the proof of 2.10) $Y$ possesses a small $\ell'$-envelope
$X \rightarrow Y$ with $X$ smooth. Therefore the claim follows from Lemma 7.8 below.

\v\noi
\textbf{Lemma 5.19}
Let $A$ be any commutative ring, let $\mathcal{SP}=\mathcal{SP}_k$ be the category of smooth projective varieties over $k$,
and let $A\mathcal{SP}$ be the category with the same objects as $\mathcal{SP}$, but with morphism sets the free $A$-modules
on $Hom_{\mathcal{SP}}(X,Y)$ (the $A$-analogue of the category $\z\mathbb{V}$ in \cite{GS} 1.2).
Finally let $\cM=\mathcal{CHM}^{eff}(k,A)$ be the category of effective Chow motives over $k$
with coefficients in $A$, and let $R_{q,\ast}(X)$ the Gersten complex of algebraic $K$-groups of a variety $X$ over $k$
defined in \cite{GS} 1.1. If
$$
\ldots \rightarrow X_2 \rightarrow X_1 \rightarrow X_0
$$
is a chain complex in $A\mathcal{SP}$ such that for all $q\geq 0$ and all $V \in \cSP$ the total complex
of the induced double complex
$$
\ldots \rightarrow R_{q,\ast}(V\times X_2) \rightarrow R_{q,\ast}(V\times X_1) \rightarrow R_{q,\ast}(V\times X_0)
$$
becomes acyclic after tensoring with $A$. Then the complex in $\cM$
$$
\ldots \rightarrow M(X_2) \rightarrow M(X_1) \rightarrow M(X_0)
$$
has a contracting homotopy.

\v\noi
The proof is literally the same as the proof of Theorem 1 in \cite{GS}, which treats the case $A=\z$.
The following is the $\ell'$-analogue of \cite{GS} Proposition 1 (in a simpler form).

\v\noi
\textbf{Lemma 5.20}
If $Z. \rightarrow X$ is a proper $\ell'$-hyperenvelope of varieties over $k$, then, for all $q\geq 0$,
the associated morphism of complexes
$$
tR_{q,\ast}(Z.) \rightarrow R_{q,\ast}(X)
$$
becomes an isomorphism after tensoring with $\z_{(\ell)}$, where $tR_{q,\ast}(Z.)$ is the total
complex associated to the natural double complex associated to the simplicial complex $R_{q,\ast}(Z.)$.
Here we set $\z_{(\ell)}=\z$ for $\ell=\infty$.

\v\noi
This follows from the following more general descent lemma.

\v\noi
\noi\textbf{Lemma 5.21} Let the assumption be as in Theorem 5.13,
and assume that
$$
C: ({\mathcal V}_k)_\ast \longrightarrow \mathcal C_{\geq 0}(\mathcal A)
$$
is equipped with the following additional data:
\begin{itemize}
\item[(i)] For every open immersion $j: U \hookrightarrow Z$ in
$(\mathcal V_k)_\ast$ there is a morphism $j^\ast: C(Z)
\rightarrow C(U)$, associated to $j$ in a functorial way.

\item[(ii)] If $i: Y \hookrightarrow Z$ is a closed immersion in
$(\mathcal V_k)_\ast$, with open complement $j: U \hookrightarrow
Z$, then there is a short exact sequence of complexes
$$
0 \rightarrow C(Y) \mathop{\rightarrow}\limits^{i_\ast} C(X)
\mathop{\rightarrow}\limits^{j^\ast} C(U) \rightarrow 0\,.
$$
This sequence is functorial with respect to proper morphisms and
open immersions, in an obvious way.
\end{itemize}
Then the following holds:

\smallskip\noi
(a) If $Z_. \rightarrow X$ is a hyperenvelope
of a variety in ${\mathcal V}_k$, then the canonical morphism
$ tC(Z_.) \longrightarrow C(X) $
(induced by the morphism $Z_0 \rightarrow X$) is a quasi-isomorphism. Here
$tC(Z_.)$ is defined as in Lemma 5.20.

\smallskip\noi
(b) Assume that the conditions of Theorem 5.13 (b) hold, for the prime $\ell$.
Assume moreover the following conditions:
\begin{itemize}
\item[(iii)]
The functor $C$ is also contravariant for finite flat morphisms $\pi: X \rightarrow Y$.

\item[(iv)]
If $\pi_\ast\cO_X$ is a free $\cO_Y$-module of rank $d$, then for the contravariant
morphism $\pi^\ast$ the composition $\pi_\ast \pi^\ast$ is the multiplication with $d$ .

\item[(v)]
For a cartesian diagram of finite flat morphisms
$$
\begin{CD}
X' @>>{\pi'}> Y' \\
@VV{\rho'}V   @VV{\rho}V   \\
X  @>>{\pi}>   Y
\end{CD}
$$
one has $\pi^\ast \rho_\ast = (\rho')_\ast (\pi')^\ast: C(Y') \rightarrow C(X)$.
\end{itemize}

\noi
Then the canonical morphism $ tC(Z_.) \longrightarrow C(X) $ is an isomorphism for
any $\ell'$-hypercovering $Z. \rightarrow X$.

\vspace{0,5cm}

\noi
We note that the Gersten complexes $R_{q,\ast}(-)$ considered above are known to
satisfy the conditions (i) to (v) above. Moreover, conditions hold for for the
Bloch-Ogus complexes considered in section 4 and the complexes $\overline{C}(-,\qzl)$
considered in Theorem 4.18. Finally, as noted in \cite{GS2}, all complexes constructed
from Rost cycle modules satisfy (i) to (v).

\v
\noi\textbf{Proof of Lemma 5.21} This follows in a similar way as in the descent theorem
[Gi] Theorem 4.1. Let me very briefly recall the three steps:
All considered varieties will be in ${\mathcal V}_k$.

(I) First consider the natural morphism $Z_. = cosk_0^X(Z) \rightarrow X$
for a morphism $Z \rightarrow X$.

If $Z \rightarrow X$ has a section, then $Z_.$ is homotopy equivalent
to the constant simplicial variety $X$, and the claim follows via
the convergent spectral sequence
$$
E^1_{p,q}(Z_.) = H_p(H_q(C(Z_.))) \Rightarrow
H_{p+q}(tC(Z_.))\,, \leqno(5.22)
$$
which is the spectral sequence for the filtration with respect to the
`simplicial' degree of the bi-complex $C(Z_.)$).

If (b) holds and $\pi: Z \rightarrow X$ is finite flat, with $\pi_\ast\cO_X$
a free $\cO_Y$-module of rank $d$ prime to $\ell$, then one still gets a
homotopy equivalence $ tC(Z_.) \longrightarrow C(X) $ since $d$ is invertible in $\z_{(\ell)}$,
see (the proof of) \cite{GS2} Lemma 3.6.

(II) If we consider $Z_. = cosk_0^X(Z) \rightarrow X$ for an arbitrary $\ell'$-envelope
$Z \rightarrow X$, then by the definition of $\ell'$-envelopes and by
localization in $X$, i.e., by the exact sequence 5.16 (ii) and
the induced one for $Z.$, and by noetherian induction, we may reduce to the
case that $Z \rightarrow X$ fulfills the conditions of (I).

(III) Then, to extend this to the general case it suffices to show
that the morphism of simplicial schemes
$$
f: Z[n+1]. := cosk^X_{n+1}sk_{n+1}(Z_.) \longrightarrow cosk^X_n
sk_n(Z_.) =: Z[n]. \leqno(5.23)
$$
induces a quasi-isomorphism $ tC(Z[n+1].) \rightarrow
C(Z[n].)$ for all $n \geq 0$, because $Z[n]_j = Z_j$ for $j\leq
n$ and hence $H_j(tC(Z_.)) \mathop{\rightarrow}\limits^{\sim}
H_j(Z[n].)$ for $j<n$ by the spectral sequence (5.22). To show
that the morphism (5.23), abbreviated $f: X'_. \rightarrow X_.$,
induces a quasi-isomorphism one then follows the proof of [SGA 4,
Vbis] (3.3.3.2). In fact, as noted in [Gi], the reasoning of loc.
cit. (3.3.3.3) shows that all morphisms $F_i: X_i' = Z[n+1]_i
\rightarrow X_i = Z[n]_i$ are $\ell'$-envelopes. Note that the diagram
in loc. cit. 3.3.3.3 should read
$$
\begin{array}{ccccc}
K'_\iota & \longrightarrow & \prod\;X'_r & \mathop{\mathop{}\limits^{\longrightarrow}\limits_{\longrightarrow}}\limits^{pr_i}\limits_{X'(\iota)pr_j} & X'_i \\
\downarrow & & \downarrow & & \downarrow \\
K_\iota & \longrightarrow & \prod\; X_r & \mathop{\mathop{}\limits^{\longrightarrow}\limits_{\longrightarrow}}\limits^{pr_i}\limits_{X(\iota)pr_j} & X_i
\end{array}
$$
for the morphism $\iota: i \rightarrow j$ in \raisebox{1mm}{$\mathop{\Delta^+}\limits_{n+1[p]}$}
(notations as in \S7 below), where the product is over all objects $r$ in \raisebox{1mm}{$\mathop{\Delta^+}\limits_{n+1[p]}$}.
Moreover, the object $\cap\,K_\iota$ should rather read $\prod\,K_\iota$, which is $\times_X K_\iota$ here.

By looking at the
bi-simplical scheme $ [X'_./X_.]$ with components $[X'_./X_.]_p :=
X'_.\times_{X_.}\ldots \times_{X_.}X'_.$ ($(p+1)$ times) and its
base change with $X'_. \rightarrow X_.$, one sees that it suffices
to replace $X'_. \rightarrow X_.$ by its base change with
$[X'_./X_.]_p$ for all $p\geq 0$. This is again of the form
(5.19), and has a section $s$. So $fs = id$, and because $sf$ is
the identity on $sk_n(Z_.)$, it is homotopic to the identity ([SGA
4, Vbis], (3.0.2.4)). Therefore $f$ in this situation is a homotopy
equivalence, and hence induces a quasi-isomorphism by (5.22).
This finishes the proof of Lemma 5.20.

\v

We are now ready to give the

\vspace{0,5cm}
\noi\textbf{Proof of Theorem 5.13:} (cf. the reasoning in \cite{GS} 5.3)
In this proof, we will only
consider varieties in ${\mathcal V}_k$ and call them simply varieties.
If $(RS)_k$ holds then we set $\ell=1$ and consider hyperenvelopes and motives with coefficients in $\z$.
If $\ell$ is a prime invertible in $k$ and we are in case (b) we consider $\ell'$-hyperenvelopes and consider motives with coefficients in $\z_{(\ell)}$.

\medskip
By Lemma 5.18, any proper variety $Z$
(recall: in ${\mathcal V}_k$) has a smooth projective $\ell'$-hyperenvelope,
and every morphism $f: Y \rightarrow X$ of proper varieties has a smooth projective
$\ell'$-hyperenvelope, i.e., there is a commutative square
$$
\begin{CD}\widetilde Y. @>{\widetilde f}>>  \widetilde X. \\
@V{h_Y}VV  @VV{h_X}V \\
Y @>>{f}> X
\end{CD}
$$
in which $h_Y$ and $h_X$ are smooth projective $\ell'$-hyperenvelopes
(Let $\widetilde X.$ be a smooth projective $\ell'$-hypercovering of $X$, and
let $\widetilde Y.$ be a smooth projective $\ell'$-hypercovering of $Y\times_X \widetilde X.$).
If $Z$ is an arbitrary reduced variety over $k$, there is an open
embedding $j: Z\subset \overline{Z}$ into a proper variety. Fixing
such an embedding $j_Z$ for every reduced variety $Z$, and a
smooth projective $\ell'$-hypercovering $\widetilde f$ for every morphism $f$
of proper varieties, we can proceed as follows. For every smooth
projective simplicial variety $Z.$ define the complex $tC(Z.)$
in $\mathcal C_{\geq 0}(\mathcal A)$ as the total complex
associated to the bi-complex
$$
C(Z.): \quad\quad
\ldots \rightarrow C(Z_r)
\mathop{\longrightarrow}\limits^{\sum_{i=0}^{r}(-1)^i(d_i)_\ast}
C(Z_{r-1}) \longrightarrow \ldots \longrightarrow C(Z_0)\,,
$$
where the $d_i$ are the face morphisms of $Z.$\;, and define $H_a(Z.):= H_a(tC(Z.))$.
For every morphism $g: W. \rightarrow Z.$ of smooth projective simplicial
schemes define
$$
C(g) := Cone(tC(W.) \mathop{\longrightarrow}\limits^{g_\ast}
tC(Z.))
$$
and $H_a(g) = H_a(C(g))$.
Finally, for a variety $Z$ over $k$ define
$$
C(Z) := C(\widetilde{i_Z})\quad \mbox{and} \quad H_a(Z) :=
H_a(C(Z))\,,
$$
where $i_Z: Y \hookrightarrow \overline{Z}$ is the embedding of
the closed complement of $j_Z: Z \hookrightarrow \overline{Z}$.
Then, for any smooth projective simplicial variety $Z.$ we have a
convergent spectral sequence
$$
E^1_{p,q}(Z.) = H_q(C(Z_p)) = H^C_q(Z_p)) \Rightarrow
H_{p+q}(Z.) = H_{p+q}(tC(Z.))\,. \leqno(5.24)
$$
Here $p$ is the simplicial degree and $q$ is the complex degree, so that
$E^2_{p,q}$ is the $p$-th homology of the
complex
$$
H^C_q(Z.): \quad\quad\quad \ldots \longrightarrow H^C_q(Z_r)
\mathop{\longrightarrow}\limits^{\sum_{i=0}^{r}(-1)^i(d_i)_\ast}
H^C_q(Z_{r-1})\longrightarrow\ldots \longrightarrow H^C_q(Z_0)\,.
$$
By the assumption on the functors $H^C_q$, the latter complex is
obtained by applying the extension of $H^C_q$ to Chow motives to
the complex of motives
$$
M(Z.): \quad\quad\quad
M(Z_0) \longrightarrow \ldots \longrightarrow M(Z_{r-1})
\mathop{\longrightarrow}\limits^{\sum_{i=0}^{r}(-1)^i(d_i)^\ast}
M(Z_r) \longrightarrow \ldots \,,
$$
$M(X)$ denoting the Chow motive of a smooth projective variety
$X$. For every morphism $g: W. \rightarrow Z.$ of smooth
projective simplicial schemes we have functorially associated
morphisms $g_\ast: H_a(W.) \rightarrow H_a(Z.)$ extending to a morphism
$E(W.) \rightarrow E(Z.)$ of the spectral sequences (5.24). If the induced morphism
of motives
$$
g^\ast: M(Z_.) \longrightarrow M(W_.)
$$
is a homotopy equivalence, then $g$ induces an isomorphism on the
$E^2$-terms of the above spectral sequences, and thus an
isomorphism $g_\ast: H_a(W_.) \rightarrow H_a(Z_.)$ for all $a$.

Using this, we get the functoriality of our homology theory on
${\mathcal V}_k$ following the reasoning in \cite{GS} 2.2 and 2.3:
Let $Ar(\mathcal P_k)$ be the category of morphisms in the category
$\mathcal P_k$ of proper varieties over $k$. For every morphism
$g: f_1 \rightarrow f_2$ in $Ar(\mathcal P_k)$, i.e., any commutative
diagram
$$
\begin{CD}
Y_1 @>{f_1}>> X_1  \\
@V{g_Y}VV @VV{g_X}V \\
Y_2 @>{f_2}>> X_2
\end{CD}
$$
there exists a commutative diagram
$$
\xymatrix{\tilde f\ar[r]^b\ar[d]_a& \tilde f_2\ar[dd]^{h_2}\\
\tilde f_1\ar[d]_{h_1}\\
f_1\ar[r]^g & f_2\, .}\leqno{(5.25)}
$$
in which $h_1, h_2$ and $a$ are $\ell'$-hyperenvelopes. We claim that $a$ induces a quasiisomorphism
$C(\tilde{f})\rightarrow C(\tilde{f}_1)$.
In fact, fix $q$ and and a variety $V$, and write $D(Z.)$ for $tR_{q,\ast}(V\times Z.)$ for a simplicial scheme
$Z.$, and $D(f) = Cone(D(Z.') \rightarrow D(Z.)$ for a morphism $f:Z.' \rightarrow Z.$ of simplicial schemes.
Since $h_1$ and $h_1a$ are $\ell'$-hyperenvelopes of a morphism of varieties, they are also $\ell'$-hyperenvelopes
after passing to the product with $V$. Naming the morphisms the same, we conclude by Lemma 5.20
that in the composition
$$
D(\tilde{f}) \stackrel{a}\longrightarrow D(\tilde{f}_1) \stackrel{h_1}\longrightarrow D(f_1)\,,
$$
the maps $h_1$ and $h_1a$ are $\ell'$-quasi-isomorphisms, i.e., quasi-isomorphisms after tensoring with $\z_{(\ell)}$.
Hence $a$ is an $\ell'$-quasi-isomorphism.
Since this holds for all $V$ and all $q\geq 0$, it follows from Lemma 5.19 that the morphism of associated complexes of $\z_{(\ell)}$-motives $a: M(\tilde{f}) \rightarrow M(\tilde{f}_1)$
is a homotopy equivalence (compare \cite{GS} Corollary 1), which implies the claim.

\medskip
Thus we get morphisms
$$
g_\ast = (a_\ast)^{-1}b_\ast: H_a(C(\tilde{f_1})) \longrightarrow H_a(C(\tilde{f_2}))
$$
for all $a \in \mathbb Z$. It follows as in \cite{GS} 2.2 that these maps do not depend
on the choices and are functorial.

If now $f: Z_1 \rightarrow Z_2$ is a proper morphism of varieties, and the embeddings
$j_{Z_\nu}: Z_\nu \hookrightarrow \overline{Z}_\nu$ into proper varieties
are as above, as well as $i_\nu =i_{Z_\nu}: Y_\nu=\overline{Z}_\nu\setminus Z_\nu \hookrightarrow
\overline{Z_\nu}$ (for $\nu = 1,2$), let $Z_f \subset Z_1\times Z_2$ be the
graph of $f$ and $\overline{Z_f}$ be its closure.
Let $i_f: Y_f=\overline{Z_f}-Z_f \hookrightarrow \overline{Z_f}$ be
the closed immersion. Then the projections
$\pi_\nu: \overline{Z_f} \rightarrow \overline{Z}_\nu$ map $Y_f$ to $Y_\nu$ (for $\nu=2$
one needs the properness of $f$), and $\pi_1$ induces an isomorphism
$Z_f = \overline{Z_f} - Y_f \rightarrow \overline{Z_1} - Y_1 = Z_1$. This gives a
canonical diagram in the category $Ar(\mathcal P_k)$ of morphisms
in the category $\mathcal P_k$ of proper varieties over $k$
$$
\begin{CD}
i_f @>{\pi_2}>> i_2 \\
@V{\pi_1}VV @.\\
i_1 @. @.
\end{CD}\,.\leqno{(5.26)}
$$
Writing $D(Z.)= tR_{q,\ast}(V\times Z.)$ as above, we have a commutative
diagram with exact rows
$$
\begin{CD}
0 @>>> D(Y_f) @>>{i_f}> D(\overline{Z_f}) @>>> D(Z_f) @>>> 0 \\
@.  @VV{\pi_1}V  @VV{\pi_1}V   @V{\cong}V{\pi_1}V   @. \\
0 @>>> D(Y_1) @>>{i_1}> D(\overline{Z_1}) @>>> D(Z_1) @>>> 0\,.
\end{CD}
$$
Hence the map $\pi_1: D(i_f) \rightarrow D(i_1)$ is a quasi-isomorphism,
and for a smooth projective $\ell'$-hypercovering  $\tilde{\pi_1}: \tilde{i_f} \rightarrow \tilde{i_1}$ of $\pi_1$
the map $\tilde{\pi_1}: D(\tilde{i_f}) \rightarrow D(\tilde{i_1})$ is a quasi-isomorphism by Lemma 5.20.
Since this holds for all $V$ and all $q$ we conclude as before from Lemma 5.19 that
$\tilde{\pi_1}: C(\tilde{i_f}) \rightarrow C(\tilde{i_1})$ is an isomorphism.

\v
We obtain
$$
f_\ast := H_a(f) := H_a((\pi_2)_\ast)H_a((\pi_1)_\ast)^{-1}:
H_a(Z_1) = H_a(C(\tilde{i_1})) \rightarrow H_a(C(\tilde{i_2})) = H_a(Z_2)\,.
$$
This is functorial with the same argument as in \cite{GS} 2.3.

Moreover, we remark that for another choice of compactifications
$j_Z$ and hence maps $i_Z$, say $j'_Z$ and $i'_Z$, the morphisms
$H_a(id_Z)$ give canonical isomorphisms between the different
constructions of $H_a$. In particular, for a smooth projective variety
$X$ we get $H_a(X) = H^C_a(X)$. Also if we take another choice for the
$\ell'$-hyperenvelopes $\tilde f$, then the reasoning in \cite{GS} 2.2 shows,
that the resulting homology theory is canonically isomorphic to
the one for the first choice. In this sense, the homology theory
is canonical.

To obtain the properties 4.20 (a) (i) and (ii) for our homology
theory, i.e., contravariance for open immersions, the exact
localization sequences, and their functoriality, we proceed as in \cite{GS} 2.4: For a variety
$Z$ and a closed subvariety $Z'\subset Z$ with open complement $U
= Z\setminus Z'$ choose a compactification $Z \hookrightarrow
\overline{Z}$, and let $Y = \overline{Z}\setminus Z$ and $Y' =
\overline{Z}\setminus U$, so that $Z' = Y'\setminus Y$. Then we
choose smooth projective $\ell'$-hyperenvelopes $\widetilde{\overline{Z}}.
\rightarrow \overline{Z}$, $\widetilde Y. \rightarrow Y$ and
$\widetilde{Y'}. \rightarrow Y'$ such that one has morphisms
$$
\widetilde i: \widetilde{Y}. \mathop{\longrightarrow}\limits^{\widetilde k}
\widetilde{Y'}. \mathop{\longrightarrow}\limits^{\widetilde{i'}}
\widetilde{\overline{Z}}.
$$
lifting the closed immersions $ i: Y
\mathop{\hookrightarrow}\limits^{k} Y'
\mathop{\hookrightarrow}\limits^{i'} \overline{Z}$. We obtain an
obvious triangle of mapping cones
$$
C(\widetilde k) \longrightarrow C(\widetilde i) \longrightarrow
C(\widetilde{i'}) \longrightarrow C(\widetilde k)[-1]
$$
which
%represents the desired triangle
%$$
%C(Z') \longrightarrow C(Z) \longrightarrow C(U) \longrightarrow
%C(Z')[-1]
%$$
%in the derived category, which, in turn,
gives rise to the desired exact localization sequence
$$
\ldots \rightarrow H_a(Z') \rightarrow H_a(Z) \rightarrow H_a(U)
\rightarrow H_{a-1}(Z') \rightarrow \ldots
$$
and thereby also to the pull-back $j^\ast$ for the open immersion
$j: U \rightarrow Z$. The functorial properties are easily
checked.

\vspace{0,5cm}
\noi\textbf{Remark 5.27} In the situation of Theorem 5.13, let $C': \mathcal{SP}_k \rightarrow
\mathcal C_{\geq 0}(\mathcal A)$ be another functor with the same
properties, and let $\varphi: C \rightarrow C'$ be a
morphism of functors. Then it is clear from the construction, that
there is a canonical induced morphism $\varphi: H \rightarrow H'$
of the associated homology theories on $\mathcal V_k$.

\vspace{0,5cm}\noi
Next we want to show that we get nothing new, if we apply the construction
of Theorem 5.13 to the complexes $\overline{C}(-,\qzl)$. More generally we have the following
second descent lemma.

\v
\noi\textbf{Lemma 5.28} Let the assumption be as in Theorem 5.13,
and let $({\mathcal V}_k)_\ast$ be the
category of all varieties over $k$ with all
proper morphisms between them. Let
$$
C: ({\mathcal V}_k)_\ast \longrightarrow \mathcal C_{\geq 0}(\mathcal A)
$$
be a covariant functor satisfying the conditions in Lemma 5.21 - more precisely
the conditions (i) and (ii) in case 5.13 (a) and the conditions (i) to (v)
incase of 5.13 (b).
Let $H^{naive}$ be the obvious homology theory on $(\mathcal V_k)_\ast$
deduced from $C$, with $H^{naive}_a(Z) = H_a(C(Z))$, and the localization
sequences (see Definition 4.20 (ii)) induced by the exact sequences
in 5.21 (ii) above. On the other hand let $H$ be the homology theory on
$(\mathcal V_k)_\ast$ derived from $C_{\mid \mathcal{SP}_k}$ via Theorem 5.13. Then there is a
canonical isomorphism of homology theories
$$
\alpha:\; H \stackrel{\sim}\longrightarrow H^{naive}\;.
$$

\vspace{0,5cm}

\noi\textbf{Proof:}
Let the notations be as in the proof of Theorem 5.13, including the notion of $\ell'$-hyperenvelopes
which are hyperenvelopes for $\ell=\infty$, and the choice of an $\ell'$-hyperenvelope
$h_f: \tilde f \rightarrow f$ for each morphism $f$ in $\mathcal P_k$ and a compactification
$j_Z: Z \hookrightarrow X_Z := \overline{Z}$ with complement $i_Z: Y_Z := \overline{Z} - Z \hookrightarrow X_Z$
for each variety $Z$.
The crucial observation is that the constructions in the proof of 5.13 can also be
applied to the complexes $C(-)$. First of all we can define $C(g) := Cone(tC(Y.) \mathop{\longrightarrow}\limits^{g_\ast} tC(X.))$
for any morphism $g: Y. \rightarrow X.$ of smooth projective simplicial varieties, where $tC(Y.)$
is the total complex associated to the bi-complex $C(Y.)$.
However we can also define $C(f) := Cone(C(Y) \mathop{\longrightarrow}\limits^{f_\ast} C(X))$
for any proper morphism $f: Y \rightarrow X$ of varieties, and property (ii) above gives a canonical morphism
$$
\beta_Z: C(i_Z) = Cone(C(Y_Z) \mathop{\longrightarrow}\limits^{(i_Z)_\ast} C(X_Z)) \longrightarrow C(Z)\,,
$$
which is a quasi-isomorphism. Thus the descent lemma below gives a canonical morphism
$$
\alpha_Z: C(\tilde{i_Z}) \mathop{\longrightarrow}\limits^{(h_{i_Z})_\ast} C(i_Z) \rightarrow C(Z)
$$
which is a quasi-isomorphism and hence induces isomorphisms
$$
\alpha_Z:\;H_a(Z) = H_a(C(\tilde{i_Z})) \longrightarrow H_a(C(Z)) = H^{naive}_a(Z)\,,
$$
where the equalities holds by definition. We show the functoriality of these maps, by showing
that the corresponding constructions in the proof of Theorem 5.13 can also be carried out for $C(-)$,
compatible with the morphisms $\alpha$. For this it suffices to note that two types of morphisms
$g_1 \rightarrow g_2$ between morphisms $g_i$ of smooth projective simplicial schemes
induces quasi-isomorphisms $C(g_1) \rightarrow C(g_2)$. The first is the morphism $a$ in (5.25),
constructed for a morphism $f_1 \rightarrow f_2$ in $Ar(\mathcal P_k)$.
But $h_1$ and $ah_1$ are $\ell'$-hyperenvelopes and so induce quasi-isomorphisms
$C(\tilde{f_1}) \rightarrow C(f_1)$ and $C(\tilde{f}) \rightarrow C(f_1)$ by
the descent lemma 5.21; so
$a$ induces a quasi-isomorphism $C(\tilde{f}) \rightarrow C(\tilde{f_1})$.
The second is the morphism $\pi_1$ in (5.26), constructed for a proper morphism of varieties
$f: Z_1 \rightarrow Z_2$. Here we get a commutative diagram
$$
\xymatrix{\tilde f\ar[r]^b\ar[d]_a& \tilde i_2\ar[dd]^{h_1}\\
\tilde i_f\ar[d]_{h}\\
i_f\ar[r]^{\pi_1} & i_1\, .}
$$
where $h, h_1$ and $a$ are $\ell'$-hyperenvelopes. It induces a commutative diagram
$$
\begin{CD}
C(\tilde{f}) @>{b_\ast}>> C(\tilde{i_1}) \\
@V{(ha)_\ast}VV @VV{(h_1)_\ast}V \\
C(i_f) @>{(\pi_1)_\ast}>> C(i_1) \\
@V{(\beta_{Z_f})_\ast}VV @VV{(\beta_{Z_1})_\ast}V \\
C(Z_f) @>{(\pi_1)_\ast}>> C(Z_1)\,.
\end{CD}
$$
Here $(ha)_\ast$ and $(h_1)_\ast$ are isomorphisms by the descent lemma 5.21,
and the lower vertical maps are isomorphisms by the property (ii) for $C$.
Moreover the lower horizontal map is an isomorphism, because $\pi_1$
induces an isomorphism $Z_f \rightarrow Z_1$.
This shows that the other two horizontal maps are isomorphisms as well.
Finally the composition $(\pi_2)_\ast ((\pi_1)_\ast)^{-1}: H_a(Z_1) \rightarrow H_a(Z_2)$
coincides with the morphism $f_\ast$.
It follows from this that one gets commutative diagrams
$$
\begin{CD}
H_a(C(\tilde{i_1})) @>{b_\ast(a_\ast)^{-1}}>> H_a(C(\tilde{i_2})) \\
@V{(h_1)_\ast}VV @VV{(h_2)_\ast}V \\
H'_a(Z_1) @>{f_\ast}>> H'_a(Z_2)
\end{CD}
$$
in which the bottom map is the one following from the functoriality of $C(-)$.
It is now clear that the constructed maps $H_a(Z) \rightarrow H'_a(Z)$
are functorial. The functoriality with respect to the long exact homology sequences
follows in a similar way.

\vspace{15mm}

\centerline{\noindent\textbf{\S 6 Proof of Theorem 5.12}}

\vspace{1,0cm}

We construct the morphism of homology theories
on ${\mathcal V}_K$,
$$
\varphi:\; \overline{H}_\ast(-,\qzl) \rightarrow H^W_\ast(-,\qzl)\leqno{(6.1)}
$$
using Remark 5.27 and Proposition 5.27. We may replace $K$ by its perfect hull $K^{per}$, because both
cohomology theories only depend on the varieties $Z_{K^{per}}=Z\times_KK^{per}$, by definition.
Obviously, one has a direct sum decomposition
$\overline{C}(X_1\coprod X_2,\qzl) \cong \overline{C}(X_1,\qzl) \oplus \overline{C}(X_2,\qzl)$ for varieties $X_1$, $X_2$,
and a morphism $\overline{C}(X,\qzl) \rightarrow \overline{C}(K,\qzl) = \qzl$ induced by
the structural morphism $X \rightarrow \Spec(K)$ for any variety.
Applying this to the connected components of each smooth
projective variety, we get functorial maps
$$
\beta_X:\; \overline{C}(X,\qzl) \rightarrow \mathop{\textstyle\bigoplus}\limits_{\pi_0(X)} \qzl
$$
for all smooth projective varieties over $K$.
By 5.27 and 5.28 these induce a morphism of homology theories on ${\mathcal V}_K$
$$
\beta:\;H^{\overline C}(-) \longrightarrow H^W(-,\qzl)\;,
$$
where $H^{\overline C}$ is the homology theory associated to the functor $X \longmapsto \overline{C}(X,\qzl)$
on $\mathcal{SP}_K$ by Theorem 5.13. Now we get the morphism (5.29) as the composition
$$
\overline{H}(-,\qzl) = H(\overline{C}(-,\qzl)) \stackrel{\alpha^{-1}}\longrightarrow H^{\overline{C}}(-)
\stackrel{\beta}\longrightarrow H^W(-,\qzl)\;,
$$
where $\alpha$ is the isomorphism of homology theories existing by 5.16.

\medskip\noi
Write $\oH(-) = \oH(-,\qzl)$ and $\WH(-) = \WH(-,\qzl)$. We want to show properties 4.22 (i), (ii) and (ii)
for $\varphi: \oH \rightarrow \WH$, i.e., for any integral variety $Z$ of dimension $e$ over $K$ we want to show
\begin{itemize}
\item[(i)] $\oH_a(k(Z)) = 0$ for $a\neq e$.

\item[(ii)] $H^W_a(k(Z))=0$ for $a\neq e$.

\item[(iii)] $\varphi: \oH_e(k(Z)) \rightarrow H^W_e(k(Z))$ is an isomorphism for $\dim(Z)=e$.
\end{itemize}
Property (i) holds by Remark 4.23, and we show (ii) and (iii) by induction.
The case $d=0$ is easy. Now assume that $Z$ is integral of dimension $d>0$, and that (ii) and (iii) hold for $e<d$.

\smallskip
First assume condition 5.12 (A). Then (ii) holds by Theorem 5.10, and
it remains to show property 4.22 (iii) for the morphism (5.29).
Let $Z$ be an integral variety of dimension $d$ over $K$ and let
$V\subset Z$ be any non-empty smooth open subvariety. Let $U\subset V$
and $U\subset X$, $Y\subset X$ be as in property $(RS)_K$ (which
holds by assumption). By possibly removing a suitable further smooth hyperplane
section we may assume that $X\setminus Y_1$ is affine. For a finite field $K$ we
argue as in the proof of theorem 5.10.

If $W$ is integral of dimension $d$, then $\overline{H}_d(W,\qzl)$ is equal to
%\begin{small}
$$
\ker(H^d(K(W)\otimes_K\overline K,\qzl(d))_{G_K}\rightarrow
\mathop{\textstyle\bigoplus}\limits_{x\in
W^1}H^{d-1}(K(x)\otimes_K\overline K,\qzl(d-1))_{G_K}),
$$
%\end{small}
and if $W$ is irreducible and smooth of dimension $d$ then the
Bloch-Ogus niveau spectral sequence (see the end of section 4) gives a canonical
edge morphism
$$
\gamma_W: H^d(\overline W,\qzl(d))_{G_K} \longrightarrow \oH_d(W,\qzl)\,,
$$
which is given by the restriction map
$H^d(\overline W,\qzl(d))_{G_K} \rightarrow H^d(K(W)\otimes_K\overline K,\qzl(d))_{G_K}$.

\vspace{0,5cm}

\noindent\textbf{Lemma 6.2} (a) There is a commutative diagram
$$
\begin{matrix}
0 & \rightarrow & H^d(\overline{U},\qzl(d))_{G_K} &
\mathop{\rightarrow}\limits^{e} &
H^0(\overline{Y^{[d]}},\qzl(0))_{G_K} &
\mathop{\rightarrow}\limits^{d_2} &
H^2(\overline{Y^{[d-1]}},\qzl(1))_{G_K} \\
&& \downarrow\rlap{$\gamma_U$} && \downarrow\rlap{$\gamma_{Y^{[d]}}$} && \downarrow\rlap{$\gamma'$} \\
0 & \rightarrow & \overline{H}_d(U,\qzl) & \mathop{\rightarrow}\limits^{e} &
\overline{H}_0(Y^{[d]},\qzl) & \mathop{\rightarrow}\limits^{d_2} &
\overline{H}_0(Y^{[d-1]},\qzl) \\
&& \downarrow\rlap{$\varphi_U$} && \downarrow\rlap{$\varphi_{Y^{[d]}}$} && \downarrow\rlap{$\varphi_{Y^{[d-1]}}$} \\
0 & \rightarrow & H^W_d(U,\qzl) &
\mathop{\rightarrow}\limits^{e} & H^W_0(Y^{[d]},\qzl) &
\mathop{\rightarrow}\limits^{d_2} &  H^W_0(Y^{[d-1]},\qzl)\,.
\end{matrix}
$$
Here the maps $e$ and $d_2$ in the first row are those occurring in
Theorem 3.1. The maps in the second and third row are the
homological analogues: Writing $H_a$ for either $\overline{H}_a$ or $H^W_a$, $e$ is the composition of the morphisms
$$
H_d(U,\qzl) \mathop{\rightarrow}\limits^{\delta}
H_{d-1}(Y_{i_d}\smallsetminus(\mathop{\cup}\limits_{i\neq i_d}
Y_i),\qzl) \mathop{\rightarrow}\limits^{\delta} \ldots
$$
$$
\ldots \mathop{\rightarrow}\limits^{\delta}
H_1(Y_{i_2,\ldots,i_d}\smallsetminus(\mathop{\cup}\limits_{i\neq
i_2,\ldots,i_d} Y_i),\qzl) \mathop{\rightarrow}\limits^{\delta}
H_0(Y_{i_1,\ldots,i_d},\qzl)
$$
where each $\delta$ is the connecting morphism for the obvious
localization sequence for the homology theory, and $d_2=
\sum_{\mu=1}^{d}(-1)^\mu\delta_\mu$, where $\delta_\mu$ is induced
by the push-forward morphisms for the inclusions
$Y_{i_1,\ldots,i_d} \hookrightarrow
Y_{i_1,\ldots,\hat{i_\mu},\ldots,i_d}$. Finally,
$\gamma'=\gamma_{Y^{[d-1]}}\circ tr'$ where
$$
tr':
H^2(\overline{Y^{[d-1]}},\qzl(1))_{G_K}
\mathop{\rightarrow}\limits^{\sim}
H^0(\overline{Y^{[d-1]}},\qzl(0))_{G_K}
$$
is the isomorphism induced by the trace maps $H^2(C,\qzl(1))\stackrel{\sim}\rightarrow \qzl$ for
the smooth projective curves $C$ which are the irreducible components of $Y^{[d-1]}$.

(b) The composition
$$
H^d(\overline{U},\qzl(d))_{G_K}
\mathop{\longrightarrow}\limits^{\gamma_U} \oH_d(U,\qzl)
\mathop{\longrightarrow}\limits^{\varphi_U} H^W_d(U,\qzl)
$$
figuring in the left column of the above diagram is an isomorphism.
\vspace{0,5cm}

\noindent\textbf{Proof} (a): For any smooth irreducible variety $W$ of
dimension $d$ over $K$, and any smooth irreducible divisor $i:
W'\hookrightarrow W$, we have a commutative diagram
$$
\begin{matrix}
H^d(\overline{W\smallsetminus W'},\qzl(d)) &
\mathop{\rightarrow}\limits^{\delta} &
H^{d-1}(\overline{W'},\qzl(d-1)) \\
\downarrow && \downarrow \\
H^d(K(W)\otimes_K\overline{K},\qzl(d)) &
\mathop{\rightarrow}\limits^{\delta} &
H^{d-1}(K(W')\otimes_K\overline{K},\qzl(d-1))
\end{matrix}\,,
$$
where the $\delta$ in the top line is the connecting morphism for
the Gysin sequence for $W'\subseteq W \supseteq W\smallsetminus
W'$, and the $\delta$ in the bottom line is the residue map for
the point in $W^1$ corresponding to $W'$. The latter induces the
connecting morphism
$$
H_d(W\smallsetminus W',\qzl) \mathop{\rightarrow}\limits^{\delta}
H_{d-1}(W',\qzl)
$$
for the localization sequence  for $(W',W,W\smallsetminus W')$.
This shows the commutativity of the top left square in 6.2, by
definition of the maps $e$.

On the other hand, for an irreducible smooth projective curve $C$ over
$\overline{K}$, and a closed point $P: \Spec(\overline{K})
\rightarrow C$, the composition
$$
H^0(\overline{K},\qzl(0)) \mathop{\longrightarrow}\limits^{P_\ast}
H^2(C,\qzl(1)) \mathop{\longrightarrow}\limits^{tr}
H^0(\overline{K},\qzl(0))
$$
is the identity. This shows the commutativity of the top right
square in 6.2. The two bottom squares commute because $\varphi$
is a morphism of homology theories.

\medskip
(b): The compositions of the vertical maps in the middle column and
the right column of the diagram in 6.2 are isomorphisms,
and the top row is exact by Theorem 3.1 (and our assumption on $U$).
But the bottom line is exact as well: This
follows in a similar (but simpler) way as in the proof of Theorem
3.1, by noting that $H^W_a(T,\qzl) = 0$ for $a\neq 0$ if $T$
is smooth and projective of dimension $>0$, by definition. It can also be deduced
from the fact that $H^W_a(U,\qzl)$ is computed as the $a$-th
homology of the complex
$$
\mathop{\textstyle\bigoplus}\limits_{\pi_0(Y^{[d]})}\,\qzl
\rightarrow
\mathop{\textstyle\bigoplus}\limits_{\pi_0(Y^{[d-1]})}\,\qzl
\rightarrow \ldots \quad \ldots \rightarrow
\mathop{\textstyle\bigoplus}\limits_{\pi_0(Y^{[0]})}\,\qzl
$$
as noted before. This shows (b).

\vspace{0,5cm}
We can now prove property 4.23 (iii) for $\varphi$:
 Because the
subvarieties $U$ as constructed above form a cofinal family in the
set of open subvarieties of $Z$, by passing to the limit we get an
isomorphism
$$
H^d(K(Z)\otimes_K\overline{K},\qzl(d))_{G_K}
\mathop{\longrightarrow}\limits^{\gamma} H_d(K(Z),\qzl)
\mathop{\longrightarrow}\limits^{\varphi} H^W_d(K(Z),\qzl)
$$
in which the first map $\gamma$ is an isomorphism by definition.
Therefore $\varphi$ is an isomorphism as wanted.

\medskip
We note that the above arguments can be applied for an irreducible smooth projective curve $X$,
where, for an open subvariety $U$, one has an exact sequence of weight homology
$$
0 = H^W(X) \rightarrow H^W_1(U) \rightarrow \mathop{\textstyle{\bigoplus}}\limits_{x\in X\setminus U} H^W(x) \rightarrow H^W_0(X) \rightarrow 0\,,
$$
and all other weight homology groups vanish. This shows (ii) for $K(X)$, and (iii) follows with Lemma 6.2.

\bigskip
Now assume condition 5.12 (B), i.e., that $\ell$ is invertible in $K$.
Again we assume that $Z$ is integral of dimension $d>0$, and that (ii) and (iii) hold for $e<d$.
By the preceding remark we may assume $d>1$. First we show property (ii) for $K(Z)$.

\medskip
Quite generally the induction assumption on (ii) and the spectral sequence
$$
E^1_{p,q}(Y)=\textstyle\bigoplus\limits_{x\in Y_p}
H_{p+q}(k(x))\Rightarrow H_{p+q}(Y)\leqno{(6.3)}
$$
(see (4.25)) show that

\medskip\noi
(6.4) \quad $H^W_a(Y)=0$ for $a>e$ if $Y$ is a variety of dimension $e<d$.

\medskip\noi
Moreover the induction assumption on (i), (ii) and (iii) and the niveau spectral sequences
for $H^W$ and $\oH$ show that

\medskip\noi
(6.5) \quad $\oH_a(Y) \mathop{\longrightarrow}\limits^{\sim} H^W_a(Y)$ is an isomorphism for all $a$ if dim($Y)<d$.

\medskip\noi
(See the proof of Lemma 4.22.)
Now we want to show (ii) and (iii) for our integral variety $Z$ of dimension $d>1$.

\medskip
First assume that the function field $k(Z)$ has a smooth projective model X, i.e.,
$k(Z) = k(X)$. Let $U \subset X$ be a non-empty open, and let $Y=X-U$. Then the long exact homology sequence
$$
\ldots \rightarrow H^W_a(Y) \rightarrow H^W_a(X) \rightarrow H^W_a(U) \rightarrow H^W_{a-1}(Y) \rightarrow \ldots \leqno{(6.6)}
$$
shows that $H^W_a(U)=0$ for $a>d$, because $H^W_a(X)=0$ for $a>0$ by the definition of
weight homology for smooth projective varieties. Since this holds for all open dense $U\subset X$, we get

\medskip\noi
(6.7) \quad $H^W_a(k(X))=0$ for $a>d$.

\medskip\noi
To show (ii) for $k(X)$ we have to show $H^W_a(k(X))= 0$ for $a<d$. Now we claim that for the homology
theories $H^W$ and $\oH$ the following holds: For any commutative diagram
$$
\begin{CD}
Z @>{i}>> X @<{j}<< U = X\setminus Z \\
@V{k}VV   @V{id}VV   @A{h}AA    \\
Z'  @>{i'}>>   X   @<{j'}<<   U' = X\setminus Z'
\end{CD}\leqno{(6.8)}
$$
with closed immersions $i, k$, and $i'$ and (hence) open immersions $j, h$, and $j'$ we have
a morphism of long exact homology sequences
$$
\begin{CD}
\ldots @>>> H_a(Z)  @>{i_\ast}>>  H_a(X)  @>{j^\ast}>>     H_a(U)  @>{\delta}>>   H_{a-1}(Z) @>>> \ldots \\
@.       @V{k_\ast}VV            @V{id}VV           @V{h^\ast}VV                 @V{k_\ast}VV        @.  \\
\ldots @>>> H_a(Z') @>{i'_\ast}>> H_a(X)  @>{(j')^\ast}>>  H_a(U') @>{\delta}>>   H_{a-1}(Z') @>>> \ldots \,.
\end{CD}\leqno{(6.9)}
$$
This is one of the general axioms of a homology theory, and it holds for the homology theories
constructed by Theorem 5.13 (hence for $H^W$) and for the homology theories considered in Lemma 5.28
(hence for $\oH$). In fact, in these cases the diagram (6.9) arises from a morphism of exact triangles
of complexes
$$
\begin{CD}
C(Z)  @>{i_\ast}>>   C(X)    @>>>  cone(C(Z) \rightarrow C(X)) @>>>  \\
@V{k_\ast}VV       @V{id}VV          @VV{(k,id)}V          \\
C(Z') @>{i'_\ast}>>  C(X)    @>>>  cone(C(Z') \rightarrow C(X))  @>>> .
\end{CD}
$$
Therefore we can pass to the inductive limit over all closed subvarieties $Y \subset X$ in (6.6),
and this gives a long exact sequence
$$
\ldots \rightarrow H^W_{a+1}(X) \rightarrow H^W_{a+1}(K(X)) \rightarrow \mathop{lim}_{\rightarrow,Y} H^W_a(Y) \rightarrow H_a^W(X) \ldots \,. \leqno{(6.10)}
$$
From this we deduce: If we can show

\medskip\noi
(6.11) \quad  $\mathop{lim}\limits_{\rightarrow,Y} H^W_a(Y) = 0$ for $1 \leq a \leq d-2$\,,

\medskip\noi
then we deduce $H^W_a(K(X))=0$ for $2 \leq a \leq d-1$, because $H^W_a(X)=0$ for $a>0$.
If in addition we can show

\medskip\noi
(6.12) \quad  $\mathop{lim}\limits_{\rightarrow,Y} H^W_0(Y) \rightarrow H^W_0(X)$ is an isomorphism,

\medskip\noi
then we also get $H^W_1(K(X)) = 0 = H^W_0(K(X))$ and hence (ii) for $K(X)$, viz., $H^W_a(K(X))= 0$
for $a\neq d$.

By (6.5) it suffices to show (6.11) for $\oH$ in place of $H^W$.
For this we will use the following result, generalizing the facts used in the case of (3.19).

\v
\noi\textbf{Lemma 6.13} Let $K$ be a perfect field, let $X$ and $X'$ be irreducible smooth varieties of dimension $d$ over $K$,
let $Y\subset X$ be a closed subscheme, and let $n,r, s \in \z$, with $n$ invertible in $K$. Let $\overline{K}$ be the algebraic closure
of $K$, let $G_K = Gal(\overline{K}/K)$, and let $\overline{V}=V\times_K\overline{K}$ for any variety $V$ over $K$.

(a) There are canonical isomorphisms for all $a\in\z$
$$
\begin{CD}
H^{a}_Y(X,\ocH_X^{r+d}(\zn(s+d)))  @>{\sim}>>  H_{d-a}(C^{r,s}(\overline{Y},\zn)_{G_K})\,,
\end{CD}
$$
where $\ocH_X^i(\zn(j))$ is the Zariski sheaf associated to the presheaf $U \mapsto H^i(\overline{U},\zn(j))_{G_K}$
($U$ open in $X$) and $C^{r,s}(\overline{Y},\zn)$ is the Bloch-Ogus-Kato complex as defined in section 4.

(b) For a closed subset $Z \subset Y$ the following diagram of canonical maps commutes:
$$
\begin{CD}
H^{a}_Y(X,\ocH_X^{i+d}(\zn(j+d)))  @>{\sim}>>  H_{d-a}(C^{i,j}(\overline{Y},\zn)_{G_K}) \\
@AAA      @AAA   \\
H^{a}_Z(X,\ocH_X^{i+d}(\zn(j+d)))  @>{\sim}>>  H_{d-a}(C^{i,j}(\overline{Z},\zn)_{G_K})
\end{CD}
$$

(c) Let $\pi: X' \rightarrow X$ be a proper dominant morphism such that $[K(X'):K(X)]$ is prime to $n$,
and let $Y' = \pi^{-1}(Y) \subset X'$. Then the pull-back maps
$$
\pi^\ast:\; H^a_Y(X,\ocH_X^r(\zn(s))) \rightarrow  H^a_{Y'}(X',\ocH_{X'}^r(\zn(s)))
$$
are injective.

\v\noi\textbf{Proof:}
Let $\ocC_X^{r,s}(\zn)$ be the complex of Zariski sheaves on $X$ associated to the presheaf $U \mapsto C^{r,s}(\overline{U},\zn)_{G_K}$.
It can be written as a complex of flasque sheaves
$$
(i_\eta)_\ast H^{r+d}(K(X)\otimes_K\overline{K},\zn(s+d))_{G_K} \stackrel{\partial}{\rightarrow}
\bigoplus\limits_{x\in X_{d-1}} (i_x)_\ast H^{r+d-1}(k(x)\otimes \overline{K},\zn(s+d-1))_{G_K} \stackrel{\partial}{\rightarrow} \ldots
$$
$$
\ldots \rightarrow \bigoplus\limits_{x\in X_0} (i_x)_\ast H^r(k(x)\otimes\overline{K},\zn(s))_{G_K}\,,
$$
where $\eta = \Spec(K(X))$ is the generic point of $X$, $i_x: x \hookrightarrow X$ is the inclusion of a point $x\in X$,
and the group $H^i(k(x)\otimes \overline{K},\zn(s))$ is regarded as a Zariski sheaf on $x$.
By the Bloch-Ogus theory and the universal exactness discussed after (3.19), there is a natural inclusion
$\ocH_X^{r+d}(\zn(s+d)) \hookrightarrow (i_\eta)_\ast H^{r+d}(K(X)\otimes_K\overline{K},\zn(s+d))_{G_K}$
which makes
$$
\ocH_X^{r+d}(\zn(s+d)) \rightarrow \ocC_X^{r,s}(\zn)
$$
a flasque resolution of $\ocH_X^{r+d}(\zn(s+d))$. This implies (a) and (b),
because $H^0_Y(\ocC_X^{r,s}(\zn)) = C^{r,s}(Y,\zn)$.

As for (c), the pull-back map can be described as the natural map
$$
\pi^\ast:\; H^a_Y(X,\ocH_X^i(\zn(j))) \rightarrow  H^a_Y(X,R\pi_\ast\ocH_{X'}^i(\zn(j)))\,,\leqno{(6.14)}
$$
induced by the pull-back morphism $\ocH_X^i(\zn(j)) \rightarrow \pi_\ast\ocH_{X'}^i(\zn(j))$.
On the other hand, $R\pi_\ast\ocH_{X'}^i(\zn(j))$ is represented by $\pi_\ast \ocC_{X'}^{r,s}(\zn)$,
by the flasqueness of $\ocC_{X'}^{r,s}(\zn)$. Now there is a natural push-forward morphism
$$
\pi_\ast: \pi_\ast \ocC_{X'}^{r,s}(\zn) \longrightarrow  \ocC_{X}^{r,s}(\zn)
$$
induced by the maps $\pi_\ast: C^{r,s}(\pi^{-1}(U),\zn) \rightarrow C^{r,s}(U,\zn)$ for all open $U \subset X$,
and we obtain a pushforward map
$$
\pi_\ast: H^a_Y(X,H^a_Y(X,R\pi_\ast\ocH_{X'}^i(\zn(j)))) \rightarrow H^a_Y(X,\ocH_X^i(\zn(j)))\,.\leqno{(6.15)}
$$
Now we claim that the composition $\pi_\ast\pi^\ast$ of (6.14) and (6.15) is the multiplication by $[K(X'):K(X)]$,
which implies (c) since $n$ is coprime to this degree. But for this it suffices to show that the resulting
morphism on the zero cohomology sheaf $\ocH_X^i(\zn(j))$ is the multiplication by $[K(X'):K(X)]$.
Since this sheaf is a subsheaf of the constant sheaf $H^i(K(X)\otimes_K\overline{K},\zn(j))$ by purity,
the claim follows from the corresponding considerations after (3.19).

\medskip
We will now show (6.11) for $\oH$ in place of $H^W$, i.e.,

\medskip\noi
(6.16) \quad  $\mathop{lim}\limits_{\rightarrow,Y} \oH_a(Y) = 0$ for $1 \leq a \leq d-2$\,,

\noi
where the inductive limit is over all proper closed subvarieties of the smooth projective variety $X$
of dimension $d$. Let $Y \subset X$ be a proper closed subscheme. Then, as in (3.19), we have a commutative diagram
$$
\begin{array}{ccccc}
U' & \subset & X' & \supset & Y' \\
\downarrow  &&  \downarrow\rlap{$\pi$}  &&  \downarrow  \\
U & \subset & X & \supset & Y
\end{array}\leqno{(6.17)}
$$
where $X'$ is an geometrically irreducible smooth and projective variety over a finite extension $K'$ of $K$
with $\ell$ not dividing $[K':K]$, $\pi$ is a proper surjective morphism which is generically finite of degree prime to $\ell$,
$U' = \pi^{-1}(U)$, and $Y'=\pi^{-1}(Y)$ is a divisor with strict normal crossings on $X'$. Moreover, as
as in the situation of (3.19), we can assume that $K$ is a perfect field, so that $X$ and $X'$ are smooth
over $K$.

\medskip
Let $T \subset X'$ be a smooth hyperplane section intersecting $Y'$ transversally, so that $X\setminus T$ is affine.
If $K$ is a finite field, such a section exists after passage to a finite extension $L/K$ with $\ell$ not dividing
$[L:K]$, which does not matter for our purposes. Let $Y_0'=Y'\bigcup T$, which is again a divisor with normal crossings.
Moreover let $Y_0=\pi(Y_0')$, and let $Y_0''=\pi^{-1}(Y_0)$, so that $Y \subset Y_0$ and $Y' \subset Y_0' \subset Y_0''$.

\medskip
By Lemma 6.13 we have injective pull-back maps $\pi^\ast: \oH_a(X) \rightarrow \oH_a(X')$, $\pi^\ast: \oH_a(Y) \rightarrow \oH_a(Y')$,
and $\pi^\ast: \oH_a(Y_0) \rightarrow \oH_a(Y_0'')$ for all $a$, and a commutative diagram
$$
\begin{CD}
\oH_a(Y') @>>> \oH_a(Y_0'') @>>> \oH_a(X') \\
@AA{\pi^\ast}A     @AA{\pi^\ast}A     @AA{\pi^\ast}A  \\
\oH_a(Y)  @>>> \oH_a(Y_0)   @>>> \oH_a(X)\,,
\end{CD}  \leqno{(6.18)}
$$
in which all vertical maps are injective. Moreover we have a factorization
$\oH_a(Y') \rightarrow \oH_a(Y_0') \rightarrow \oH_a(Y_0'')$, with

\medskip\noi
(6.19) \quad $\oH_a(Y_0')=0$ for $1\leq a \leq d-2$.

\medskip\noi
In fact, we have an isomorphism $\oH_a(Y_0') \cong H^W_a(Y_0')$ by (6.5), and a long exact sequence
$$
\ldots \rightarrow H^W_{a+1}(X') \rightarrow H^W_{a+1}(U_0') \rightarrow H^W_a(Y_0') \rightarrow H^W_a(X') \rightarrow \ldots\,,\leqno{(6.20)}
$$
where $U_0' = X'\setminus Y_0'$. The claim follows, since $H^W_a(X') = 0$ for $a>0$ by definition,
and $H^W_a(U_0') = 0$ for $a<d$ by Proposition 5.4. From diagram (6.18) we now deduce that
$\oH_a(Y) \rightarrow \oH_a(Y_0)$ is the zero map for $1\leq a \leq d-2$. Since $Y$ was arbitrary,
we obtain (6.16).

\medskip\noi
The sequence (6.20) also implies an isomorphism $H^W_0(Y_0') \mathop{\rightarrow}\limits^{\sim} H^W_0(X')$,
and we have a commutative diagram
$$
\begin{CD}
H^W_0(Y_0') @>{\sim}>> H^W_0(X')\\
@VV{tr}V               @V{\cong}V{tr}V   \\
H^W_0(\Spec(K)) @>{=}>> H^W_0(\Spec(K))\,,
\end{CD} \leqno{(6.21)}
$$
in which the vertical maps are the trace morphisms, given by the
push-forward $tr=f_\ast: H^W_0(V) \longrightarrow H^W_0(\Spec(K))$
for a proper variety $f: V \rightarrow \Spec(K)$. Here the right-hand trace
map is an isomorphism by definition, so the left-hand trace is an isomorphism, too.
By (6.5) we see that $tr: \oH_0(Y_0') \rightarrow \oH_0(\Spec(K))$ is an isomorphism.
On the other hand, we have a commutative diagram
$$
\begin{CD}
\oH_0(Y_0') @>{tr}>>  \oH_0(\Spec(K)) \\
@V{\pi_\ast}VV             @V{=}VV     \\
\oH_0(Y_0)  @>{tr}>>  \oH_0(\Spec(K))\,,
\end{CD} \leqno{(6.22)}
$$
in which $\pi_\ast$ is a surjection by (the proof of) Lemma 6.3.
We conclude that all maps in (6.22) must be isomorphisms. By (6.5) we
conclude that $H^W_0(Y_0) \mathop{\rightarrow}\limits^{tr} H^W_0(\Spec(K))$
is an isomorphism, and by the diagram
$$
\begin{CD}
H^W_0(Y_0)  @>>>  H^W_0(X) \\
@VV{tr}V           @V{\cong}V{tr}V \\
H^W_0(\Spec(K)) @>{=}>> H^W_0(\Spec(K))
\end{CD}
$$
we conclude that $H^W_0(Y_0) \rightarrow H^W_0(X)$ is an isomorphism.
Since $Y$ was arbitrary, we obtain (6.12), and this finishes the proof of (ii)
for $K(X)$.

\medskip
We will now show (iii) for $K(X)$. We have a commutative diagram with exact rows
\begin{footnotesize}
$$
\begin{CD}
0 = H^W_d(X) @>>> H^W_d(K(X)) @>>> \mathop{lim}\limits_{\rightarrow,Y} H^W_{d-1}(Y) @>>> H^W_{d-1}(X) @>>> H^W_{d-1}(K(X))=0 \\
@AAA   @AA{\varphi_{K(X)}}A   @AA{\cong}A   @AAA   @AAA  \\
0 = \oH_d(X) @>>> \oH_d(K(X)) @>{\delta}>> \mathop{lim}\limits_{\rightarrow,Y} \oH_{d-1}(Y) @>{I}>> \oH_{d-1}(X) @>>> \oH_{d-1}(K(X))=0
\end{CD}.\leqno{(6.23)}
$$
\end{footnotesize}\noi
Here the groups $H^W_d(X)$ and $H^W_{d-1}(X)$ vanish by definition of the weight homology
(note $d>1$), we have $\oH_d(X)=0$ by Theorem 3.17,
and the vanishing of the groups on the right is obvious for $\oH$ and follows from (ii) for $H^W$.
If we show the vanishing of $\oH_{d-1}(X)$ or, equivalently, the vanishing of $I$,
we have shown that $\varphi_{K(X)}$ is an isomorphism, i.e., (iii) for $K(X)$.

\medskip
We consider the situation of diagram (6.17) again. With the notations above, we have a commutative
diagram with exact rows
$$
\begin{CD}
H^W_d(U_0') @>{\delta^W}>> H^W_{d-1}(Y_0') @>>> H^W_{d-1}(X') = 0 \\
@AA{\varphi_{U_0'}}A   @AA{\cong}A   @AAA  \\
\oH_d(U_0') @>{\overline{\delta}}>>  \oH_{d-1}(Y_0') @>{\alpha}>> \oH_{d-1}(X')
\end{CD},
$$
in which $\varphi_{U_0'}$ is a surjective by Lemma 6.2 (b), because $U_0'$
satisfies the assumptions of that lemma. Since $\delta^W$ is surjective as well,
we get the surjectivity of $\overline{\delta}$, i.e., the vanishing of $\alpha$.

\medskip
Now we have a commutative diagram
$$
\begin{CD}
\oH_{d-1}(Y') @>{i_\ast'}>>  \oH_{d-1}(X') \\
@AAA   @AAA   \\
\oH_{d-1}(Y)  @>{i_\ast}>>   \oH_{d-1}(X)
\end{CD},
$$
in which the vertical maps are injective by Lemma 6.13 (c), and $i_\ast'$ factors
through $\alpha$ and hence is zero. This shows that $i_\ast$ is zero. Since $Y$ was
arbitrary, we conclude $I=0$ in (6.23) as wanted. This proves (i), (ii) and (iii)
if $K(Z)$ has a smooth proper model.

\bigskip
In general, we shall use the following lemma.

\v
\noi\textbf{Lemma 6.24}
Let $Z. \rightarrow Z$ be a proper $\ell'$-hypercovering, where $\ell$ is a prime or $\ell=\infty$
(see the convention above 5.16). Let $C$ be a functor with values in a $\z_{(\ell)}$-linear
abelian category satisfying the assumptions of Theorem 5.13 (where we let $\z_{(\infty)}=\z$, see 5.20),
and let $H$ be the homology theory associated to it via Theorem 5.13.

(a) For any open $U\subseteq Z$ let $U. \rightarrow U$ be the base change of $Z. \rightarrow Z$
with $U\hookrightarrow Z$. There is a canonical spectral sequence
$$
E^1_{p,q}(H,U) =  H_q(X^U_p)  \Rightarrow  H_{p+q}(U)
$$
generalizing the spectral sequence (5.24) which covers the case where $U.$ is proper and smooth.

(b) For an open immersion $\alpha: V\subseteq U$ and the base change $V. \rightarrow V$ of $Z. \rightarrow Z$
with $V \rightarrow U$ there is a morphism of spectral sequences
$$
\alpha^\ast: E^\ast_{p,q}(H,U)  \longrightarrow  E^\ast_{p,q}(H,V)\,,
$$
functorially in $V$.

(c) If $\varphi: C \rightarrow C'$ is a morphism of functors as above, and $\varphi$ also denotes
the morphism $H \rightarrow H'$ of associated homology theories, one has an associated morphism of spectral sequences
$$
E^\ast_{p,q}(H,U) \longrightarrow E^\ast_{p,q}(H',U)\,,
$$
functorial in $U$.

\v\noi\textbf{Proof:} (a) Let $Y. \rightarrow Y$ be the (componentwise) closed complement of $U. \rightarrow U$
in $Z. \rightarrow Z$. Then $Y. \rightarrow Z.$ is an $\ell'$-hypercovering of the morphism $Y \rightarrow X$.
By  7.8 there exists a smooth projective $\ell'$-hypercovering $\tilde{Y}. \rightarrow \tilde{X}.$ of $Y \rightarrow Z$,
and a simplicial smooth projective $\ell'$-hypercovering $Y.. \rightarrow X..$ of $Y. \rightarrow Z.$\;.
Then $Y_p. \rightarrow X_p.$ is a smooth projective $\ell'$-hypercovering of $Y_p \rightarrow Z_p$,
so that $cone(C(Y_p.) \rightarrow C(X_p.))$ computes the homology of $U_p$,
and the techniques used in the proof of Theorem 5.13 show that the total complex associated to
$cone(C(Y..) \rightarrow C(X..))$ is quasi-isomorphic to the total complex associated to
$cone(C(\tilde{Y}.) \rightarrow C(\tilde{X}.))$, which computes the homology of $U$. The spectral sequence
now arises by the filtration of $cone(C(Y..) \rightarrow C(X..))$ with respect to the first degree.

The functorial properties (b) and (c) can easily be seen from this.

\bigskip
We now prove (ii) and (iii) for $K(Z)$ where $Z$ is any irreducible proper variety of dimension $d$, assuming
(ii) and (iii) hold for function fields of smaller degree of transcendence.
Let $X. \rightarrow Z$ be some strict proper smooth $\ell'$-hypercovering, and, for any open $U\subset Z$
let $X^U. \rightarrow U$ be the base change with $U \hookrightarrow Z$. Then we get a morphism of spectral sequences
$$
E^\ast_{p,q}(\oH,U) \rightarrow  E^\ast_{p,q}(H^W,U)
$$
By passing to the inductive limit over all $U$ we get a morphism of spectral sequences
for the base change $X^\eta.$ of $X.$ with the generic point $\eta$ of $X$, between
$$
E^1_{p,q}(\oH,\eta) = \oH_q(X^\eta_p) \Rightarrow \oH_{p+q}(\eta)
$$
and
$$
E^1_{p,q}(H^W,\eta) = H^W_q(X^\eta_p) \Rightarrow H^W_{p+q}(\eta)\,.
$$
Now, by the definiton of small coverings and the explicit description (7.2) of coskeleta,
one shows that, for each $n$, $(X^\eta_n)_{red}$ is a product of spectra of finite field extensions of $K(\eta)$,
which are function fields of smooth projective varieties of dimension $d$.
Thus, by what was proved above, for function fields of smooth projective varieties,
the morphisms of initial terms are isomorphisms, so that we get isomorphisms
$\oH_a(K(Z)) \mathop{\rightarrow}\limits^{\sim} H^W_a(K(Z)$ for all $a$,
and in particular (iii) for $K(Z)$.
Now (ii) for $K(Z)$ follows from (i) for $K(Z)$.

\vspace{1,5cm}

\centerline{\noindent\textbf{\S 7 Appendix: Hyperenvelopes and $\ell'$-hypercoverings}}

\vspace{1,0cm}

We collect some results on hypercoverings and hyperenvelopes which are
scattered or not explicit enough in the literature,
and we introduce so-called simplicial hypercoverings.
\\

Let $\cC$ be a category in which finite inverse (not necessarily filtered) limits exist, i.e., in which
finite products and difference kernels exist. Let $\Delta^o\cC=Hom(\Delta^o,\cC)$ be the category
of simplicial objects in $\cC$, i.e., contravariant functors from the category
$\Delta$ of the sets $[m]=\{0,\ldots,m\}$ with monotone (non-decreasing) maps.
As usual we write the objects $X$ in $\Delta^o\cC$ as $X.$, where $X_n$ stands for $X([n])$;
similarly for the truncated case.
Let $\Delta_n \subset \Delta$ be the full subcategory formed by the objects $[m]$
with $m\leq n$, and let $\Delta^o_n\cC=Hom(\Delta^o_n,\cC)$ be the category of
$n$-truncated simplicial objects in $\cC$, see \cite{SGA 4.2}, expos\'e V bis.
In these categories again finite inverse limits exist, and are formed `componentwise',
i.e., the functors $X. \mapsto X_m$ commute with the limits.
The functor of restriction
$$
sk_n : \Delta^o\cC \longrightarrow \Delta^o_n\cC\,,
$$
called $n$-skeleton, commutes with finite inverse limits and has a right adjoint
$$
cosk_n: \Delta^o_n\cC \longrightarrow \Delta^o\cC\,,
$$
called $n$-coskeleton, which commutes with finite inverse limits again. One has
$$
cosk_n(X)_m = \mathop{\mathrm{lim}}\limits_{\mathop{\leftarrow}\limits_{\mathop{\Delta}\limits_{n[m]}}} X_q\,,
\leqno{(7.1)}
$$
where
%$\Delta(\leq n)\subset \Delta$
$\mathop{\Delta}\limits_{n}$ is the full subcategory of $\Delta$ formed
by the objects $[q]$ with $q\leq n$ and $\mathop{\Delta}\limits_{n[m]}$ is the category of
arrows $[q]\rightarrow [m]$ in $\mathop{\Delta}\limits_{n}$.  From this it folllows that the
adjunction morphisms $Z. \rightarrow sk_n(cosk_n(Z.))$ are isomorphisms for
any $n$-truncated simplicial object $Z.$,
i.e., one can identify $cosk_n(Z)_m = Z_m$ for $m \leq n$ (because then $\mathop{\Delta}\limits_{n[m]}$
has the final object $id_{[m]}$).\\

If we apply this to the subcategory $\cC/S$ of objects over a fixed base object $S$, we
obtain functors which we call $sk^S_n$ and $cosk^S_n$, respectively.\\
%(If we do not mention $S$ we can always take $S = \Spec(\z)$.)

An explicit formula is $cosk^S_n(X) = \prod_S   X_\iota$ where
$$
X_\iota  \longrightarrow  \prod_S\; X_r  \mathop{\mathop{}\limits^{\longrightarrow}\limits_{\longrightarrow}}\limits^{pr_i}\limits_{X(\iota)pr_j}  X_i\,,\leqno{(7.2)}
$$
is the difference kernel. Here $r$ runs through the objects and $\iota:\,i\rightarrow j$ through the maps in $\mathop{\Delta}\limits_{n[m]}$.\\

We may write $sk_n$ because this construction does not depend on $S$. More generally,
fix a simplicial $S$-object $T.$ and consider the categories $\Delta^0\cC/T.$ and
$\Delta^0_n\cC/sk_nT.$ of simplicial obects over $T.$ and $n$-truncated simplicial
objects over $sk_n T.$, respectively. Then $sk_n$ induces a functor
$$
sk_n: \Delta^0\cC/T. \longrightarrow \Delta^0_n\cC/sk_nT. \,,
$$
and this functor has a right adjoint
$$
cosk^{T.}_n: \Delta^0_n\cC/sk_nT. \longrightarrow \Delta^0\cC/T.
$$
which we may call $n$-coskeleton over $T.$\,. In fact, it is easily checked that
we have
$$
cosk^{T.}_n(Z.) = cosk^S_n(Z.)\times_{cosk^S_n(sk_nT.)}T.
\leqno{(7.3)}
$$
for an $n$-truncated simplicial object $Z.$ where the morphism $T. \rightarrow cosk^S_n(sk_nT.)$
comes from adjunction.
In other words, we have $cosk^{T.}_n(Z.)_m = cosk^S_n(z.)_m\times_{cosk^S(T.)_m}T_m$
for all $m$. This does not depend on $S$ but only on $T.$, because we can always
replace $S$ by any other $S'$ over which $T.$ lies: If $T.$ is a simplicial $S$-object, then any simplicial
object over $T.$ is canonically a simplicial $S$-object, and the right hand side of (7.3) is the
same for $S'$ and $S$. If $\underline{T}.$ is the constant simplicial object associated to
an $S$-object $T$ (which is the constant functor $\Delta \rightarrow \cC/S$ with value $T$),
then a simplicial object $X$ over $\underline{T}.$ is the same as a simplicial $T$-object,
and one has $cosk^{\underline{T}.}_n(X) = cosk^T_n(X)$.  \\

Now let $\cP$ be a class of morphisms in $\cC$ which contains all isomorphisms and
is closed under composition and base change.

\vspace{0,5cm}

\noindent\textbf{Definition 7.4}
(a) A morphism $f: X \rightarrow Y$ in $\cC$ is called a $\cP$-cover, if it is in $\cP$.

\smallskip
(b) For an object $S$ in $\cC$, a simplicial $S$-object $X.$ is called a $\cP$-hypercover of $S$, if
for all $n \geq 0$ the canonical morphism
$$
can: X_n \rightarrow cosk^{S}_{n-1}(sk_{n-1}(X.))_n
$$
coming from adjunction is in $\cP$ (where $cosk^S_{-1}(X.) := \underline{T}.$).

\smallskip
(c) More generally, a morphism $f: X. \rightarrow Y.$ of simplicial objects in $\cC$ is
called a $\cP$-hypercover, if for all $n\geq 0$ the canonical adjunction morphism
$$
can: X_n \rightarrow cosk^{Y.}_{n-1}(sk_{n-1}(X.))_n
$$
is in $\cP$ (where $cosk^{Y.}_{-1}(X.) := Y.$).

\vspace{0,5cm}

Here (b) is a special case of (c), because (b) means that $X. \rightarrow \underline{S}.$
is a $\cP$-hypercover. For the following we note that, by (7.3), we have
$$
cosk^{Y.}_m(sk_m(X.) = Cosk_m(X.)\times_{Cosk_m(Y.)}Y.\,,
\leqno{(7.5)}
$$
where we define $Cosk_m(X.) = cosk_m(sk_m(X.))$.

\vspace{0,5cm}

\noindent\textbf{Lemma 7.6} (a) If $X. \rightarrow Y.$ is a $\cP$-hypercover and
$Y'. \rightarrow Y.$ is any morphism, then $X.\times_{Y.}Y'. \rightarrow Y'.$
is a $\cP$-hypercover.

\smallskip
(b) If $f: X. \rightarrow Y.$ and $g: Y. \rightarrow Z.$ are $\cP$-hypercovers, then
$gf: X. \rightarrow Z.$ is a $\cP$-hypercover.

\vspace{0,5cm}

\noi\textbf{Proof} By assumption the class $\cP$ is closed under isomorphisms, pullbacks
and composition. Therefore both claims follow from (7.5) and the fact
that $sk_n$ and $cosk_n$ and hence $Cosk_n$ commute with (fiber) products.
Indeed, for (a) one uses the formula
$$
Cosk_m(X.\times_{Y.}Y'.)\times_{Cosk_m(Y'.)}Y'.
= Cosk_m(X.)\times_{Cosk_m(Y.)}Y'.
$$
and for (b) the factorization of the morphism $can: X. \rightarrow Cosk_m(X.)\times_{Cosk_m(Z.)}Z.$
into
$$
X. \mathop{\longrightarrow}\limits^{can} Cosk_m(X.)\times_{Cosk_m(Y.)}Y.
\mathop{\longrightarrow}\limits^{id\times can} Cosk_m(X.)\times_{Cosk_m(Y.)}Cosk_m(Y.)\times_{Cosk_m(Z.)}Z.
$$
and the obvious contraction isomorphism.

\v\noi
We need the following generalization.

\v\noi
\textbf{Definition 7.7} Let $\cC^\Delta$ be the category of simplicial objects in $\cC$.
(a) A morphism $f: X. \rightarrow Y.$ in $\cC^\Delta$ is called
a simplicial $\cP$-cover, if all morphisms $f_n: X_n \rightarrow Y_n$ are $\cP$-covers.

\smallskip
(b) A simplicial $\cP$-hypercover of an object in $\cC^\Delta$ is a $s\cP$-hypercover
where $s\cP$ is the class of simplicial $\cP$-covers. (Note that this class contains
all isomorphisms and is closed under composition and base change.)

\v\noi
We note the following result. Let $\cC' \subset \cC$ be a full subcategory.

\v\noi
\textbf{Lemma 7.8} Assume that the following conditions hold for $\cC' \subset \cC$ and the class $\cP$.

\smallskip\noi
(i) Finite sums (coproducts) exist in $\cC$, and $\cC'$ is closed under finite sums.

\smallskip\noi
(ii) A morphism $f: \coprod_{i=1}^m X_i \rightarrow Y$ on a finite sum lies in $\cP$ if for one
$i$ the restriction $X_i \rightarrow Y$ of $f$ lies in $\cP$.

\smallskip\noi
(iii) Any object $Y$ in $\cC$ possesses a $\cP$-cover $f: X \rightarrow Y$ with $X$ in $\cC'$.

(a) Then any object $X$ in $\cC$ possesses a (split) $\cP$-hypercover $f: X. \rightarrow Y$ with $X.$ a simplicial object in $\cC'$, and any
simplicial object $Y.$ in $\cC$ possesses a (split) $\cP$-hypercover $f: X. \rightarrow Y.$ with $X.$ a
simplicial object in $\cC'$.

(b) Similarly, any object $Y. \in \cC^\Delta$ possesses a simplicial $\cP$-cover $X. \rightarrow Y.$ in $\cC'$,
and a simplicial $\cP$-hypercover $X.. \rightarrow Y.$ with a bisimplicial objects in $\cC'$.

\v\noi
(a) The proof follows by the usual constructions, see \cite{SGA 4.2} V bis, 5.1.3,
and \cite{GS} Lemma 2: If $sk_{n-1}(X.)$ has already been defined, one sets
$$
X_n = Z_n \; \textstyle{\coprod} \;\mathop{\textstyle{\coprod}}\limits_{i=0}^{n-1} X_i
$$
where $\alpha_n: Z_n \rightarrow (cosk_{n-1}^{Y.}sk_{n-1}X.)_n$ is a $\cP$-cover, with the
natural degeneration and face maps between $X_n$ and $X_{n-1}$, and the obvious
morphism $f_n: X_n \rightarrow Y_n$ which on the component $Z_n$ is $\alpha_n$
followed by the natural projection $(cosk_{n-1}^{Y.}sk_{n-1}X.)_n \rightarrow Y_n$.
Condition (ii) implies that $f_n$ is in $\cP$.

(b) The class $s\cP$ of morphisms in $\cC^\Delta$ satisfies conditions (i), (ii) and (iii),
the latter because any $\cP$-hypercover $Y. \rightarrow X.$ is also a simplicial $\cP$-cover
(see \cite{GS2} Cor. 2.7),
so that we can apply (a) to $\cC^\Delta$, ${\cC'}^\Delta$ and $s\cP$.

\v\noi
The following definition can be found in \cite{Gi} p. 89 and \cite{GS} 1.4.1.

\vspace{0,5cm}

\noindent\textbf{Definition 7.6}
A morphism $f: X \rightarrow Y$ of schemes is called an envelope, if it is
proper, and for any $y\in Y$ there is a point $x \in X$ such that $f(x) = y$
and the morphism $k(x) \rightarrow k(y)$ is an isomorphism.

\v
It is easy to see that these morphisms form a class $\cP$ as above, and the corresponding
hypercoverings are called hyperenvelopes. We are interested in the following variant (see also
\cite{GS2} Definition 5.14).

\v\noi
\textbf{Definition 7.7} (a) Let $\ell$ be a prime. A morphism $f: X \rightarrow Y$ of schemes is
called an $\ell'$-cover, if it is proper, and if for any point $y\in Y$ there is a point $x\in X$
such that $f(x)=y$ and $k(x)/k(y)$ is a finite field extension of degree prime to $\ell$. Call
$f$ small if $X\times_Y\eta$ is finite over $\eta$ for any generic point $\eta \in Y$.

\smallskip
(b) A morphism $f: X. \rightarrow Y.$ of simplicial schemes is called a (small) simplicial $\ell'$-cover,
if every $f_n: X_n \rightarrow Y_n$ is a (small) $\ell'$cover.

\v
It is easy to see that in each case these morphisms form classes $\cP$ as in the beginning (in case
(b) in the category of simplicial schemes).
We call the corresponding hypercoverings $\ell'$-hypercoverings in case (a) and simplicial $\ell'$-hypercoverings
in case (b).

\vspace{0.5cm}

\noi
Uwe Jannsen\\
Fakult\"at f\"ur Mathematik\\
Universit\"at Regensburg\\
93040 Regensburg\\
GERMANY\\
uwe.jannsen@mathematik.uni-regensburg.de

\end{document}